\documentclass[letterpaper, 10 pt, conference]{ieeeconf}%
\usepackage{graphicx}
\usepackage{epsfig}
\usepackage{mathptmx}
\usepackage{times}
\usepackage{amsmath}
\usepackage{amssymb}
\usepackage[section]{placeins}
\usepackage{placeins}
\usepackage{amsmath}
\usepackage{multirow}
\usepackage{graphicx}
\usepackage[normalem]{ulem}
\usepackage{hyperref}
\usepackage{subcaption}
\usepackage{algorithm,tabularx}
\usepackage{algpseudocode}
\usepackage{amsfonts}
\usepackage{cases}%
\usepackage{romannum}
\usepackage{epstopdf}
\usepackage{dirtytalk}
\providecommand{\U}[1]{\protect\rule{.1in}{.1in}}
\IEEEoverridecommandlockouts
\newcommand{\R}{\mathbb{R}}

\newcommand{\tsup}[1]{\textsuperscript{#1}}
\newcommand{\mb}[1]{\mathbf{#1}}
\newtheorem{assumption}{Assumption}
\newtheorem{theorem}{Theorem}

\newtheorem{lemma}{Lemma}
\newtheorem{remark}{Remark}
\useunder{\uline}{\ul}{}
\makeatletter
\newcommand{\multiline}[1]{  \begin{tabularx}{\dimexpr\linewidth-\ALG@thistlm}[t]{@{}X@{}}
#1
\end{tabularx}
}
\makeatother
\usepackage{xcolor}
\begin{document}

\title{{\LARGE \textbf{Distributed Non-convex Optimization of Multi-agent Systems Using Boosting Functions to Escape Local Optima: Theory and Applications}}}
\author{Shirantha Welikala and Christos G. Cassandras
\thanks{$^{\star}$Supported in
part by NSF under grants ECCS-1931600, DMS-1664644, CNS-1645681, by AFOSR
under grant FA9550-19-1-0158, by ARPA-E's NEXTCAR program under grant
DE-AR0000796 and by the MathWorks.}
\thanks{Authors are with the Division of Systems Engineering and Center for Information and Systems Engineering, Boston University, Brookline, MA 02446 \texttt{{\small \{shiran27,cgc\}@bu.edu}}.}}
\maketitle

\begin{abstract}
We address the problem of multiple local optima arising due to non-convex objective functions in cooperative multi-agent optimization problems. To escape such local optima, we propose a systematic approach based on the concept of \emph{boosting functions}. The underlying idea is to temporarily transform the gradient at a local optimum into a \emph{boosted}\ \emph{gradient} with a non-zero magnitude. We develop a Distributed Boosting Scheme (DBS) based on a gradient-based optimization algorithm using a novel optimal variable step size mechanism so as to guarantee convergence. Even though our motivation is based on the coverage control problem setting, our analysis applies to a broad class of multi-agent problems. Simulation results are provided to compare the performance of different boosting functions families and to demonstrate the effectiveness of the boosting function approach in attaining improved (still generally local) optima.

\end{abstract}

\thispagestyle{empty} \pagestyle{empty}

\section{Introduction}

A cooperative multi-agent system is a collection of interacting subsystems
(also called \textit{agents}), where each agent controls its \textit{local
state} so as to collectively optimize a common \textit{global objective}
subject to various constraints. In a \emph{distributed} optimization approach, each agent controls its state using only locally available information. The goal is to drive all agents to a globally optimal set of states. This can be a challenging task depending on the nature of: 
(\romannum{1}) the agents (which may be sensor nodes, vehicles, robots, supply sources, or processors of a multi-core computer), (\romannum{2}) the constraints on their decision space,
(\romannum{3}) the inter-agent interactions, and, (\romannum{4}) the global objective function. Therefore, a large number of optimization methods can be
found in the literature specifically developed to address different classes of
multi-agent systems.

\subsection{Literature Review}

Cooperative multi-agent system optimization arises in coverage control
\cite{Zhong2011}, formation control \cite{Lin2014}, monitoring \cite{Zhou2018}%
, flocking \cite{Ghapani2016}, resource allocation \cite{Marden2014,Sun2020}, learning
\cite{Panait2005,Xu2015}, consensus \cite{SunH2017,Su2015}, smart grid
\cite{Gonzalez-Briones2018,Molzahn2017,Pipattanasomporn2009}, transportation
\cite{Burmeister1997,Dotoli2013} and smart cities \cite{Anagnostopoulos2018}.
In these applications, gradient-based techniques are typically used due to their simplicity (see the survey paper \cite{Nedic2018}). However, more computationally complex schemes, e.g., using the Alternating Direction Method of Multipliers (ADMM) \cite{Boyd2010},\cite{Bastianello2018CDC}, are also gaining
popularity due to their greater generality.

In some multi-agent systems, properties of the associated objective function,
such as convexity, can be exploited to achieve a global optimum. For example,
the Relaxed-ADMM approach in \cite{Bastianello2018CDC} converges to the global
optimum for convex objective functions (along with a few minor additional
conditions). On the other hand, there are many settings where the objective
function takes a non-convex form making it difficult to attain a global
optimum (e.g., see \cite{Zhong2011,Sun2014,Nedic2018,Boyd2010}). In such
situations, one often resorts to global optimization techniques such as
simulated annealing \cite{Kirkpatrick1983,Chiu}, genetic algorithms
\cite{Holland1984}, or particle swarm algorithms \cite{Kennedy1995} (see the
survey papers \cite{Floudas2008,Arora1995}). The common feature in these
approaches is to introduce an element of randomness in the process of
controlling agents. Along the same lines, the ladybug exploration method
proposed in \cite{Schwager2008} tries to hover over probable local optima
solutions aiming to find a better optimum. These methods are computationally
intensive and usually infeasible for on-line optimization.

The issue of non-convexity in the objective functions has recently attracted
renewed attention for specific classes of multi-agent systems by exploiting
properties that the objective function may possess. For example, when the
objective function is submodular, tight performance bound guarantees may be
found \cite{Sun2017}. Methods like local optima smoothing \cite{Addis2005} and
balanced detection \cite{Zhong2011} trade-off local approximations and global
exploration of the objective function to achieve a better optimum. In
\cite{Sun2014}, the concept of a \textquotedblleft \textit{boosting function}\textquotedblright\ is used to escape local optima and seek better optima solutions through an exploration of the search space which exploits the objective function's structure. However, none of these methods so far is designed to function in a generic \textit{distributed} multi-agent setting and convergence
guarantees are also lacking.

\subsection{Background work}
In this paper, we propose a distributed approach to solve general non-convex
multi-agent optimization problems, inspired by the centralized boosting function
approach proposed in \cite{Sun2014}. The key idea behind boosting function approach is to temporarily alter the local
objective function of an agent whenever an equilibrium is reached, by defining
a new auxiliary local objective function. This is done indirectly by
transforming the local gradient (of the local objective) to get a new
\textit{boosted gradient} (which corresponds to the gradient of the unknown
auxiliary local objective). Therefore, a boosting function, formally, is a
transformation of the local gradient, whenever it becomes zero; the result of
the transformation is a non-zero \textit{boosted gradient}. After following
the boosted gradient, when a new equilibrium point is reached, we revert to the original objective function, and the gradient-based algorithm converges to a new (potentially better) equilibrium point. In contrast to randomly perturbing the gradient components (e.g., \cite{Kirkpatrick1983}),
boosting functions provide a systematic way to force each agent to move in a
well-chosen direction that further explores the feasible space based on
structural properties of the objective function and knowledge of both the feasible space and of the current agent states. 

Typically, when an agent follows the boosted gradient direction, it is said to be in the \textit{boosting mode}, and otherwise, it is said to be in a \textit{normal mode}. Further, the underlying technique of computing the boosted gradient is called as a \textit{boosting function family}. Furthermore, a \textit{boosting scheme} defines how each agent switch between following boosting mode and normal mode.

Details on the design of boosting functions and their use in the distributed optimization framework are discussed in this work using the class of multi-agent coverage control problems. In coverage control problems, the objective is to determine the best arrangement for a set of agents (e.g., sensor nodes) in a given mission space to maximize the probability of detecting randomly occurring events over this space. Typically the associated objective function has a non-convex form \cite{Zhong2011} resulting in multiple locally optimal configurations. Therefore, for the coverage control problems, the use of boosting functions approach as a means of escaping local optima is justified.

\subsection{Contributions}
The contribution of this paper is first to provide a formal analysis of the original centralized boosting scheme (CBS) \cite{Sun2014} (in a more generic problem setting) so as to establish convergence and then to develop a generic distributed boosting scheme (DBS) whereby each agent may
asynchronously switch between a boosting and a normal mode independent of other agents. We show that the latter scheme also converges, i.e., the asynchronous distributed boosting process reaches a terminal point where a new (generally local but improved) optimum is reached. Central to this process is a method for selecting \emph{optimal variable step sizes} in the underlying distributed gradient-based optimization algorithm. These theoretical contributions have been independently discussed in authors' paper \cite{Welikala2019P1}. Although our motivation for the contributions mentioned above comes from the coverage control problem setting, they apply to a broad class of multi-agent systems, beyond coverage or consensus-like problems. 

To provide specific details on the process of boosting functions design, we consider the class of multi-agent coverage control problems. With regard to that, the conventional coverage control framework in \cite{Zhong2011} is first enhanced to incorporate the effect of possible discontinuities in the sensing functions employed by agents. Next, two new boosting function families are introduced for the class of coverage control problems (termed \emph{Arc-Boosting} and \emph{V-Boosting}). Finally, based on the previously obtained theoretical results (on the convergence of generic DBS), a novel convergence guaranteed DBS is proposed for the coverage control problems. 


\subsection{Organization}

The general cooperative multi-agent optimization problem is introduced in Section \ref{Sec:GeneralProblemFormulation} along with the concept of boosted gradients and associated boosting schemes. Section \ref{Sec:VariableStepSizeScheme} formally develops the proposing optimal variable step size selection mechanism based on which we show the convergence of general distributed boosting schemes. Then, Section \ref{Sec:CoverageControlApplication} is dedicated to illustrating an application of developed boosting concepts for the class of multi-agent coverage control problems. 

Specifically, Section \ref{SubSec:CoverageControlProblemFormulation} revisits the multi-agent coverage control problem and Section \ref{SubSec:CoverageControlDistributedOptimization} presents its distributed gradient-ascent-based solution technique. Next, Section \ref{SubSec:CoverageControlBoostingFunctions} introduces the concept of boosting functions and boosting function families for coverage control application. The proposing DBS is discussed in Section \ref{SubSec:CoverageControlBoostingSchemes}. Section \ref{SubSec:CoverageControlDBSConvergence} describes the application of developed convergence guaranteeing optimal variable step size scheme to the coverage control problem. Finally, Section \ref{SubSec:SimulationResults}
presents simulation results illustrating the effectiveness of the introduced distributed boosting framework and Section \ref{Sec:Colclusion} concludes the paper stating some interesting future research directions.

\section{Problem formulation}\label{Sec:GeneralProblemFormulation}

We consider cooperative multi-agent optimization problems of the general
form,
\begin{equation}
\mathbf{s}^{\ast}=\arg\max_{\mathbf{s}\in\mathbf{F}}H(\mathbf{s}),
\label{Eq:GeneralizedOptimizationProb}%
\end{equation}
where, $H:\mathbb{R}^{mN}\rightarrow\mathbb{R}$ is the \textit{global
objective} function and $\mathbf{s}=[s_{1},s_{2},\ldots,s_{N}]\in
\mathbb{R}^{mN}$ is the controllable \textit{global state}. Here, for any
$i\in\{1,2,\ldots,N\}$,\ $s_{i}\in\mathbb{R}^{m}$ represents the \textit{local
state} of agent $i$. Further, $\mathbf{F}\subseteq\mathbb{R}^{mN}$ represents
the feasible space for $\mathbf{s}$. In this work, linearity or
convexity-related conditions are not imposed on the global
objective function $H(\mathbf{s})$.

In order to model the inter-agent interactions, an undirected graph denoted by
$\mathcal{G}=(\mathcal{V},\mathcal{A})$ is used where $\mathcal{V}%
=\{1,2,\ldots,N\}$ is a set of $N$ agents, and, $\mathcal{A}$ is the set of
communication links between those agents. The set of \textit{neighbors} of an
agent $i\in\mathcal{V}$ is denoted by $B_{i}=\{j:j\in\mathcal{V}%
,\ (i,j)\in\mathcal{A}\}$. The \textit{closed neighborhood} of agent $i$ is
defined as $\bar{B}_{i}=B_{i}\cup\{i\}$ and $\vert \bar{B}_{i}\vert$ denotes the cardinality of the set $\bar{B}_{i}$. It is assumed that each agent $i$
shares its local state information $s_{i}$ with its neighbors in $B_{i}$. As a
result, agent $i$ has knowledge of its \textit{neighborhood state}
$\bar{s}_{i}$, where $\bar{s}_{i}=\{s_{j}:j\in\bar{B}_{i}\}$.

In this problem setting, an agent $i$ is also assumed to have a \textit{local
objective} function $H_{i}(\bar{s}_{i})$ where $H_{i}:\mathbb{R}^{m|\bar
{B}_{i}|}\rightarrow\mathbb{R}$. Note that $H_{i}(\bar{s}_{i})$ only depends
on agent $i$'s neighborhood state $\bar{s}_{i}$. 
The relationship between local and global objective functions is not restricted to any specific form except for the condition: 
\begin{equation}\label{Eq:TechnicalCondition}
\frac{\partial H_i(\bar{s}_i)}{\partial s_i}=0,\ \forall i \in\mathcal{V} \implies \nabla H(\mathbf{s})=0.
\end{equation}
This condition clearly holds for any problem with a \emph{separable} form \cite{Sun2014} 
$H(\mathbf{s})=H_{i}(\bar{s}_{i})+H_{i}^{c}(s_{i}^{c})$ where
$H_{i}^{c}:\mathbb{R}^{m(N-1)}\rightarrow\mathbb{R}$ and 
$s_{i}^{c}=[s_{1},s_{2},\ldots,s_{i-1},s_{i+1},\ldots,s_{N}]$.
Note that cooperative multi-agent systems which are inherently distributed (e.g., \cite{Zhong2011}) naturally have separable objective functions. Moreover, many problems of interest with an \emph{additive} form \cite{Bastianello2018CDC} $H(\mathbf{s})=\sum_{i=1}^{N}H_{i}(\bar{s}_{i})$ also satisfy this condition


\subsection{Distributed gradient-ascent method}

Due to the versatile nature of $H$ and $\mathbf{F}$ in
\eqref{Eq:GeneralizedOptimizationProb}, applicable solving techniques are
limited to global optimization methods. Even-though many such techniques are
available \cite{Floudas2008,Arora1995}, in this paper we consider a simple
gradient-ascent scheme so as to take advantage of its simplicity in terms of
analysis, computation, and on-line implementation despite the obvious
limitation of attaining only local optima. We are also interested in solving
\eqref{Eq:GeneralizedOptimizationProb} through \emph{distributed} schemes so
that each agent $i$ updates its local state $s_{i}$ according to
\begin{equation}
s_{i,k+1}=s_{i,k}+\beta_{i,k}d_{i,k},\label{Eq:GradientAscentGeneral}%
\end{equation}
where, $\beta_{i,k}\in\mathbb{R}$ is a step size, and $d_{i,k}=\frac{\partial
H_{i}(\bar{s}_{i,k})}{\partial s_{i}}\in\mathbb{R}^{m}$ denotes the locally
available gradient.

\subsection{Escaping local optima using boosting functions}

\label{SubSec:IntroductionToBoosting}

Converging to a local optimum is the main drawback of using a gradient-based
method like \eqref{Eq:GradientAscentGeneral}, when the global objective
function $H$ is non-convex and/or the feasible space $F$ is non-convex. In
\cite{Sun2014}, this problem has been addressed by introducing the concept of
\textit{boosting functions} as an effective systematic method of escaping
local optima.

\subsubsection{Boosting functions}

The main idea here is to temporarily alter the local objective function $H_i(\bar{s}_i)$ whenever an equilibrium is reached with a newly defined auxiliary objective function $\hat{H}_{i}(\bar{s}_i)$. However, we are interested in the boosted gradient $\hat{d}_{i,k} = \frac{\partial \hat{H}_i(\bar{s}_{i,k})}{\partial s_i}$ rather than $\hat{H}_i(\bar{s}_i)$. A boosted gradient is a transformation of the associated local gradient $d_{i}$ taking place at an equilibrium point (where its value is zero); the result of the transformation is a non-zero $\hat{d}_{i}\neq0$ which, therefore, forces the agent to move in a direction determined by the boosting function and to explore the feasible space further. When a new equilibrium point is reached, we revert to the original objective function, and then the gradient-based algorithm converges to a new (potentially better and never worse) equilibrium point. 

The key to
boosting functions is that they are selected to exploit the structure of the
objective functions $H(\mathbf{s})$ and $H_{i}(\bar{s}_{i})$, of the feasible space $\mathbf{F}$, and of the agent state trajectories. Unlike various forms of randomized state perturbations away from their current equilibrium \cite{Kirkpatrick1983,Chiu}, boosting functions provide a formal rational systematic transformation process of the form $\hat{d}_{i}=f(d_{i},\bar{s}%
_{i})$ where the boosting function $f$ depends on the specific problem type. Details on boosting functions and their design process are discussed in
Section \ref{SubSubSec:BoostingFunctionConstruction} and \ref{Sec:CoverageControlApplication}. In what follows, we present a general-purpose boosting function choice
to provide insight into boosting functions in a generic setting.

In many multi-agent optimization problems, local optima arise
when a cluster of agents provides a reasonably high performance by maintaining
their local states in close proximity while completely ignoring globally
dispersed state configurations. In such a case, a boosting function that enhances a separation among local states is a natural choice, especially suited for applications like coverage control, formation control, monitoring, consensus and transportation. In fact, for coverage control problems, such a boosting function has already been proven to be effective in \cite{Sun2014}. Therefore, in a generic setting, a candidate boosted gradient $\hat{d}_{i}=f_{i}(d_{i})$ for agent $i$ can be obtained by letting $\psi_{ij}=(s_{i}-s_{j})$ and defining $\hat{d}_{i}=\nabla_{\psi_{ij}}H_{i}(\bar{s}_{i})$ where its $l$\textsuperscript{th}  component is
\begin{equation}
\hat{d}_{i}^{l}=\frac{\partial H_{i}(\bar{s}_{i})}{\partial\psi_{ij}^{l}%
}=\underbrace{\frac{\partial H_{i}(\bar{s}_{i})}{\partial s_{i}^{l}}%
}_{=\ d_{i}^{l}}\frac{\partial s_{i}^{l}}{\partial\psi_{ij}^{l}}%
+\underbrace{\frac{\partial H_{i}(\bar{s}_{i})}{\partial s_{j}^{l}}%
}_{\triangleq\ d_{ji}^{l}}\frac{\partial s_{j}^{l}}{\partial\psi_{ij}^{l}}.
\end{equation}
Now, by replacing $\frac{\partial s_{i}^{l}}{\partial\psi_{ij}^{l}}$ and
$\frac{\partial s_{j}^{l}}{\partial\psi_{ij}^{l}}$ with scalar parameters
$\alpha_{ij}$ and $\eta_{ij}$, an entire \emph{family of boosting functions}
can be obtained as $\hat{d}_{i}=f_{i}(d_{i})=\alpha_{ij}d_{i}+\eta_{ij}d_{ji}$
where $d_{ji}=\frac{\partial H_{i}(\bar{s}_{i})}{\partial s_{j}}$ (see also
\eqref{Eq:WeightTransforms} and \eqref{Eq:BoostedGradientGeneral}). Note that
setting $\alpha_{ij}=1$ and $\eta_{ij}=-1$ gives an interesting boosting
function choice of the form $\hat{d}_{i}=f_{i}(d_{i})=d_{i}-d_{ji}$. Since
$d_{ji}$ represents the direction towards which agent $j$ should move to
increase $H_{i}$, this is clearly an intuitive general choice for a boosting
function at agent $i$. 


\subsubsection{A boosting scheme} 

\begin{figure}[!t]
\centering
\includegraphics[width=3.4in]{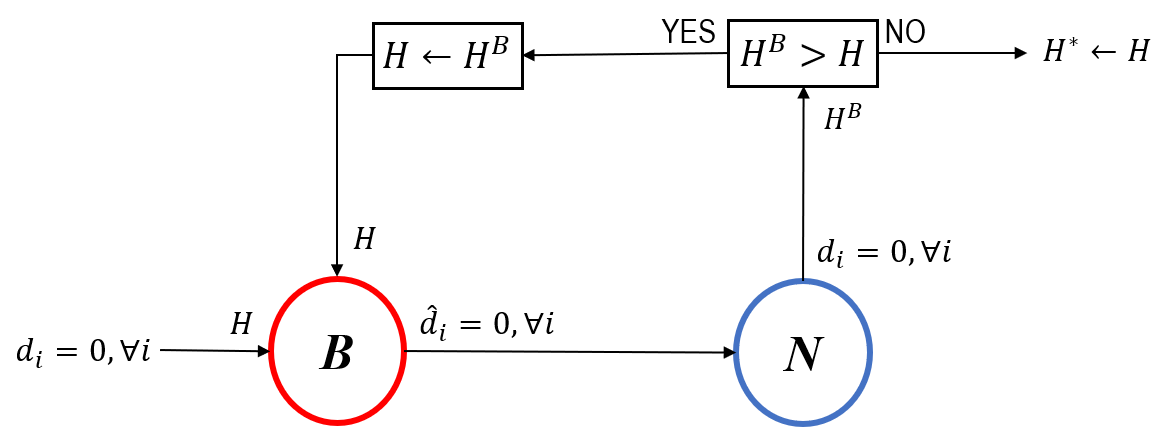} \caption{A
centralized boosting scheme (CBS) (see also Section \ref{SubSubSec:TerminationConditionsForImplementations})}%
\label{Fig:CentralizedBoostingBlockDiagram}%
\end{figure}

When an agent $i$ is following the boosted gradient
direction $\hat{d}_{i,k}$, it is said to be in the \textit{Boosting Mode} and
its state updates take the form
\begin{equation}
s_{i,k+1}=s_{i,k}+\beta_{i,k}\hat{d}_{i,k}.\label{Eq:BoostedGradientAscent}%
\end{equation}
Similarly, when an agent $i$ is following the \textquotedblleft
normal\textquotedblright\ gradient direction $d_{i,k}$ as in
\eqref{Eq:GradientAscentGeneral}, it is said to be in the \textit{Normal
Mode}. When developing an optimization scheme to solve
(\ref{Eq:GeneralizedOptimizationProb}), we need a proper mechanism, referred
to as a \textit{Boosting Scheme, }to switch the agents between normal and
boosting modes. A centralized boosting scheme (CBS) is outlined in Fig.
\ref{Fig:CentralizedBoostingBlockDiagram}, where the normal mode is denoted by
\textbf{\textit{N}} and the boosting mode is denoted by \textbf{\textit{B}}.
In a CBS, all agents are synchronized to operate in the same mode. In Fig.
\ref{Fig:CentralizedBoostingBlockDiagram}, $H$ denotes the global objective
function value which is initially stored by all agents the first time mode
\textbf{\textit{B}} is entered when $d_{i}=0$ for all $i\in\mathcal{V}$. After
$\hat{d}_{i}=0$ for all $i\in\mathcal{V}$, the agents re-enter mode
\textbf{\textit{N}} and, when a new equilibrium is reached, the new
post-boosting value of the global objective function $H(\mb{s})$ 
is denoted by $H^{B}$. If $H^{B}>H$, an
improved equilibrium point is attained and the process repeats by re-entering
mode \textbf{\textit{B}} with the new value $H^{B}$. The process is complete
when this centralized controller fails to improve $H(\mb{s})$, i.e., when $H^{B} \leq H$.

This CBS was used in \cite{Sun2014} with appropriately defined boosting functions in mode \textit{\textbf{B}} to obtain improved performance for a variety of multi-agent coverage control problems (A more detailed version this CBS is discussed in Section \ref{Sec:CoverageControlApplication}). However, there
has been no formal proof to date that this process converges.
Moreover, our goal is to develop a \emph{Distributed Boosting Scheme} (DBS)
where each agent can independently switch between modes \textbf{\textit{B}}
and \textbf{\textit{N }}at any time. Such a scheme 
(\romannum{1}) improves the scalability of the system, 
(\romannum{2}) eliminates the requirement of a centralized controller, 
(\romannum{3}) reduces computational and communication costs, and, 
(\romannum{4}) can potentially improve convergence times. Furthermore, in problems such as coverage control \cite{Zhong2011}, the original problem is inherently
distributed and makes a DBS a natural approach.

\begin{figure}[!t]
\centering
\includegraphics[width=3.4in]{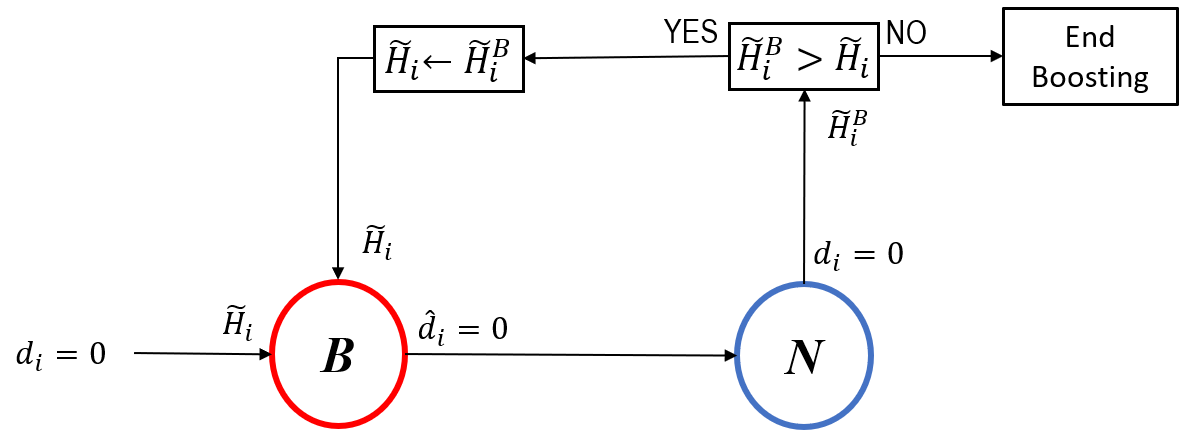}
\caption{A distributed boosting scheme (DBS) asynchronously
applied by each agent $i=1,\ldots,N$ (see also Section \ref{SubSubSec:TerminationConditionsForImplementations})}
\label{Fig:DecentralizedBoostingBlockDiagram} 
\end{figure}

A simple DBS version of Fig. \ref{Fig:CentralizedBoostingBlockDiagram} is shown in Fig. \ref{Fig:DecentralizedBoostingBlockDiagram} where local use of the global objective $H$ is now replaced by a local estimate of $H$, denoted by $\bar{H}_i$, which will be formally introduced later (in \eqref{Eq:NeighborhoodCostFunctionTheory}). An application-specific variation of this DBS is discussed in Section \ref{Sec:CoverageControlApplication}, but they do not affect the convergence analysis that follows.  

One can see that convergence of the DBS is far from obvious since agents may be at
different modes at any time instant and, as their states change, the
interaction among agents could lead to oscillatory behavior. Note that the
notion of convergence involves not only the existence of equilibria such that
$d_{i}=0$ or $\hat{d}_{i}=0$, but also a guarantee that the condition
$H^{B}\leq H$ is eventually satisfied. We will show that a key
to guaranteeing convergence is a process for \textit{optimally selecting a variable step size} $\beta_{i,k}$ in (\ref{Eq:GradientAscentGeneral})
and (\ref{Eq:BoostedGradientAscent}).

\subsubsection{Guidelines on selecting boosting functions}\label{SubSubSec:BoostingFunctionConstruction}
As mentioned before, the main objective of a boosting function $f$ (recall $\hat{d}_i = f(d_i,\bar{s}_i)$) is to provide a meaningful boosted gradient direction $\hat{d}_i$ for any agent $i$ to follow upon reaching a local optimum (where $d_i = 0$). The meaningfulness of boosted gradients comes from the fact that underlying boosting function $f$ is constructed considering one or few of the \textit{source factors}: (\romannum{1}) Nature of the objective functions $H, H_i$, (\romannum{2}) Nature of the gradient expression $d_i$, (\romannum{3}) Nature of the feasible space $\mb{F}$, and, (\romannum{4}) State trajectory of the agent $i$. As a result of this, boosted gradients are expected to achieve the capabilities to: (\romannum{1}) Drive agents beyond local optima, (\romannum{2}) Explore potentially good search directions, and, (\romannum{3}) Explore undiscovered regions in the feasible space. Thus, compared to the methods where agents are randomly perturbed to escape from local optima, boosting function approach should be more effective.

Although developing a closed-form equation for the boosted gradient $\hat{d}_{i,k}$ (i.e., for $f$) heavily depends on the application, a few basic guidelines can be proposed as follows.
\begin{enumerate}
    \item In the analytical expression of $d_{i,k}$, identify the components which control its direction and magnitude. Also, identify the dependence of those components on the aforementioned source factors.
    \item A candidate expression for $\hat{d}_{i,k}$ can be obtained by transforming (altering) a sub-component of the expression of $d_{i,k}$, such that the resulting $\hat{d}_{i,k}$ satisfies the basic boosting function requirement $\hat{d}_{i,k} \neq 0$ when $d_{i,k} = 0$.
    \item Also, for an agent with $d_{i,k}=0$, try to formulate a rationale (or a list of rationals) which gives a globally beneficial new travel direction (i.e. a $\hat{d}_{i,k}$). An example rational used in the coverage control application is: `Move away from the closest neighbor', which, when followed, motivates the agents to spread out from each other - leading to a better solution with higher coverage. 
    \item Now, using the knowledge acquired in step (1), to achieve one or few of the rationals identified in step (3), appropriately design the transformation mentioned in step (2).  
\end{enumerate}
Section \ref{Sec:CoverageControlApplication} provides a complete application example of the aforementioned steps with more details.

\subsubsection{Convergence criteria}

When a DBS is considered, unlike the case of a
CBS in Fig. \ref{Fig:CentralizedBoostingBlockDiagram}, the decentralized
nature of agent behavior causes agents to switch between modes
(normal/boosting) independently and asynchronously from each other. As a
result, at a given time instant, a subset of the agents will be in normal mode
(following \eqref{Eq:GradientAscentGeneral}) while others are in boosting mode
(following \eqref{Eq:BoostedGradientAscent}). This partition of the complete
agent set $\mathcal{V}$ leads to two agent sets henceforth denoted by
$\mathcal{N}$ and $\mathcal{B}$ respectively. 

Let us define the \emph{extended neighborhood} of an agent $i$ as $\tilde{B}_i \triangleq \cup_{j\in\bar{B}_i} \bar{B}_j$. For any agent $i\in\mathcal{V}$, the following conditions are defined as the \textit{convergence criteria}:
\begin{align}
\lim_{k\rightarrow\infty}d_{i,k} &= 0\mbox{ when }\tilde{B}_{i}\subseteq
\mathcal{N},\label{Eq:RequiredConvergence1}\\
\lim_{k\rightarrow\infty}d_{i,k} &  =0\mbox{ when }i\in\mathcal{N}, \tilde{B}_{i}%
\cap\mathcal{B}\neq\emptyset,\label{Eq:RequiredConvergence2}\\
\lim_{k\rightarrow\infty}\hat{d}_{i,k} &  =0\mbox{ when }i\in\mathcal{B}%
,\tilde{B}_{i}\cap\mathcal{B}\neq\emptyset.\label{Eq:RequiredConvergence3}%
\end{align}
These convergence criteria enforce the capability of an agent $i$ to escape its current mode (normal or boosting) irrespective of the surrounding neighbor mode partitions $\tilde{B}_i\cap\mathcal{N}$ and $\tilde{B}_i\cap\mathcal{B}$. Since boosting will only continue as long as there is a gain from the boosting stages (i.e., \textquotedblleft$\bar{H}_{i}^{B}>\bar{H}_{i}$\textquotedblright\ in Fig. \ref{Fig:DecentralizedBoostingBlockDiagram}), it is clear how these criteria can lead all agents to terminate their boosting stages (i.e., to reach the \textquotedblleft End Boosting\textquotedblright\ state).

Upon this termination, criterion \eqref{Eq:RequiredConvergence1} reapplies and guarantees achieving $d_{i,k}\rightarrow0$ for all $i\in\mathcal{V}$, which will directly imply $\nabla H(\mathbf{s}_{k}) \rightarrow 0$ (from \eqref{Eq:TechnicalCondition}). Therefore, convergence to a stationary solution of \eqref{Eq:GeneralizedOptimizationProb} is achieved (again, not necessarily a global optimum). Finally, note that the criterion \eqref{Eq:RequiredConvergence1} applies to the convergence of any gradient-based method where boosting is not used.




\section{Convergence analysis through optimal variable step sizes}
\label{Sec:VariableStepSizeScheme}

This section proposes a variable step size scheme which guarantees the convergence criteria \eqref{Eq:RequiredConvergence1}-\eqref{Eq:RequiredConvergence3} required when a general problem of the form \eqref{Eq:GeneralizedOptimizationProb} is solved using \eqref{Eq:GradientAscentGeneral} and \eqref{Eq:BoostedGradientAscent}. Our main results are Theorem \ref{Th:Convergence} (which guarantees \eqref{Eq:RequiredConvergence1}) and Theorem \ref{Th:ConvergenceHat} (which guarantees \eqref{Eq:RequiredConvergence2} and \eqref{Eq:RequiredConvergence3}). These depend on some assumptions which are presented first, starting with the nature of the local objective functions.


\begin{assumption}
\label{As:PropertiesofH_i} Any local objective function $H_{i}(\bar{s}%
_{i}),\ i\in\mathcal{V}$, satisfies the following conditions:

\begin{enumerate}
\item $H_{i}(\cdot)$ is continuously differentiable and its gradient $\nabla
H_{i}(\cdot)$ is Lipschitz continuous (i.e., $\exists K_{1i}$ such that
$\forall x,y\in\mathbb{R}^{m|\bar{B}_{i}|}$, $\Vert\nabla H_{i}(x)-\nabla
H_{i}(y)\Vert\leq K_{1i}\Vert x-y\Vert$).

\item $H_{i}(\cdot)$ is a non-negative function with a finite upper bound
$H_{UB}$, i.e., $H_{i}(x)<H_{UB}<\infty,$ $x\in\mathbb{R}^{m|\bar{B}_{i}|}$.
\end{enumerate}
\end{assumption}

Through the global and local objective function relationship, this assumption
will propagate to the global objective function. However, for this work,
Assumption \ref{As:PropertiesofH_i} is sufficient.

We begin by developing an optimal variable step size scheme for agents
$i\in\mathcal{V}$ such that $\bar{B}_{i}\subseteq\mathcal{N}$ (i.e., all
neighboring agents are also in normal mode - following
\eqref{Eq:GradientAscentGeneral}). The respective convergence criterion for
this case is \eqref{Eq:RequiredConvergence1}.

\subsection{Convergence for agents $i\in\mathcal{V}$ such that $\tilde{B}%
_{i}\subseteq\mathcal{N}$}
\label{subsec:ConvergenceWhenAllNormal} 
We begin by developing an optimal variable step size scheme for agents $i\in\mathcal{V}$ such that $\tilde{B}_{i}\subseteq\mathcal{N}$, i.e., all agents in the extended neighborhood are in normal mode - following \eqref{Eq:GradientAscentGeneral}. The respective convergence criterion for
this case is \eqref{Eq:RequiredConvergence1}. For notational convenience, let
$\mathbf{q}_{i}=\{1,2,\ldots,q_{i}\}$ with $q_{i}=|\bar{B}_{i}|$ represent an
ordered (re-indexed) version of the closed neighborhood set $\bar{B}_{i}$. For
this situation, agent $i$'s neighborhood state update equation can be
expressed as $\bar{s}_{i,k+1}=\bar{s}_{i,k}+\mathbf{\bar{\beta}}%
_{i,k}\mathbf{\bar{d}}_{i,k}$ by combining \eqref{Eq:GradientAscentGeneral}
for all $j\in\bar{B}_{i}$. Here, $\bar{s}_{i,k+1},$ $\bar{s}_{i,k}$ and
$\mathbf{\bar{d}}_{i,k}$ are $mq_{i}$-dimensional column vectors;
equivalently, they may be thought of as $q_{i}\times1$ block-column matrices
with their $j$\textsuperscript{th} block (of size $\mathbb{R}^{m\times1}$, and
$j\in\mathbf{q}_{i}$) being, $s_{j,k+1},$ $s_{j,k}$ and $d_{j,k}$
respectively. Accordingly, $\mathbf{\bar{\beta}}_{i,k}$ is a $q_{i}\times
q_{i}$ block-diagonal matrix, where its $j$\textsuperscript{th} block on the
diagonal (of size $m\times m$ and $j\in\mathbf{q}_{i}$) is $\beta_{j,k}I_{m}$;
$I_{m}$ is the $m\times m$ identity matrix and $\beta_{j,k}\in\mathbb{R}$ is
the (scalar) step size of agent $j$.

Following lemma provides a modified version of the widely used descent lemma \cite{Bertsekas2016nonlinear} so that it can be used to analyze maximization problems such as \eqref{Eq:GeneralizedOptimizationProb}.
\begin{lemma}\label{Lm:ModifiedDescentLemma}
For a function $f:\mathbb{R}^{n}\rightarrow\mathbb{R}$, if the Lipschitz
continuity constant of $\nabla f$ is $L$, then, $\forall x,y\in\mathbb{R}^{n}$,
\begin{equation}
f(x+y) \geq f(x) + y^{T}\nabla f(x) - \frac{L}{2} \Vert y \Vert^{2}. 
\end{equation}
\end{lemma}
\emph{Proof: } Consider a function $g=-f:\mathbb{R}^{n}\rightarrow\mathbb{R}$. Then, the Lipschitz
continuity constant of $\nabla g$ will also be $L$. Now, we can apply the usual descent lemma \cite{Bertsekas2016nonlinear} to the function $g$ (to compare $g(x+y)$ and $g(x)$). Then, $\forall x,y\in\mathbb{R}^{n}$,  
\begin{align*}
g(x+y) &\leq g(x) + y^{T}\nabla g(x) + \frac{L}{2} \Vert y\Vert^{2}, \mbox{ multiplying by $(-1)$,}\\
-g(x+y) &\geq -g(x) - y^{T}\nabla g(x) - \frac{L}{2} \Vert y\Vert^{2}, \mbox{ using $-g = f,$}\\
f(x+y) &\geq f(x) + y^{T}\nabla f(x) - \frac{L}{2} \Vert y\Vert^{2}.
\end{align*}
\hfill $\blacksquare$

Now, under Assumption \ref{As:PropertiesofH_i}, the Lemma \ref{Lm:ModifiedDescentLemma} can be applied to a local objective function $H_{i}(\bar{s}_{i,k})$ for the
aforementioned neighborhood state update $\bar{s}_{i,k+1}=\bar{s}%
_{i,k}+\mathbf{\bar{\beta}}_{i,k}\mathbf{\bar{d}}_{i,k}$ as follows:
\begin{align}
H_{i}(\bar{s}_{i,k+1}) &  \geq H_{i}(\bar{s}_{i,k})+(\mathbf{\bar{\beta}%
}_{i,k}\mathbf{\bar{d}}_{i,k})^{T}\nabla H_{i}(\bar{s}_{i,k})-\frac{K_{1i}}%
{2}\Vert\mathbf{\bar{\beta}}_{i,k}\mathbf{\bar{d}}_{i,k}\Vert^{2}%
\nonumber\label{Eq:DescentLemmaForH_i}\\
&  =H_{i}(\bar{s}_{i,k})+\sum_{j\in\bar{B}_{i}}\left[  \beta_{j,k}d_{j,k}%
^{T}d_{ji,k}-\frac{K_{1i}}{2}\beta_{j,k}^{2}\Vert d_{j,k}\Vert^{2}\right]
\nonumber\\
&  =H_{i}(\bar{s}_{i,k})+\sum_{j\in\bar{B}_{i}}\Delta_{ji,k},%
\end{align}
with
\begin{eqnarray}
\Delta_{ji,k} &\triangleq& \beta_{j,k}d_{j,k}^{T}d_{ji,k}-\frac{K_{1i}}{2}\beta_{j,k}%
^{2}\Vert d_{j,k}\Vert^{2}\in\mathbb{R},\label{Eq:Contributionji}\\
d_{ji,k} &\triangleq& \nabla_{j}H_{i}(\bar{s}_{i,k})=\frac{\partial H_{i}(\bar{s}_{i,k}%
)}{\partial s_{j}}\in\mathbb{R}^{m}.\label{Eq:CrossDerivativedji}%
\end{eqnarray}
The term $d_{ji,k}$ in \eqref{Eq:CrossDerivativedji} gives the sensitivity of
agent $i$'s local objective $H_{i}$ to the local state $s_{j}$ of agent
$j\in\bar{B}_{i}$. Also, $K_{1i}$ is the Lipschitz constant corresponding to
$\nabla H_{i}$. Note that the term $\Delta_{ji,k}$ in
\eqref{Eq:Contributionji} depends on the step size $\beta_{j,k}$ which is
selected by agent $j\in\bar{B}_{i}$.

In \eqref{Eq:DescentLemmaForH_i}, each $\Delta_{ji,k}$ term can be thought of
as a contribution coming from neighboring agent $j$ to agent $i$, so as to
improve (increase) $H_{i}$. However, in order for an agent $i$ to know its
contribution to agent $j\in\bar{B}_{i}$ (i.e., $\Delta_{ij,k}$) the following
assumption is required.

\begin{assumption}
\label{As:LocalAvailabilityOfd_ij} Any agent $i\in\mathcal{V}$ has knowledge
of the cross-gradient terms $d_{ij,k}$, the local Lipschitz constants $K_{1j}%
$, and the objective function values $H_{j}(\bar{s}_{j,k})$ at the $k$th
update instant.
\end{assumption}

This assumption is consistent with our concept of neighborhood, where
neighbors share information through communication links. Thus, any agent $i$
has access to the parameters it requires: $d_{ij,k}(=\partial H_{j}/\partial
s_{i})$ and $K_{1j}$ from all its neighbors $j\in\bar{B}_{i}$. Note that when
the form of the local objective functions $H_{i}$ is identical and all pairs
$(H_{i},H_{j})$, $j\in B_{i}$, have a symmetric structure, Assumption
\ref{As:LocalAvailabilityOfd_ij} holds without any need for additional
communication exchanges. Many cooperative multi-agent
optimization problems have this structure, including the class of multi-agent
coverage control problems (see Lemma
\ref{Lm:LocalAvailabilityd_ij} in Section \ref{Sec:CoverageControlApplication}).

We now define a \textit{neighborhood objective} function $\tilde{H}_{i}%
(\tilde{s}_{i,k})$ for any $i\in\mathcal{V}$, where $\tilde{H}_{i}%
:\mathbb{R}^{m|\tilde{B}_{i}|}\rightarrow\mathbb{R}$, $\tilde{B}_{i}%
=\cup_{j\in\bar{B}_{i}}B_{j}$, and, $\tilde{s}_{i,k}=\{s_{j}:j\in\tilde{B}%
_{i}\}$, as follows:
\begin{equation}
\tilde{H}_{i}(\tilde{s}_{i,k})=\sum_{j\in\bar{B_{i}}}H_{j}(\bar{s}_{j,k}).
\label{Eq:NeighborhoodCostFunctionTheory}%
\end{equation}
This neighborhood objective function can be viewed as agent $i$'s estimate of
the total contribution of agents in $\bar{B}_{i}$ towards the global objective function.

\begin{remark}
\label{Rm:WightedNeighborhoodObjective} In some problems, if the global and
local objective functions are not directly related in an additive manner, then
$\tilde{H}_{i}(\tilde{s}_{i,k})=\sum_{j\in\bar{B_{i}}}w_{ij}H_{j}(\bar
{s}_{j,k})$ can be used as a candidate for the neighborhood objective
function. Here, $\{w_{ij}\in\mathbb{R}^{\geq0}:j\in\bar{B}_{i}\}$ represents a
set of weights (scaling factors). All results presented in this section can be
generalized to such neighborhood objective functions as well.
\end{remark}

\begin{remark}
\label{Rm:NeighborhoodObjectiveForDistributedBoosting} The
neighborhood objective functions play an important role in DBS because a distributed scheme comes at the cost of each agent losing the global
information $H(\mathbf{s})$. In contrast, in the CBS of Fig.
\ref{Fig:CentralizedBoostingBlockDiagram}, $H(\mathbf{s})$ plays a crucial
role in the \textquotedblleft$H^{B}>H$\textquotedblright\ block. As a remedy,
in a DBS each agent $i$ uses a neighborhood objective function $\tilde{H}_{i}$
as a means of locally estimating the global objective function value (see \textquotedblleft$\tilde{H}_{i}^B>\tilde{H}_{i}$\textquotedblright\ block in Fig. \ref{Fig:DecentralizedBoostingBlockDiagram}). However,
as seen in the ensuing analysis, the form of $\tilde{H}_{i}$ is not limited to
\eqref{Eq:NeighborhoodCostFunctionTheory} - it can take any appropriate form
(see Remark
\ref{Rm:WightedNeighborhoodObjective} above and Remark
\ref{Rm:NeighborhoodObjectiveForCoverageControlDBS} in Section
\ref{Sec:CoverageControlApplication} for details).
\end{remark}


Enabled by the fact that $\tilde{B}_i \subseteq \mathcal{N} = \emptyset$, applying \eqref{Eq:DescentLemmaForH_i} to any agent $j\in \bar{B}_i$ gives $H_{j}(\bar{s}_{j,k+1})\geq H_{j}(\bar{s}_{j,k})+\sum_{l\in\bar{B}_{j}}\Delta_{lj,k}$.
Summing both sides of this relationship over all $j\in\bar{B}_{i}$ and using the definition in
\eqref{Eq:NeighborhoodCostFunctionTheory} yields
\begin{equation}
\tilde{H}_{i}(\tilde{s}_{i,k+1})\geq\tilde{H}_{i}(\tilde{s}_{i,k}%
)+(\tilde{\Delta}_{i,k}+Q_{i,k}),\label{Eq:DescentLemmaForH_Bi}%
\end{equation}
where we define
\begin{eqnarray}
\tilde{\Delta}_{i,k} &\triangleq& \sum_{j\in\bar{B}_{i}}\Delta_{ij,k}%
,\label{Eq:DeltaTildai}\\
Q_{i,k} &\triangleq& \sum_{j\in B_{i}}(\Delta_{jj,k}+\Delta_{ji,k}+\sum_{l\in B_{j}%
-\{i\}}\Delta_{lj,k}).\label{Eq:Q}%
\end{eqnarray}
Note that $\tilde{\Delta}_{i,k}$ in \eqref{Eq:DeltaTildai} is a function of
terms $\Delta_{ij,k}$ (and not $\Delta_{ji,k}$)  which are locally available to and controlled by agent $i$, i.e., $\beta_{i,k},d_{i,k}$ and $d_{ij,k},\forall
j\in\Bar{B}_{i}$. In contrast, agent $i$ does not have any control over
$Q_{i,k}$ in \eqref{Eq:Q}, as this strictly depends through
\eqref{Eq:Contributionji} on the step sizes of agent $i$'s \textit{extended
neighborhood}, i.e., $\beta_{j,k},\forall j\in\tilde{B}_{i}-\{i\}$.


Nonetheless, \eqref{Eq:DescentLemmaForH_Bi} implies that the neighborhood
objective function $\tilde{H}_{i}(\tilde{s}_{i,k})$ can be increased by at
least $(\tilde{\Delta}_{i,k}+Q_{i,k})$ at any update instant $k$. Thus, to
maximize the gain in $\tilde{H}_{i}(\tilde{s}_{i,k})$, agent $i$'s step size
$\beta_{i,k}$ is selected according to the following \textit{auxiliary
problem}:%

\begin{equation}
\label{Eq:AuxiliaryOptimizationProblem1}\begin{aligned} \beta_{i,k}^{\ast} =\ & \underset{\beta_{i,k}}{\arg\max} & & \tilde{\Delta}_{i,k} \\ & \text{subject to} & &\tilde{\Delta}_{i,k} > 0. \end{aligned}
\end{equation}

\begin{lemma}
\label{Lm:AuxProb1} The solution to the auxiliary problem
\eqref{Eq:AuxiliaryOptimizationProblem1} is
\begin{equation}
\beta_{i,k}^{\ast}=\frac{1}{\sum_{j\in\bar{B}_{i}}K_{1j}}\frac{d_{i,k}^{T}%
\sum_{j\in\bar{B}_{i}}d_{ij,k}}{\Vert d_{i,k}\Vert^{2}}.
\label{Eq:OptimalStepSize}%
\end{equation}

\end{lemma}

\emph{Proof: } Using \eqref{Eq:Contributionji} and \eqref{Eq:DeltaTildai},
$\tilde{\Delta}_{i,k}$ can be written as
\[
\tilde{\Delta}_{i,k}=\beta_{i,k}d_{i,k}^{T}\sum_{j\in\bar{B}_{i}}%
d_{ij,k}-\beta_{i,k}^{2}\Vert d_{i,k}\Vert^{2}\frac{\sum_{j\in\bar{B}_{i}%
}K_{1j}}{2}.%
\]
Note the quadratic and concave nature of $\tilde{\Delta}_{i,k}$ with respect
to agent $i$'s step size $\beta_{i,k}$. Therefore, using the KKT conditions \cite{Bertsekas2016nonlinear},
the optimal $\beta_{i,k}$ to the problem
\eqref{Eq:AuxiliaryOptimizationProblem1} can be directly obtained as
\eqref{Eq:OptimalStepSize}. Let us denote the optimal objective function value of the problem
\eqref{Eq:AuxiliaryOptimizationProblem1} as $\tilde{\Delta}_{i,k}^{\ast}$. It
is easy to show that $\beta_{i,k}^{\ast}$ in \eqref{Eq:OptimalStepSize} is
feasible (i.e., $\tilde{\Delta}_{i,k}^{\ast}>0$) as long as $\beta_{i,k}%
^{\ast}\neq0$. \hfill $\blacksquare$

\begin{remark}\label{Rm:PathologicalSituation}
The extreme situation where $\beta_{i,k}^{\ast}=0$ occurs is when
$\sum_{j\in\bar{B}_{i}}d_{ij,k}=0$. However, since this \textquotedblleft
pathological situation\textquotedblright\ can be detected by agent $i$, if it
occurs, the agent can consider two options: 1) Use a reduced neighborhood $\bar
{B}_{i}^{1}\subset\bar{B}_{i}$ to calculate $\beta_{i,k}^{\ast}$ so that
$\beta_{i,k}^{\ast}\neq0$, hence $\tilde{\Delta}_{i,k}^{\ast}>0$, or 2) Use
the weighted form of \eqref{Eq:NeighborhoodCostFunctionTheory} (see Remark
\ref{Rm:WightedNeighborhoodObjective}) and manipulate the weight factors
$\{w_{ij}:j\in\bar{B}_{i}\}$ so as to get a step size $\beta_{i,k}^{\ast}%
\neq0$ (e.g., enforcing $w_{ij}=0,\forall j\ni d_{i,k}^{T}d_{ij}<0$ will give
$\beta_{i,k}>0$, hence $\tilde{\Delta}_{i,k}^{\ast}>0$).
\end{remark}

Regarding the term $Q_{i,k}$ in \eqref{Eq:Q} over which agent $i$ does not
have any control, let us first establish the following property.

\begin{lemma}
\label{Lm:QProperties1} The term $Q_{i,k}$ can be expressed as
\begin{equation}
Q_{i,k}=\sum_{j\in B_{i}}(\tilde{\Delta}_{j,k}+\sum_{l\in B_{j}-\{i\}}\left[
\Delta_{lj,k}-\Delta_{jl,k}\right]  ). \label{Eq:SimplifiedQ_i}%
\end{equation}
Further, if $B_{i}=\bar{B}_{j}-\{i\}$, then under \eqref{Eq:OptimalStepSize},
$Q_{i,k}>0$.
\end{lemma}
\emph{Proof: } In \eqref{Eq:Q}, let us add and subtract $\sum_{l\in
B_{j}-\{i\}}\Delta_{jl,k}$ to the inner terms of the main summation. Then,
using the definition \eqref{Eq:DeltaTildai}, the expression in
\eqref{Eq:SimplifiedQ_i} is obtained. To prove the second part, note that the
first inner term of the main summation of \eqref{Eq:SimplifiedQ_i} (i.e.,
$\tilde{\Delta}_{j,k}$) is always positive under the optimal step size given
in \eqref{Eq:OptimalStepSize}. Let us then consider the net effect of the
second inner term of $Q_{i,k}$, denoted by $Q_{i,k}^{\prime}$, where we have%
\[
Q_{i,k}^{\prime}=Q_{i,k}-\sum_{j\in B_{i}}\tilde{\Delta}_{j,k}=\sum_{j\in
B_{i}}\sum_{l\in B_{j}-\{i\}}\left[  \Delta_{lj,k}-\Delta_{jl,k}\right]  .
\]
Using the fact that $\Delta_{lj,k}-\Delta_{jl,k}=0$ when $l=j$, we can add a
dummy term into the inner summation to get
\[
Q_{i,k}^{\prime}=\sum_{j\in B_{i}}\sum_{l\in\bar{B}_{j}-\{i\}}\left[
\Delta_{lj,k}-\Delta_{jl,k}\right]  =\sum_{j\in B_{i}}\sum_{l\in B_{i}}\left[
\Delta_{lj,k}-\Delta_{jl,k}\right]  ,
\]
where the last step follows from the assumption $B_{i}=\bar{B}_{j}-\{i\}$.
Observing that the two running variables $l,j$ in the summations above are
interchangeable, we get $Q_{i,k}^{\prime}=0$. This implies that under
\eqref{Eq:OptimalStepSize}, $Q_{i,k}=\sum_{j\in B_{i}}\tilde{\Delta}_{j,k}>0$.
\hfill$\blacksquare$

We now make the following assumption regarding $Q_{i,k}$.

\begin{assumption}
\label{As:PositivityOfQTilde} Consider the sum,
\begin{equation}
\tilde{Q}_{i,k}=\sum_{l=k-T_{i}}^{k}Q_{i,l}, \label{Eq:QTilde}%
\end{equation}
such that $0\leq T_{i}\leq k$. Then, $\exists T_{i}<\infty$ such that
$\tilde{Q}_{i,k}\geq0$.
\end{assumption}

When the graph $\mathcal{G}(\mathcal{V},\mathcal{A})$ is complete, the condition $B_{i}=\bar{B}_{j}-\{i\}$ in Lemma \ref{Lm:QProperties1} is true for all $i\in\mathcal{V}$. In such cases, Assumption \ref{As:PositivityOfQTilde} is immediately satisfied with $T_{i}=1,\forall i\in\mathcal{V}$. On the other hand, when the graph $\mathcal{G}(\mathcal{V},\mathcal{A})$ is sparse enough, it can be considered as a collection of fully connected sub-graphs (exploiting the partitioned nature of local objective functions $H_{i}(\bar{s}_{i})$). Then, Assumption \ref{As:PositivityOfQTilde} also holds with $T_{i}=1,\forall i\in\mathcal{V}$. More generally, when each agent selects its step size according to \eqref{Eq:OptimalStepSize}, it ensures that $\tilde{\Delta}_{i,k}^{\ast}>0$. In addition, $\Delta_{ii,k}>0$ whenever the step size $\beta_{i,k}$ is positive. The assumption is further supported by the fact that each $Q_{i,k}$ in $\tilde{Q}_{i,k}$ is also a summation of $\Delta_{jj,k},$ $\Delta_{ji,k}$ and $\Delta_{lj,k}$ terms (noting in particular the positive first terms in \eqref{Eq:Q} as well as in \eqref{Eq:SimplifiedQ_i}). Moreover, it is locally verifiable if the agent communicates with its neighbors. In practice, we have never seen this assumption violated over extensive simulation examples (see Fig. \ref{Fig:PositivityOfQ_i} in Section \ref{Sec:CoverageControlApplication} and accompanying discussion).


Before establishing the convergence proof in Theorem \ref{Th:Convergence}, we
need one final technical condition.
\begin{assumption}
\label{As:Psi} For all $i\in\mathcal{V}$, there exists a function $\Psi_{i,k}$
such that $0<\Psi_{i,k}$ and \begin{numcases}{}
\label{Eq:Psi1}
0 \leq\Psi_{i,k}\Vert d_{i,k}\Vert^{2}\leq\tilde{\Delta}_{i,k}^{\ast}+\tilde{Q}_{i,k},
\mbox{ when } 0<\tilde{\Delta}_{i,k}^{\ast}+\tilde{Q}_{i,k},\ \ \ \   \\
\label{Eq:Psi2}
0 \leq\Psi_{i,k}\Vert d_{i,k}\Vert^{2}\leq\tilde{\Delta}_{i,k}^{\ast},
\mbox{ when } 0<\tilde{\Delta}_{i,k}^{\ast}.
\end{numcases}
\end{assumption}


This assumption is trivial because whenever the optimal step size in
\eqref{Eq:OptimalStepSize} is used, $0<\tilde{\Delta}_{i,k}^{\ast}$, hence,
for some $1<K_{2}$, $\Psi_{i,k}=\tilde{\Delta}_{i,k}^{\ast}/(K_{2}\Vert
d_{i,k}\Vert^{2})$ is a candidate function. Moreover, at time instants when
$0<\tilde{\Delta}_{i,k}^{\ast}+\tilde{Q}_{i,k}$ occurs, for some $1<K_{2}$,
$\Psi_{i,k}=(\tilde{Q}_{i,k}+\tilde{\Delta}_{i,k}^{\ast})/(K_{2}\Vert
d_{i,k}\Vert^{2})$ can be used as a candidate function for $\Psi_{i,k}$.

We can now state our main convergence theorem.

\begin{theorem}
\label{Th:Convergence} For all $i\in\mathcal{V}$ such that $\bar{B}%
_{i}\subseteq\mathcal{N}$, under Assumptions \ref{As:PropertiesofH_i}%
,\ref{As:LocalAvailabilityOfd_ij},\ref{As:PositivityOfQTilde}, and
\ref{As:Psi}, the step size selection in \eqref{Eq:OptimalStepSize} guarantees
the convergence criterion \eqref{Eq:RequiredConvergence1}, i.e.,
$\lim_{k\rightarrow\infty}d_{i,k}=0$.
\end{theorem}

\emph{Proof:} By Assumption \ref{As:PositivityOfQTilde}, a $T_{i}$ value can
be defined for $\tilde{Q}_{i,k}$ at each $k$. Consider a sequence of
consecutive discrete update instants $\{k_{1}+1,\ldots,k_{1}^{\prime}\}$ (for
short, we use the notation $(k_{1},k_{1}^{\prime}]$), where, $T_{i}%
=k_{1}^{\prime}-k_{1}$ is associated with $\tilde{Q}_{i,k_{1}^{\prime}}$ and
$T_{i}>k-k_{1}$ applies to all $\tilde{Q}_{i,k}$, $k\in(k_{1},k_{1}^{\prime
}-1]$. This means $0<\sum_{k=k_{1}+1}^{k^{\prime}}Q_{i,k}$ and $0\geq
\sum_{k=k_{1}+1}^{k}Q_{i,k},$ $\forall k\in(k_{1},k_{1}^{\prime}-1]$. In
addition, by Lemma \ref{Lm:AuxProb1}, $0<\tilde{\Delta}_{i,k}^{\ast}$ $\forall
k\in(k_{1},k_{1}^{\prime}]$. Thus, $0<\sum_{k=k_{1}+1}^{k_{1}^{\prime}}%
(\tilde{\Delta}_{i,k}^{\ast}+Q_{i,k})$. Now, by summing up both sides of
\eqref{Eq:DescentLemmaForH_Bi} over all update steps $k\in(k_{1},k_{1}%
^{\prime}]$ yields
\begin{equation}
\tilde{H}_{i}(\tilde{s}_{i,k_{1}^{\prime}+1})\geq\tilde{H}_{i}(\tilde
{s}_{i,k_{1}+1})+\sum_{k=k_{1}+1}^{k_{1}^{\prime}}(\tilde{\Delta}_{i,k}^{\ast
}+Q_{i,k}).\label{Eq:DescentLemmaFromk1tok1p}%
\end{equation}
Similarly, using Assumption \ref{As:Psi} and summing both sides of
\eqref{Eq:Psi2} over all $k\in(k_{1},k_{1}^{\prime}-1]$ and using
\eqref{Eq:Psi1} for $k=k_{1}^{\prime}$ yields
\begin{equation}
0\leq\sum_{k=k_{1}+1}^{k_{1}^{\prime}}\Psi_{i,k}\Vert d_{i,k}\Vert^{2}\leq
\sum_{k=k_{1}+1}^{k_{1}^{\prime}}(\tilde{\Delta}_{i,k}^{\ast}+Q_{i,k}%
).\label{Eq:Assumption4Fromk1tok1p}%
\end{equation}
By Assumption \ref{As:PositivityOfQTilde}, the length of the chosen interval
$(k_{1},k_{1}^{\prime}]$ is always finite. Therefore, any $\{1,\ldots,k_{2}\}$
with $k_{2}<\infty$ can be decomposed into a sequence of similar
sub-intervals: $\{(k_{11},k_{11}^{\prime}],(k_{12},k_{12}^{\prime}%
],\ldots,(k_{1L},k_{1L}^{\prime}]\}$ where $k_{11}=0,$ $k_{1i}^{\prime
}=k_{1(i+1)}$ $\forall i\in(0,L]$. If $k_{2}$ is such that $k_{1L}^{\prime
}<k_{2}$ (which happens if $0>\sum_{k=k_{1L}^{\prime}+1}^{k_{2}}Q_{i,k}$), Assumption \ref{As:PositivityOfQTilde} implies that there exists some
$k_{2}^{\prime}$ such that $k_{2}<k_{2}^{\prime}<\infty$ which satisfies
$0<\sum_{k=k_{1L}^{\prime}+1}^{k_{2}^{\prime}}Q_{i,k}$ (i.e., $(k_{1L}%
^{\prime},k_{2}^{\prime}]$ is the new last sub-interval of $(0,k_{2}^{\prime
}]$). Then, by writing the respective expressions in
\eqref{Eq:DescentLemmaFromk1tok1p} and \eqref{Eq:Assumption4Fromk1tok1p} for
each such sub-interval of the complete interval $(0,k_{2}^{\prime}]$ and
summing both sides over all $k$ yields
\begin{equation}
\tilde{H_{i}}(\tilde{s}_{i,k_{2}^{\prime}+1})\geq\tilde{H}_{i}(\tilde{s}%
_{i,1})+\sum_{k=1}^{k_{2}^{\prime}}(\tilde{\Delta}_{i,k}^{\ast}+Q_{i,k}%
),\label{Eq:DescentLemmaFrom0tok2}%
\end{equation}%
\begin{equation}
0\leq\sum_{k=1}^{k_{2}^{\prime}}\Psi_{i,k}\Vert d_{i,k}\Vert^{2}\leq\sum
_{k=1}^{k_{2}^{\prime}}(\tilde{\Delta}_{i,k}^{\ast}+Q_{i,k}%
),\label{Eq:Assumption4From0tok2}%
\end{equation}
respectively. Using Assumption \ref{As:PropertiesofH_i} in
\eqref{Eq:DescentLemmaFrom0tok2} gives $|\bar{B}_{i}|H_{UB}\geq\tilde{H_{i}%
}(\tilde{s}_{i,k_{2}^{\prime}+1})-\tilde{H}_{i}(\tilde{s}_{i,1})\geq\sum
_{k=1}^{k_{2}^{\prime}}(\tilde{\Delta}_{i,k}^{\ast}+Q_{i,k})$. Combining this
with \eqref{Eq:Assumption4From0tok2} yields
\begin{equation}
\sum_{k=1}^{k_{2}^{\prime}}\Psi_{i,k}\Vert d_{i,k}\Vert^{2}\leq|\bar{B}%
_{i}|H_{UB}.\label{Eq:Convergence1FinalStep}%
\end{equation}
By Assumption \ref{As:PropertiesofH_i}, the term $|\bar{B}_{i}|H_{UB}$ in
\eqref{Eq:Convergence1FinalStep} is a finite positive number. Also, by
Assumption \ref{As:Psi}, $\Psi_{i,k}>0,$ $\forall k$. Therefore, taking limits
of the above expression as $k_{2}^{\prime}\rightarrow\infty$ implies the
convergence criterion in \eqref{Eq:RequiredConvergence1} as long as the
optimal step sizes given by \eqref{Eq:OptimalStepSize} are used.
\hfill $\blacksquare$

\subsection{Convergence for agents $i$ such that $\tilde{B}_{i}\cap\mathcal{B}\neq\emptyset$}

In this case, at least some of the agents in $\tilde{B}_{i}$ are in boosting
mode, following \eqref{Eq:BoostedGradientAscent}. Following the same approach
as in Section \ref{subsec:ConvergenceWhenAllNormal}, we seek an optimal
variable step size selection scheme similar to \eqref{Eq:OptimalStepSize} so
as to ensure the convergence criteria given in \eqref{Eq:RequiredConvergence2}
and \eqref{Eq:RequiredConvergence3}. Compared to \eqref{Eq:DescentLemmaForH_i}
the ascent lemma relationship for $H_{i}({\bar{s}_{i,k}})$ takes the form:
\begin{equation}
H_{i}(\bar{s}_{i,k+1})\geq H_{i}(\bar{s}_{i,k})+\sum_{j\in\bar{B}_{i}%
\cap\mathcal{N}}\Delta_{ji,k}+\sum_{j\in\bar{B}_{i}\cap\mathcal{B}}\hat
{\Delta}_{ji,k},
\end{equation}
where $\Delta_{ji,k}$ for $j\in\mathcal{N}$ is the same as
\eqref{Eq:Contributionji} and we set
\begin{equation}
\hat{\Delta}_{ji,k}=\beta_{j,k}\hat{d}_{j,k}^{T}d_{ji,k}-\frac{K_{1i}}{2}%
\beta_{j,k}^{2}\Vert\hat{d}_{j,k}\Vert^{2}\in\mathbb{R}.
\end{equation}
Then, the ascent lemma for neighborhood objective function
$\tilde{H}_{i}(\tilde{s}_{i,k})$ can be expressed as
\begin{equation}
\tilde{H}_{i}(\tilde{s}_{i,k+1})\geq\tilde{H}_{i}(\tilde{s}_{i,k}%
)+(\tilde{\Delta}_{i,k}+Q_{i,k}),\label{Eq:DescentLemmaForH_BiHat}%
\end{equation}
with
\begin{equation}
\tilde{\Delta}_{i,k} \triangleq 1_{\{i\in\mathcal{N}\}}[\sum_{j\in\bar{B_{i}}}%
\Delta_{ij(k)}]+1_{\{i\in\mathcal{B}\}}[\sum_{j\in\bar{B_{i}}}\hat{\Delta
}_{ij(k)}],\label{Eq:DeltaTildaHat}%
\end{equation}
\begin{equation}%
\begin{split}
Q_{i,k} &  \triangleq \sum_{j\in B_{i}}(1_{\{j\in\mathcal{N}\}}[\Delta_{jj,k}%
+\Delta_{ji,k}]+1_{\{j\in\mathcal{B}\}}[\hat{\Delta}_{jj,k}+\hat{\Delta
}_{ji,k}]\\
&  +\sum_{l\in\{B_{j}-\{i\}\}}[1_{\{l\in\mathcal{N}\}}\Delta_{lj,k}%
+1_{\{l\in\mathcal{B}\}}\hat{\Delta}_{lj(k)}]),
\end{split}
\label{Eq:Qhat}%
\end{equation}
where $1_{\{\cdot\}}$ is the usual indicator function. Under this new
$\tilde{\Delta}_{i,k}$ in \eqref{Eq:DeltaTildaHat}, the same auxiliary problem
as in \eqref{Eq:AuxiliaryOptimizationProblem1} is used to determine the step
size $\beta_{i,k}^{\ast}$ to optimally increase the neighborhood cost function
$\tilde{H}_{i}(\tilde{s}_{k})$.

\begin{lemma}
\label{Lm:AuxProb1hat} The solution to the auxiliary problem
\eqref{Eq:AuxiliaryOptimizationProblem1} with $\tilde{\Delta}_{i,k}$ given in
\eqref{Eq:DeltaTildaHat} is
\begin{equation}
\beta_{i,k}^{\ast}=%
\begin{cases}
\frac{1}{\sum_{j\in\bar{B}_{i}}K_{1j}}\frac{d_{i,k}^{T}(\sum_{j\in\bar{B}_{i}%
}d_{ij,k})}{\Vert d_{i,k}\Vert^{2}}\mbox{ when }i\in\mathcal{N},\\
\frac{1}{\sum_{j\in\bar{B}_{i}}K_{1j}}\frac{\hat{d}_{i,k}^{T}(\sum_{j\in
\bar{B}_{i}}d_{ij,k})}{\Vert\hat{d}_{i,k}\Vert^{2}}\mbox{ when }i\in
\mathcal{B}.
\end{cases}
\label{Eq:OptimalStepSizehat}%
\end{equation}

\end{lemma}

\emph{Proof: } The proof follows the same steps as that of Lemma
\ref{Lm:AuxProb1} and is, therefore, omitted. \hfill $\blacksquare$

Note that the step size selection criterion in \eqref{Eq:OptimalStepSizehat}
for agent $i$ does not depend on its neighbors' modes. Thus, it offers a
generalization of \eqref{Eq:OptimalStepSize}. However, note that $\beta
_{i,k}^{\ast}$ now depends on its own mode. This is due to the fact that the
selection of $\beta_{i,k}^{\ast}$ allows agent $i$ to maximize the increment
in the neighborhood objective function $\tilde{H}_{i}(\tilde{s}_{i})$ which is
defined in \eqref{Eq:NeighborhoodCostFunctionTheory} independently from
boosting. Thus, the use of $\beta_{i,k}^{\ast}$ provides a regulation
mechanism for state update steps. 

To establish the convergence criteria \eqref{Eq:RequiredConvergence2} and
\eqref{Eq:RequiredConvergence3}, Assumptions \ref{As:PropertiesofH_i},
\ref{As:LocalAvailabilityOfd_ij} and \ref{As:PositivityOfQTilde} are still
required. Note that Assumption \ref{As:PositivityOfQTilde} should now be
considered under the new expression for $Q_{i,k}$ in \eqref{Eq:Qhat}; its
justification is similar as before. Moreover, a generalized version of Lemma
\ref{Lm:QProperties1} is given below.

\begin{lemma}
\label{Lm:QProperties2} The term $Q_{i,k}$ in \eqref{Eq:Qhat} can be expressed
as,
\begin{equation}%
\begin{split}
Q_{i,k}=\sum_{j\in B_{i}}(\tilde{\Delta}_{j,k}+ &  \sum_{l\in B_{j}%
-\{i\}}[1_{\{l\in\mathcal{N}\}}\Delta_{lj,k}-1_{\{j\in\mathcal{N}\}}%
\Delta_{jl,k}\\
&  \ \ \ \ \ \ +1_{\{l\in\mathcal{B}\}}\hat{\Delta}_{lj,k}-1_{\{j\in
\mathcal{B}\}}\hat{\Delta}_{jl,k}]).
\end{split}
\label{Eq:SimplifiedQ_ihat}%
\end{equation}
Further, if $B_{i}=\bar{B}_{j}-\{i\}$, then under
\eqref{Eq:OptimalStepSizehat}, $Q_{i,k}>0$.
\end{lemma}

\emph{Proof: } The proof follows the same steps as that of Lemma
\ref{Lm:QProperties1} and is, therefore, omitted. \hfill $\blacksquare$

Finally, Assumption \ref{As:Psi} needs to be modified into the following form
to incorporate the possibility that $i\in\mathcal{B}$.

\begin{assumption}
\label{As:PsiHat} For all $k$, there exists a function $\Psi_{i,k}$ such that
$0<\Psi_{i,k}$ and, if $i\in\mathcal{N}$:
\[%
\begin{cases}
0\leq\Psi_{i,k}\Vert d_{i,k}\Vert^{2}\leq\tilde{\Delta}_{i,k}^{\ast}+\tilde
{Q}_{i,k} & \mbox{ when }0<\tilde{\Delta}_{i,k}^{\ast}+\tilde{Q}_{i,k},\\
0\leq\Psi_{i,k}\Vert d_{i,k}\Vert^{2}\leq\tilde{\Delta}_{i,k}^{\ast} &
\mbox{ when }0<\tilde{\Delta}_{i,k}^{\ast},
\end{cases}
\]
otherwise, if $i\in\mathcal{B}$:
\[%
\begin{cases}
0\leq\Psi_{i,k}\Vert\hat{d}_{i,k}\Vert^{2}\leq\tilde{\Delta}_{i,k}^{\ast
}+\tilde{Q}_{i,k} & \mbox{ when }0<\tilde{\Delta}_{i,k}^{\ast}+\tilde{Q}%
_{i,k},\\
0\leq\Psi_{i,k}\Vert\hat{d}_{i,k}\Vert^{2}\leq\tilde{\Delta}_{i,k}^{\ast} &
\mbox{ when }0<\tilde{\Delta}_{i,k}^{\ast}.
\end{cases}
\]
Here, $\tilde{Q}_{i,k}$ is evaluated from \eqref{Eq:QTilde} using
\eqref{Eq:Qhat} and, $\tilde{\Delta}_{i,k}^{\ast}$ from
\eqref{Eq:DeltaTildaHat} using \eqref{Eq:OptimalStepSizehat}.
\end{assumption}

The following theorem can now be established.

\begin{theorem}
\label{Th:ConvergenceHat} Under Assumptions \ref{As:PropertiesofH_i}%
,\ref{As:LocalAvailabilityOfd_ij},\ref{As:PositivityOfQTilde}, and
\ref{As:PsiHat}, the step size selection in \eqref{Eq:OptimalStepSizehat}
guarantees the convergence conditions stated in
\eqref{Eq:RequiredConvergence1}-\eqref{Eq:RequiredConvergence3}: if
$i\in\mathcal{N},\ $then $\lim_{k\rightarrow\infty}d_{i,k}=0$, and, if
$i\in\mathcal{B}$, then$\ \lim_{k\rightarrow\infty}\hat{d}_{i,k}=0$.
\end{theorem}

\emph{Proof: } The proof uses the same steps as in that of Theorem
\ref{Th:Convergence}. The only difference lies in the use of new terms for
$\tilde{\Delta}_{i,k}$, $\tilde{\Delta}_{i,k}^{\ast}$ and $Q_{i,k}$, given by
\eqref{Eq:DeltaTildaHat}, \eqref{Eq:OptimalStepSizehat} and \eqref{Eq:Qhat}.
Then, the final step of the proof is
\begin{equation}
\sum_{k=1}^{k_{2}^{\prime}}\Psi_{i,k}[1_{\{i\in\mathcal{N}\}}\Vert
d_{i,k}\Vert^{2}+1_{\{i\in\mathcal{B}\}}\Vert\hat{d}_{i,k}\Vert^{2}]\leq
|\bar{B}_{i}|H_{UB}.\label{Eq:Convergence1HatFinalStep}%
\end{equation}
By Assumption \ref{As:PropertiesofH_i}, the R.H.S. of the above expression is
finite and positive. Taking limits when $k_{2}^{\prime}\rightarrow\infty$
yields convergence criteria given in \eqref{Eq:RequiredConvergence2} and
\eqref{Eq:RequiredConvergence3}. Further, noting that Theorem
\ref{Th:ConvergenceHat} is a generalization of Theorem \ref{Th:Convergence}
with the step size selection scheme \eqref{Eq:OptimalStepSizehat} replacing
\eqref{Eq:OptimalStepSize}, \eqref{Eq:RequiredConvergence1} is also satisfied.
\hfill $\blacksquare$

\subsection{Discussion}

\subsubsection{Extending to dynamic graphs} \label{SubSubSec:FixedGraphsToVaryingGraphs}
Both the considered main problem \eqref{Eq:GeneralizedOptimizationProb} and the formulated variable step size method \eqref{Eq:OptimalStepSizehat} assumes that agents are inter-connected (i.e., inter agent communications occur) according to a fixed graph topology $\mathcal{G}$. As pointed out earlier, the nature of $\mathcal{G}$ can affect Assumption \ref{As:PositivityOfQTilde} (specifically through the $T_i$ value) of the convergence proof. Nevertheless However, due to the nature of the used convergence proof, it is reasonable to expect that the developed variable step size method (i.e., the Theorem \ref{Th:ConvergenceHat}) is extendable to cases where the graph $\mathcal{G}$: 
(\romannum{1}) Varies sufficiently slower than the convergence rate, and, (\romannum{2}) Converges to an asymptotic graph configuration. 
In fact, the coverage control application which will be used to demonstrate the proposed solution technique (in Section \ref{Sec:CoverageControlApplication}) belongs to the latter case. Moreover, since the variable step sizes in \eqref{Eq:OptimalStepSizehat} leads each agent to maximize the improvement of their neighborhood objective function, we can expect them to converge even when the graph $\mathcal{G}$ varies rapidly (however without showing any oscillatory behavior).

\subsubsection{Feasible space constraint} \label{SubSubSec:EffectOfTheFeasibleSpace}
The considered main problem in \eqref{Eq:GeneralizedOptimizationProb} includes a feasible space ($\mb{F}$) constraint for the global state $\mb{s}$. However, to simplify the convergence analysis process (discussed above), it has not taken into account so far. This omission is further justifiable because, even if an agent hits a constraint during its state update process (\eqref{Eq:GradientAscentGeneral} or \eqref{Eq:BoostedGradientAscent}), it can always resort to a standard gradient projection method \cite{Bertsekas2016nonlinear}. Moreover, for a such situation, the following lemma presents an additional condition which needs to be satisfied in order to guarantee the convergence of the proposed variable step size method \eqref{Eq:OptimalStepSizehat}. 

\begin{lemma}\label{Lm:Eq:ConvergenceConditionForProjectedGradient}
If feasible space $\mb{F}$ is convex, and if an agent $i$'s local and cross gradients satisfy the conditions,
\begin{equation}\label{Eq:ConvergenceConditionForProjected}
\begin{aligned}
\vert d_{i,k}^T \sum_{j \in B_i} d_{ij,k} \vert &< \Vert d_{i,k} \Vert^2
\mbox{ when } i \in \mathcal{N},\\
\vert \hat{d}_{i,k}^T (\sum_{j \in B_i} d_{ij,k} + (d_{i,k}-\hat{d}_{i,k}) ) \vert &< \Vert \hat{d}_{i,k} \Vert^2
\mbox{ when } i \in \mathcal{B},
\end{aligned}
\end{equation}
the step sizes $\beta_{i,k} = \beta_{i,k}^*$ given by \eqref{Eq:OptimalStepSizehat} when used in \eqref{Eq:GradientAscentGeneral} or \eqref{Eq:BoostedGradientAscent} with appropriate gradient projections (onto $\mb{F}$), will lead the state $s_{i,k}$ to a stationary point (i.e., to the convergence). 
\end{lemma}

\emph{Proof: } Consider the problem where the neighborhood objective function $\tilde{H}_i(\tilde{s}_{i,k})$ needs to be maximized using the projected state updates of $s_{i,k}$ on the convex feasible space $\mb{F}$. For this situation, according to \cite{Bertsekas2016nonlinear}, the convergence condition on the step sizes $\beta_{i,k}$ is $0<\beta_{i,k}<\frac{2}{K_i}$, where $K_i$ is the Lipschitz constant of $\nabla \tilde{H}_i$. Note that we can write $K_i = \sum_{j\in \bar{B}_i}K_{1j}$ due to \eqref{Eq:NeighborhoodCostFunctionTheory}. Also, for $i \in \mathcal{N}$, $\beta_{i,k}^*$ expression given in \eqref{Eq:OptimalStepSizehat} can be modified into the form,
\begin{equation}\label{Eq:SimilifiedOptimumStepSize}
    \beta_{i,k}^* = \frac{1}{K_i} \left[ 1 + \frac{d_{i,k}^T \sum_{j\in B_i}d_{ij,k}}{\Vert d_{i,k} \Vert^2} \right].
\end{equation}
Now, by enforcing the convergence condition: $0 < \beta_{i,k}^* <\frac{2}{K_i}$ yields the first condition in \eqref{Eq:ConvergenceConditionForProjected}. Similarly the second condition in \eqref{Eq:ConvergenceConditionForProjected} can be obtained when the $\beta_{i,k}^*$ expression for $i\in \mathcal{B}$ in \eqref{Eq:OptimalStepSizehat} is considered. \hfill $\blacksquare$

From a practical standpoint, during the gradient ascent, if the projections play a major role, it is better to check the conditions stated in Lemma \ref{Lm:Eq:ConvergenceConditionForProjectedGradient}. If they are being violated, the neighborhood reduction technique and/or the weight factor manipulation techniques mentioned in Remark \ref{Rm:PathologicalSituation} can be used to change the $B_i$ and/or $\tilde{H}_i$ respectively so that the conditions in Lemma \ref{Lm:Eq:ConvergenceConditionForProjectedGradient} are satisfied.

In fact, the main reason behind the inclusion of the feasible space constraint $\mb{s} \in \mb{F}$ in \eqref{Eq:GeneralizedOptimizationProb} is that it can play an important role in designing boosting functions $f_i$. For example, a boosted gradient can be constructed such that $\hat{d}_i = f_i(d_i,\mb{F})$ with some some special features of $\mb{F}$ are being utilized. For more details see the V-Boosting and Arc-Boosting methods introduced for coverage problems discussed in Section \ref{Sec:CoverageControlApplication}.

\subsubsection{Variable step sizes compared to fixed step sizes}
Typically, in a centralized setting, using a fixed step size for the gradient descent is computationally cheap, and, if done correctly, it should deliver a higher convergence rate compared to variable step size methods. However, in a distributed setting where agents are independently and intermittently alter the followed gradient direction (\eqref{Eq:GradientAscentGeneral} and \eqref{Eq:BoostedGradientAscent}), using a fixed step size (typically $\beta_{i,k} = \frac{1}{K_i}$) might not lead to good overall convergence properties. Further, establishing the convergence for such fixed step size approach is a challenging task - without making restrictive (non-trivial) assumptions. In contrast, the proposed variable step
size method has the following advantages: (\romannum{1}) It is designed so as
to account for the distributed and cooperative nature of the underlying
problem, (\romannum{2}) Its convergence has been established by making only a
few locally verifiable assumptions, (\romannum{3}) It is not computationally heavy compared
to line search methods, and, (\romannum{4}) During different modes
(boosting/normal) the step sizes are automatically adjusted.
As a result of these positive traits, in applications, the variable step size method showed better convergence results compared to fixed step size methods (see Sections \ref{SubSec:VariableStepSizesApplicationExample} and \ref{Sec:CoverageControlApplication}). 

\subsubsection{Termination conditions for modes} \label{SubSubSec:TerminationConditionsForImplementations}
In applications, the equilibrium conditions $d_i = 0$ and $\hat{d}_i = 0$ used in boosting schemes should be replaced with appropriate termination conditions \cite{Bertsekas2016nonlinear} such as $\Vert d_i \Vert \leq \epsilon_1$ and $\Vert \hat{d}_i \Vert \leq \epsilon_2$ (respectively) where $\epsilon_1,\epsilon_2$ are two chosen small positive scalars.  

\subsubsection{Escaping and converging to saddle points}
Due to the non-convexity of the objective function, saddle points may exist in the feasible space.
However, as shown in \cite{Lee2017,Panageas2019}, first-order methods \eqref{Eq:GradientAscentGeneral}
almost always avoid a large class of saddle points (called strict saddle points) inherently. Nevertheless, if boosting
functions are deployed through \eqref{Eq:BoostedGradientAscent}, clearly, saddle points are easier to escape
from compared to local minima. Moreover, even if the convergence criteria
\eqref{Eq:RequiredConvergence1} - \eqref{Eq:RequiredConvergence3} lead to a saddle point, it will have a
higher cost compared to initially attained local minima (or saddle points) as a result of the comparison stage used in
boosting schemes (e.g., see \say{$H^B > H$} block in Fig. \ref{Fig:CentralizedBoostingBlockDiagram}).

\subsection{
An application example for the variable step size method}
\label{SubSec:VariableStepSizesApplicationExample}

In this section, a simple example is provided to illustrate the operation and convergence (i.e., validity) of the proposed variable step size method. In this example, local objective functions are restricted to take a quadratic form,
\begin{equation}
    H_i(\bar{s}_i) = -\Vert \sum_{j\in \bar{B}_i} A_{ij}s_j - b_i \Vert^2_{C_i}
    = -\Vert g_i(\bar{s}_i)\Vert^2_{C_i}
\end{equation}
where $A_{ij} \in \R^{r \times m}, b_i \in \R^r$ and $C_i \in \R^{r \times r}$ for any $i \in \mathcal{V}, j \in \bar{B}_i$. The weighting matrix $C_i$ is symmetric and positive definite. The weighted norm is defined as $\Vert v \Vert^2_C = v^T C v$ with $v \in \R^r$ and $C \in \R^{r \times r}$. The parameter $r$ represents the dimension of the local cost function. Also, note that $ g_i(\bar{s}_i) = \sum_{j\in \bar{B}_i} A_{ij}s_j - b_i$. Assuming the parameters $A_{ij},b_i,C_i, \ \forall i \in \mathcal{V}, \ \forall j \in \bar{B}_i$ and the graph $\mathcal{G} = (\mathcal{V},\mathcal{E})$ are predefined (also given the specific $N, m$ and $r$ value combination), the interested optimization problem is,
\begin{equation}
    \mb{s}^* = [s_1^*,s_2^*,\ldots,s_N^*] = \arg\max_{\mb{s}} H(\mb{s}) = \sum_{i = 1}^N H_i(\bar{s}_i).
\end{equation}
Due to the quadratic nature of the associated objective functions, a closed form expression can be obtained for the global optimum $\mb{s}^*$. Also, as a result of the convexity, we do not need to use a boosting functions approach in this case. Therefore, we use this example to compare the performance of the proposed variable step size method (when used in a distributed gradient ascent), with respect to a fixed step size method (when used in a centralized gradient ascent).

For the (distributed) variable step size computation (at agent $i$ using \eqref{Eq:OptimalStepSize}), the local gradient $d_{i,k}$ is 
\begin{equation}
    d_{i} =\frac{\partial H_i(\bar{s}_i)}{\partial s_i} = -2 A_{ii}^T C_i g_i(\bar{s}_i),
\end{equation}
the cross gradients $d_{ij,k},\ \forall j \in \bar{B}_i$,
\begin{equation}
    d_{ij} = \left[\frac{\partial H_j(\bar{s}_j)}{\partial s_i}\right]_{i\leftrightarrow j} 
    = - 2 A_{ji}^T C_j (\sum_{l \in \bar{B}_j} A_{jl}s_l - b_j),
\end{equation}
and, the local Lipschitz constants $K_{1j}, j \in \bar{B}_i$,
\begin{equation}
    K_{1j} = 2\Vert A_j^T C_j A_j\Vert_\infty ,\ A_j = [\{A_{jl}\}_{l\in \bar{B}_j}] \in \R^{r\times m\vert \bar{B}_j\vert},
\end{equation}
are used. In contrast, for the use in centralized gradient ascent, the global gradient component of agent $i$, $d_{i,k}^G$ where
\begin{equation}
    d_{i}^G = \frac{\partial H(\mb{s})}{\partial s_i} = -2 \sum_{j \in \bar{B}_i} A_{ji}^T C_j g_j(\bar{s}_j)
\end{equation}
is used as a replacement for $d_{i,k}$ in \eqref{Eq:GradientAscentGeneral}. In there, the step size is kept as fixed value at $\frac{1}{K_i}$ where $K_i = \sum_{j \in \bar{B}_i} K_{1j}$ according to \cite{Bertsekas2016nonlinear} (also see Remark \ref{SubSubSec:EffectOfTheFeasibleSpace}). Finally, in order to assess the convergence, we use relative error profile $e_k$ \cite{Bastianello2018CDC} where   
\begin{equation}
    e_k = \log \left[ \frac{1}{N} \sum_{i=1}^N \frac{\Vert s_{i,k}-s_i^* \Vert}{\Vert s_i^* \Vert}\right],
\end{equation}
over a simulation (i.e, for a single realization). 

\begin{figure}[h]
\centering
\begin{subfigure}{0.49\columnwidth}
\includegraphics[width=\textwidth]{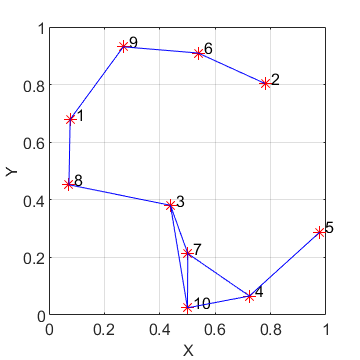}
\caption{Graph $\mathcal{G} = (\mathcal{V},\mathcal{E})$}
\end{subfigure}
\hfill
\begin{subfigure}{0.49\columnwidth}
\includegraphics[width=\textwidth]{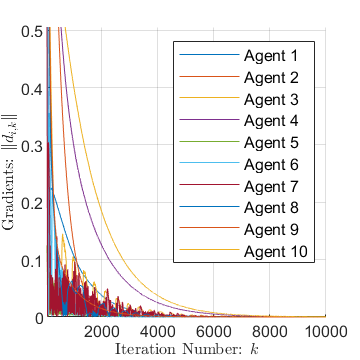}
\caption{Agent Local Derivatives}
\end{subfigure}\hfill
\begin{subfigure}{0.60\columnwidth}
\includegraphics[width=\textwidth]{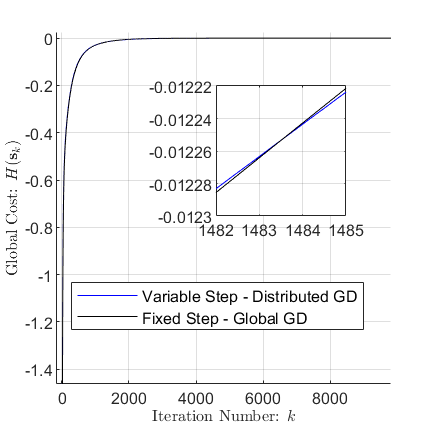}
\caption{Global Objective Value}
\end{subfigure}\hfill
\caption{Numerical Example 1}%
\label{Fig:VariableStepSizesExperiment2}%
\vspace{-3mm}
\end{figure}

\begin{figure}[h]
\centering
\begin{subfigure}{0.49\columnwidth}
\includegraphics[width=\textwidth]{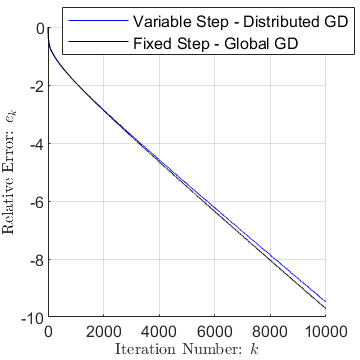}
\caption{A Relative Error Profile}
\end{subfigure}
\hfill
\begin{subfigure}{0.49\columnwidth}
\includegraphics[width=\textwidth]{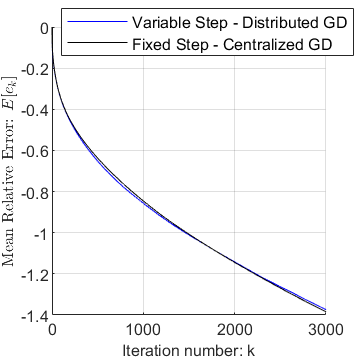}
\caption{Mean Relative Error Profile (computed over 100 realizations)}
\end{subfigure}
\caption{Relative Error Profiles for Numerical Example 1.}%
\label{Fig:VariableStepSizesExperiment2Error}%
\vspace{-3mm}
\end{figure}

\begin{figure}[!t]
\centering
\begin{subfigure}{0.49\columnwidth}
\includegraphics[width=\textwidth]{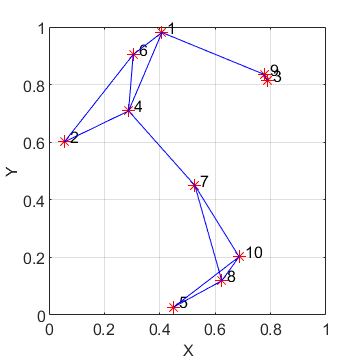}
\caption{Graph $\mathcal{G} = (\mathcal{V},\mathcal{E})$}
\end{subfigure}
\hfill
\begin{subfigure}{0.49\columnwidth}
\includegraphics[width=\textwidth]{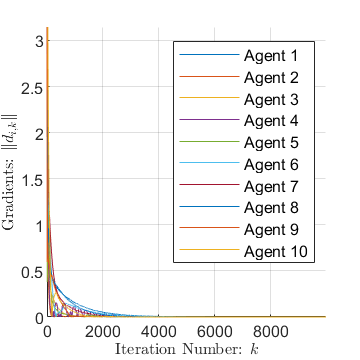}
\caption{Local Gradient Values}
\end{subfigure}
\hfill
\begin{subfigure}{0.49\columnwidth}
\includegraphics[width=\textwidth]{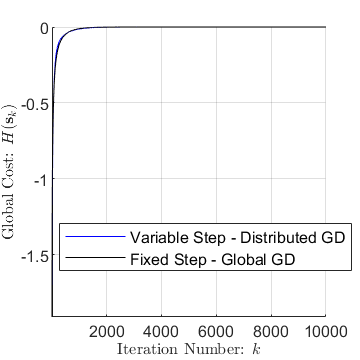}
\caption{Global Objective Value}
\end{subfigure}\hfill
\begin{subfigure}{0.49\columnwidth}
\includegraphics[width=\textwidth]{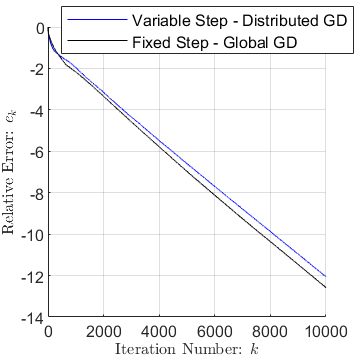}
\caption{A Relative Error Profile}
\end{subfigure}
\caption{Numerical Example 2}%
\label{Fig:VariableStepSizesExperiment3}%
\vspace{-3mm}
\end{figure}

\begin{figure}[!t]
\centering
\begin{subfigure}{0.49\columnwidth}
\includegraphics[width=\textwidth]{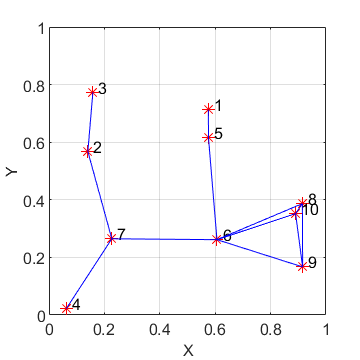}
\caption{Graph $\mathcal{G} = (\mathcal{V},\mathcal{E})$}
\end{subfigure}
\hfill
\begin{subfigure}{0.49\columnwidth}
\includegraphics[width=\textwidth]{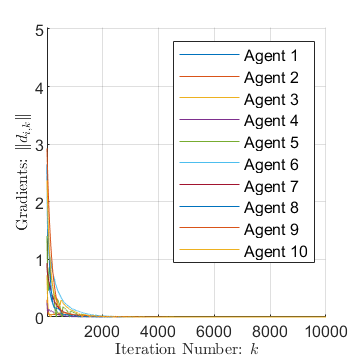}
\caption{Local Gradient Values}
\end{subfigure}
\hfill
\begin{subfigure}{0.49\columnwidth}
\includegraphics[width=\textwidth]{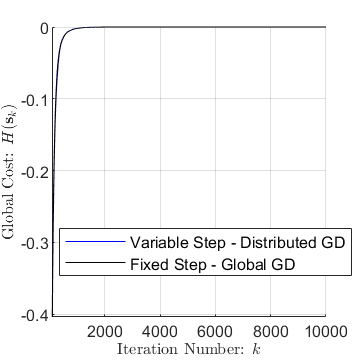}
\caption{Global Objective Value}
\end{subfigure}\hfill
\begin{subfigure}{0.49\columnwidth}
\includegraphics[width=\textwidth]{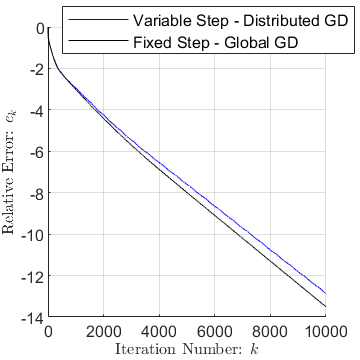}
\caption{A Relative Error Profile}
\end{subfigure}
\caption{Numerical Example 3}%
\label{Fig:VariableStepSizesExperiment4}%
\vspace{-3mm}
\end{figure}

In simulations, fixed dimensional parameters $N=10$ and $m=r=2$ are used. Note that $m=r$ is required here to guarantee the existence of a solution where $d_{i}=0,\forall i\in\mathcal{V}$. It it is easy to show that the optimal global objective function value is $H(\mathbf{s}^{\ast})=0$. To generate the inter-agent connections (i.e., the graph $\mathcal{G}$) a random geometric graph generation is used taking $0.4$ as the communication range parameter \cite{Bastianello2018CDC}. The remaining problem parameters $A_{ij},b_{i},C_{i},s_{i,0}\ \forall i\in\mathcal{V},\ \forall j\in\bar{B}_{i}$ are generated randomly (keeping the graph $\mathcal{G}$ fixed).

The experimental results shown in Figs. \ref{Fig:VariableStepSizesExperiment2}-
\ref{Fig:VariableStepSizesExperiment4} (corresponding to three different experiments) confirm our theoretical results regarding the convergence. The $H(\mathbf{s}_{k})$ profiles in Fig. \ref{Fig:VariableStepSizesExperiment2}(c) show that the proposed distributed variable step size method provides a slightly faster convergence than the centralized fixed step size method for $k\leq1483$ where at $k=1483$, the $H(\mathbf{s}_{k})$ value is $99.95\%$ closer to the optimal than the initial value $H(\mathbf{s}_{0})=26.1432$; for $k\geq1484$, the centralized fixed step method is slightly faster. This cross-over behavior can be understood as a result of local gradients $d_{i,k}$ becoming smaller as $k$ increases and adapting step sizes $\beta_{i,k}$ in \eqref{Eq:GradientAscentGeneral} when $d_{i,k}$ is very small is less effective. This cross-over behavior can also be seen in numerical examples shown in Figs. \ref{Fig:VariableStepSizesExperiment3} and \ref{Fig:VariableStepSizesExperiment4}. In all, our general observation over extensive similar examples is that the result of such a comparison (between distributed variable step and centralized fixed step methods) depends on the network topology




\section{Application to coverage control problem}

\label{Sec:CoverageControlApplication}

This section uses the class of multi-agent coverage control problems to illustrate: (\romannum{1}) the boosting functions related concepts introduced in Section \ref{Sec:GeneralProblemFormulation}, and, (\romannum{2}) the optimal variable step size selection mechanism proposed in Section \ref{Sec:VariableStepSizeScheme}.

We use the preliminary work regarding the coverage control problems presented in \cite{Zhong2011} where a distributed gradient based solution scheme has been proposed. The work in \cite{Sun2014} extends the solution proposed in \cite{Zhong2011} by adding the capability to escape local optima through a centralized boosting scheme - without a convergence analysis. In contrast, this section uses the developed theory for the class of general cooperative multi-agent optimization problems (discussed in Section \ref{Sec:GeneralProblemFormulation} and \ref{Sec:VariableStepSizeScheme}) to construct a convergence guaranteed distributed boosting scheme for the class of multi-agent coverage control problems.

Under this section, subsections \ref{SubSec:CoverageControlProblemFormulation} and \ref{SubSec:CoverageControlDistributedOptimization} presents the basic coverage control problem formulation and its distributed gradient based solution technique as proposed in \cite{Zhong2011}, along with few improvements. These improvements aim to: (\romannum{1}) Incorporate agents with limited sensing range, (\romannum{2}) Establish the convergence for a situation where boosting is not used, and, (\romannum{3}) Propose a mechanism to compute the Lipschitz constants associated with each agent's local objective function. 

In the first halves of the subsections \ref{SubSec:CoverageControlBoostingFunctions} and \ref{SubSec:CoverageControlBoostingSchemes}, boosting function families and the centralized boosting scheme proposed in \cite{Sun2014} are reviewed respectively. Then, in the second halves, two novel boosting function families and a novel distributed boosting scheme are presented, respectively. Then, Section \ref{SubSec:CoverageControlDBSConvergence} presents the convergence analysis of the proposing DBS.

\subsection{Basic coverage control problem formulation}
\label{SubSec:CoverageControlProblemFormulation} 

The coverage control problem aims to find an optimal arrangement for a given set of agents (sensor nodes) inside a given \textit{mission space} so as to maximize the probability of detecting randomly occurring events. It is assumed that the agent sensing capabilities, characteristics of the mission space, and any priori information on the spacial likelihood of random event occurrences (in the mission space) are fixed and known beforehand.

The mission space $\Omega\subset\mathbb{R}^{2}$ is modeled as a non-self-intersecting polygon - which is a polygon with no intersections between any two non-consecutive edges. The mission space may contain a finite set of non-self-intersecting polygonal obstacles denoted by $\{M_{1},M_{2},\ldots,M_{h}\}$, where, $M_{i}\subset\mathbb{R}^{2}$ represents the interior space of the $i$\textsuperscript{th} obstacle. Therefore, agent motion and deployment are constrained to a non-convex \textit{feasible space} $F = \Omega\backslash(\cup_{i=1}^{h} M_{i}$).

In order to quantify the spacial likelihood of random event occurrence in the mission space, an \textit{event density} function $R:\Omega\rightarrow \mathbb{R}$ is used. Typically, $R(x)=0,\ \forall x\not \in F$; $R(x)\geq0,\forall x\in\Omega$, and $\int_{\Omega}R(x)dx<\infty$ are assumed. Further, if no advance information is available, then $R(x)=1,\forall x\in\Omega$ is used. Furthermore, it is assumed that when an event occurs, it will emit a signal enabling it to be detected by nearby agents.

The mission space is considered to have $N$ agents. At a given update instant $k$ (discrete), the position coordinates of agent $i$ (i.e., the controllable \textit{local state}) is denoted by $s_{i,k} \in F \subset \R^2$. Therefore, the \textit{global state} of the multi-agent system is denoted by $\mb{s}_{k} = [s_{1,k},s_{2,k},\ldots,s_{N,k}]$. We write $\mb{s}_k \in \mb{F}$ to denote $s_{i,k} \in F\ \forall i$. For notational convenience, the update instant subscript $k$ is omitted unless it is important.

\begin{figure}[tbh]
\centering
\includegraphics[width=3in]{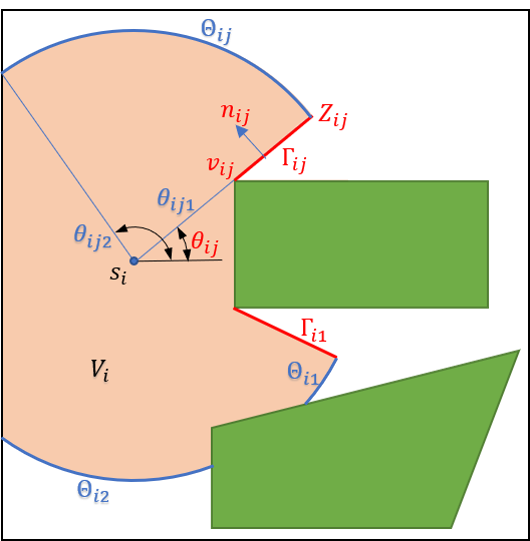} \caption{Mission space with one
agent}%
\label{Fig:MissionSpaceGeometry}%
\end{figure}

The sensing capabilities of agent $i$ is assumed to have two main physical characteristics: (\romannum{1}) Beyond a finite \textit{sensing radius} $\delta_{i} \in \R$ (from agent location $s_i$) it cannot detect any events, (\romannum{2}) Similar to a vision sensor, the presence of obstacles hinder the sensing capability of the. Considering these two factors, a \textit{visibility region} for agent $i$ is defined as $V_{i} = \{x: \Vert x-s_{i}\Vert\leq\delta_{i}, \forall\lambda\in(0,1], (\lambda x+(1-\lambda)s_{i}) \in F\}$. Fig. \ref{Fig:MissionSpaceGeometry} is provided to identify all associated geometric parameters in this model.

Consider an event denoted by $E_{i,x}$ = \textquotedblleft Agent $i$ detecting an event occurring at $x$\textquotedblright. A \textit{sensing function} $\hat{p}_{i}(x,s_{i})$ is used to quantify the probability of such an event (i.e. $Prob(E_{i,x})$). However, when the aforementioned sensing capability characteristics are incorporated, $\hat{p}_{i}(x,s_{i})$ takes the form
\begin{equation}
Prob(E_{i,x}) = \hat{p}_{i}(x,s_{i}) =
\begin{cases}
p_{i}(x,s_{i}) & \mbox{ if } x \in V_{i},\\
0 & \mbox{ if } x \not \in V_{i},
\end{cases}
\end{equation}
where, $p_{i}(x,s_{i})$ is defined so that $p_{i}:\R^2 \times \R^2 \rightarrow \R$ and is differentiable and monotonically decreasing in $D_{i}(x) \equiv \Vert x-s_{i} \Vert$. As an example, $p_{i}(x,s_{i}) = p_{0i}e^{-\lambda_{i}D_{i}(x)}$ or $p_{i}(x,s_{i}) = p_{0i}\Vert D_{i}(x) \Vert ^{-\lambda_{i}}$ represents two such typical choices for $p_{i}(x,s_{i})$. However, note that $\hat{p}_{i}(x,s_{i})$ can be strictly discontinuous w.r.t. $x$, $s_{i}$ or $D_{i}(x)$.

Now, for a given position $x \in \Omega$, assuming the set of events $E_{x} =
\{E_{i,X}: i=1,2,\ldots,N\}$ are independent from each other (i.e., assuming independently detecting agents), the probability of
happening at least one of the events in $E_{x}$ is defined as the
\textit{joint detection probability} $P(x,\mathbf{s})$ where
\begin{equation} \label{Eq:P}
P(x,\mathbf{s}) = 1 - \prod_{i=1}^{N}[1-\hat{p}_{i}(x,s_{i})].
\end{equation}

Combining the event density and joint detection probability, the objective function $H(\mb{s})$ of the coverage control problem given in \cite{Zhong2011} is 
\begin{equation}
\label{Eq:GlobalObjectiveFunction}
H(\mathbf{s}) = \int_{F} R(x)P(x,\mathbf{s})dx.
\end{equation}

\begin{remark}\label{Rm:CoverageObjectiveFunctionForm}
The form of the objective function in \eqref{Eq:GlobalObjectiveFunction} is not limited to coverage control problems. For example, consider a situation where $R(x)$ represents the value associated with the point $x$, and, $P(x,\mb{s})$ represents the interaction between the multi agent system $\mb{s}$ and the point $x$. For such paradigms, the same objective function form in \eqref{Eq:GlobalObjectiveFunction} can be used. Therefore, the ensuing discussion can be extended for such applications as well. 
\end{remark}

The underlying multi-agent optimization problem is 
\begin{equation}
\label{Eq:OptimizationProblemCoverageControl}
\mathbf{s}^{*} = \arg\max_{\mathbf{s} \in\mathbf{F}} H(\mathbf{s}),
\end{equation}
where $\mathbf{s}^\ast$ represents the optimal agent placement.

Note that the objective function in \eqref{Eq:GlobalObjectiveFunction}
is non-linear and non-convex, while the feasible space $\mb{F}$ is also non-convex. Therefore, the coverage control problem posed in
\eqref{Eq:OptimizationProblemCoverageControl} has an identical structure (and
qualities) to the considered general cooperative multi-agent optimization
problem in \eqref{Eq:GeneralizedOptimizationProb}. Therefore, \eqref{Eq:OptimizationProblemCoverageControl} can have multiple local optimal solutions (even in the simplest configurations). Hence, the application of distributed boosting functions approach can aid the agents to escape local optima while solving \eqref{Eq:OptimizationProblemCoverageControl}.

As discussed in Section \ref{Sec:GeneralProblemFormulation}, if a distributed
gradient based method is to be used to solve \eqref{Eq:OptimizationProblemCoverageControl} (i.e., to get to a local optimum), such a method can be helpful in constructing a distributed boosting functions approach (to escape local optima). Therefore, as the next step, let's discuss the distributed gradient based method used to solve \eqref{Eq:OptimizationProblemCoverageControl}.

\subsection{Distributed optimization solution}
\label{SubSec:CoverageControlDistributedOptimization}

In coverage control problems, two agents are considered to be \textit{neighbors} if their visibility regions overlap \cite{Sun2014}. According to this notion of neighbors, the usual sets representing the
\textit{neighborhood} $B_{i}$ and the \textit{closed neighborhood} $\bar
{B}_{i}$ of an agent $i$ are defined as $B_{i} = \{j:V_{j} \cap V_{i}
\neq\emptyset, i \neq j \}$ and $\bar{B}_{i} = B_{i} \cup\{i\}$
respectively. It is assumed that agents share their local state
information $s_{i}$ with their neighbors, so that each agent has knowledge
of its \textit{neighborhood state} $\bar{s}_{i} = \{s_{j}:j \in\bar{B}_{i}\}$.
We use an undirected graph $\mathcal{G}=(\mathcal{V},\mathcal{A})$ to model inter-agent interactions, where $\mathcal{V} = \{1,2,\ldots,N\}$ and $\mathcal{A} = \{(i,j): \forall i,j \in \mathcal{V},\ i \neq j,\ j \in B_i\}$ (same notations and definitions as before).

Under these definitions, it is shown in \cite{Sun2014} that the coverage
control global objective $H(\mathbf{s})$ in
\eqref{Eq:GlobalObjectiveFunction} can be partitioned as $H(\mathbf{s}) = H_{i}(\bar{s}_{i}) + H_{i}^{c}(s_{i}^{c})$, where 
\begin{equation}
\label{Eq:LocalObjectiveCoverageControl}H_{i}(\bar{s}_{i}) = \int_{V_{i}}R(x)
\prod_{j \in B_{i}} \left[  1-\hat{p}_{j}(x,s_{j})\right]  p_{i}(x,s_{i})dx,
\end{equation}
and
\begin{equation}
H_{i}^{c}(s_{i}^{c}) = \int_{F}R(x)( 1 - \prod_{j\in\mathcal{V}-\{i\}}
\left[  1-\hat{p}_{j}(x,s_{j})\right]  ) dx,
\end{equation}
with $s_{i}^{c} = \{s_{j}: j \in\mathcal{V}-\{i\} \}$. 
Thus, the $H_{i}(\bar{s}_{i})$ term only depends on the neighborhood state $\bar{s}_i$, which is locally available at agent $i$ under the assumed neighbor information sharing paradigm. Therefore, $H_{i}(\bar{s}_{i})$ is called the \textit{local objective} function of agent $i$. On the other hand, $H_{i}^{c}(s_{i}^{c})$ is independent of $s_{i}$.

As a result, the local gradient of agent $i$, defined as
$d_{i}=\frac{\partial H_{i}(\bar{s}_{i})}{\partial s_{i}}\in\mathbb{R}^{2}$,
is always equal to the global gradient component $\frac{\partial
H(\mathbf{s})}{\partial s_{i}}$. Therefore, each agent $i$ can evaluate its
global gradient component using only its own local objective function
$H_{i}(\cdot)$ and the neighborhood state $\bar{s}_{i}$. As a result, the
distributed gradient ascent scheme in \eqref{Eq:GradientAscentGeneral} (i.e., $s_{i,k+1} = s_{i,k} + \beta_{i,k}d_{i,k}$) can be used to solve the problem in \eqref{Eq:OptimizationProblemCoverageControl} in a distributed manner.

In order to execute \eqref{Eq:GradientAscentGeneral}, each agent must properly evaluate its local gradient $d_{i,k}$ and select its step size $\beta_{i,k}$. 
The next Section \ref{SubSubSec:GradientDerivation} provides the derivation of $d_{i,k}$ and analyzes its structure which is pivotal in effective designing of boosted gradients for the use of boosting functions approach. The step size selection scheme is then presented in Section \ref{SubSubSec:StepSizeConstantDerivation}.

\subsubsection{Derivation of the gradient $d_{i}$}
\label{SubSubSec:GradientDerivation} 

Observing that the gradient $d_{i}$ is a
two dimensional vector, we write $d_{i} = [d_{iX},d_{iY}]^{T}$ and use
the Leibniz's rule \cite{Flanders1973} in \eqref{Eq:GlobalObjectiveFunction} to express $d_{iX}$ as
\begin{equation}
\label{Eq:GradientXGeneral}
\begin{aligned} 
d_{iX} = \frac{\partial H_i(\bar{s}_i)}{\partial s_{iX}} 
&= \frac{\partial }{\partial s_{iX}} \int_{V_i}R(x)\Phi_i(x)p_i(x,s_i)dx\\ 
&= \int_{V_i}R(x)\Phi_i(x)\frac{\partial p_i(x,s_i)}{\partial s_{iX}}dx\\
&+ \int_{\partial V_i}R(x)\Phi_i(x)p_i(x,s_i)\underline{V_{x}}\cdot\underline{n_{x}} dl, \end{aligned}
\end{equation}
where, 
\begin{equation}\label{Eq:Phi}
\Phi_{i}(x) = \prod_{j \in B_{i}}\left[  1-\hat{p}_{j}(x,s_{j})\right].    
\end{equation}
The second term in \eqref{Eq:GradientXGeneral} is a
line integral over the boundary of the sensing region $\partial V_{i}$. The terms $\underline{V_{x}}$ and $\underline{n_{x}}$ stand, respectively, for the rate of change and the unit normal vector of $\partial V_{i}$ at $x$ due to an infinitesimal change in $s_{iX}$, where $s_i = [s_{iX},s_{iY}]^T$.

looking at Fig. \ref{Fig:MissionSpaceGeometry}, observe that the shape of a
boundary $\partial V_{i}$ is formed by: (\romannum{1}) Mission space edges, (\romannum{2}) Obstacle edges, (\romannum{3}) Obstacle vertices, and, (\romannum{4}) Sensing range. However, when $s_{iX}$ (or $s_{iY}$) is perturbed infinitesimally, $\underline{V_{x}} \neq 0$ only
when $x$ lies on $\partial V_{i}$ components formed due to latter two factors. Therefore, we label the linear segments of $\partial V_{i}$ formed due to obstacle vertices as $\Gamma_{i} = \{\Gamma_{i1},\Gamma_{i2},\ldots\}$ and the circulary shaped curves formed due to finite sensing range as $\Theta_{i} =
\{\Theta_{i1},\Theta_{i2},\ldots\}$.

Using the fact that the sensing function $p_{i}(x,s_{i})$ depends on
$D_{i}(x) (= \Vert x-s_{i} \Vert)$, the first term in
\eqref{Eq:GradientXGeneral} can be simplified. Further, considering the
behavior of $\underline{V_{x}}\cdot\underline{n_{x}}$ on the segments in
$\Gamma_{i}$ and $\Theta_{i}$ sets, the line integral part of
\eqref{Eq:GradientXGeneral} can also be simplified to get two additional terms
(one term for linear segments $\Gamma_i$ and the other one for circular segments $\Theta_i$ of $\partial V_{i}$). Omitting some details, the complete expression for $d_{iX}$ is
\begin{equation}
\label{Eq:CompleteGradientX}%
\begin{split}
d_{iX}  &  = \int_{V_{i}}w_{i1}(x,\bar{s}_{i})\frac{(x-s_{i})_{X}}{\Vert
x-s_{i} \Vert}dx\\
&  + \sum_{\Gamma_{ij} \in\Gamma_{i}} sgn(n_{ijX})\frac{\sin{\theta_{ij}}%
}{\Vert v_{ij}-s_{i} \Vert} \int_{0}^{Z_{ij}}w_{i2}(\rho_{ir}(r),\bar{s_{i}%
})rdr\\
&  + \sum_{\Theta_{ij} \in\Theta_{i}} \delta_{i}\cos{\theta}\int_{\theta
_{ij1}}^{\theta_{ij2}}w_{i3}(\rho_{i\theta}(\theta),\bar{s_{i}})d\theta,
\end{split}
\end{equation}
where, $sgn(\cdot)$ is the signum function, and we define:
\begin{align}
w_{i1}(x,\bar{s}_{i}) =  &  -R(x)\Phi_{i}(x)\frac{dp_{i}(x,s_{i})}{dD_{i}%
(x)},\\
w_{i2}(x,\bar{s}_{i}) = w_{i3}(x,\bar{s}_{i}) =  &  R(x)\Phi_{i}%
(x)p_{i}(x,s_{i}),
\end{align}
with
\begin{align*}
\rho_{ir}(r) =  &  \frac{v_{ij}-s_{i}}{||v_{ij}-s_{i}||}r+v_{ij},\\
\rho_{i\theta}(\theta) =  &  s_{i} + \delta_{i}
\begin{bmatrix}
\cos{\theta} & \sin{\theta}%
\end{bmatrix}
^{T}.
\end{align*}

In order to uniquely quantify a line segment $\Gamma_{ij} \in
\Gamma_{i}$, it should contain the following geometric parameters \cite{Zhong2011}: end point $Z_{ij}$, angle $\theta_{ij}$, obstacle vertex $v_{ij}$, and direction $n_{ij}=[n_{ijX}, n_{ijY}]^T$. Thus, each $\Gamma_{ij}$ is a 4-tuple $(Z_{ij},\theta_{ij},v_{ij},n_{ij})$. Similarly, a circular arc segment $\Theta_{ij}$ is quantified by starting angle $\theta_{ij1}$ and ending angle $\theta_{ij2}$. Therefore, each $\Theta_{ij}$ term is a pair $(\theta_{ij1},\theta_{ij2})$ (see Fig. \ref{Fig:MissionSpaceGeometry}).

The complete expression in \eqref{Eq:CompleteGradientX} can be understood as a
sum of forces acting on agent $i$ (located in $s_{i}$), generated by different
points $x \in V_{i}$. The weight function $w_{i1}(x,\bar{s}_{i})$ in the first
term represents the magnitude of the force pulling agent $i$ towards point $x \in V_{i}$. The weight function $w_{i2}(x,\bar{s}_{i})$
in the second term describes the force generated in the lateral direction to the
line $\Gamma_{ij}$ (inwards the region $V_{i}$) by a point $x\in \Gamma_{ij}$.
Similarly, the weight function $w_{i3}(x,\bar{s}_{i})$ represents the
magnitude of the attraction force generated by a point $x \in \Theta_{ij}$.

Therefore, the gradient component $d_{iX}$ can be thought of as a function of
three weight functions: $d_{iX}=d_{iX}(w_{i1},w_{i2},w_{i3})$. This
representation will be used in the construction of boosting functions (specifically in constructing an expression for boosted gradients) in
Section \ref{SubSubSec:CovConConstructionOfBoostedGradientExpression}.

In contrast to previous work \cite{Zhong2011,Sun2014}, the effect of a limited sensing range is now incorporated into the gradient derivation resulting in the third term in \eqref{Eq:CompleteGradientX}. This modification is essential when $p_{i}(x,s_{i})$ does not approach its zero lower bound as $D_{i}(x)\rightarrow\delta_{i}$. This term is critically exploited in the
construction of the new family of boosting functions named \textquotedblleft
Arc-Boosting\textquotedblright\ as described in subsection
\ref{SubSubSec:DifferentBoostingFunctionFamilies} which, as we will see,
exhibits the best possible performance in terms of escaping local optima.

By following the same procedure, $d_{iY}$ can be derived as
\begin{equation}
\label{Eq:CompleteGradientY}%
\begin{split}
d_{iY}  &  = \int_{V_{i}}w_{i1}(x,\bar{s}_{i})\frac{(x-s_{i})_{Y}}{\Vert
x-s_{i} \Vert}dx\\
&  + \sum_{\Gamma_{ij} \in\Gamma_{i}} sgn(n_{ijY})\frac{\cos{\theta_{ij}}%
}{\Vert v_{ij}-s_{i} \Vert} \int_{0}^{Z_{ij}}w_{i2}(\rho_{ir}(r),\bar{s_{i}%
})rdr\\
&  + \sum_{\Theta_{ij} \in\Theta_{i}} \delta_{i}\sin{\theta}\int_{\theta
_{ij1}}^{\theta_{ij2}}w_{i3}(\rho_{i\theta}(\theta),\bar{s_{i}})d\theta.
\end{split}
\end{equation}
Now, using \eqref{Eq:CompleteGradientX} and \eqref{Eq:CompleteGradientY} each agent can locally evaluate its gradient ascent direction $d_{i,k}$ at each update instant $k$, as required in \eqref{Eq:GradientAscentGeneral}.

\subsubsection{Derivation of Step Size $\beta_{i,k}$}
\label{SubSubSec:StepSizeConstantDerivation} 

The step size selection mechanism proposed in this section is only required when the boosting function approach is not used. However, we present a method for computing the Lipshitz constant $K_{1i}$ of the $\nabla H_{i}(\bar{s}_{i})$, which will be an integral part of the optimal step size selection mechanism discussed in Section \ref{Sec:VariableStepSizeScheme} - when those concepts are applied.

As shown in \cite{Bertsekas2016nonlinear}, when an objective function
$H(\mathbf{s})$ is assumed to have a globally Lipschitz continuous gradient
with associated Lipschitz constant $K_{1}$, a state update law $\mathbf{s}%
_{k+1}=\mathbf{s}_{k}+\beta\nabla H(\mathbf{s}_{k})$ allows $\mathbf{s}_{k}$
to converge to a stationary state $\mathbf{s}^{\ast}$ (i.e. $\lim
_{k\rightarrow\infty}\mathbf{s}_{k}=\mathbf{s}^{\ast}$) when the step size
$\beta$ is chosen such that $\beta\in(0,\frac{2}{K_{1}}]$. 

This result cannot be directly applied to the distributed coverage control problem based on \eqref{Eq:OptimizationProblemCoverageControl} and \eqref{Eq:GradientAscentGeneral} due to two reasons: (\romannum{1}) The gradient of $H(\mathbf{s})$ in \eqref{Eq:GlobalObjectiveFunction}, $\nabla H(\mathbf{s})$, is only locally Lipschitz continuous (with corresponding Lipschitz constant $K_{1}=K_{1}(\mathbf{s})$), (\romannum{2}) Both the evaluation of $K_{1}(\mathbf{s})$ and communicating it across the agent network $\mathcal{G}$ prevents the decentralization of the gradient method in \eqref{Eq:GradientAscentGeneral}.

As a remedy, the step size $\beta_{i}$ in \eqref{Eq:GradientAscentGeneral} is
chosen such that $\beta_{i}\in(0,\frac{2}{k_{1i}(\bar{s}_{i})}]$ where
$K_{1i}(\bar{s}_{i})$ is the Lipschitz constant of $\nabla H_{i}(\bar{s}_{i}%
)$, to guarantee the convergence of \eqref{Eq:GradientAscentGeneral}. Using the
formal definition of the Lipschitz constant, an estimate for $K_{1i}(\bar
{s}_{i})$ can be computed locally at each agent $i$, using only the knowledge
of $\bar{s}_{i}$, through
\begin{equation}%
\begin{split}
\Vert\nabla(\nabla H_{i}(\bar{s}_{i}))\Vert_{\infty} &  =\max\{\sum_{j\in
\bar{B}_{i}}(|d_{j_{X}i_{X}i}|+|d_{j_{Y}i_{X}i}|)\ ,\\
\sum_{j\in\bar{B}_{i}}(|d_{j_{X}i_{Y}i}|+|d_{j_{Y}i_{X}i}|)\} &  \leq
K_{1i}(\bar{s}_{i}).
\end{split}
\label{Eq:LocalLipschitzConstant}%
\end{equation}
where, a general term $d_{j_{Y}i_{X}i}$ takes the form
\begin{equation}
\label{Eq:LocalLipschitzConstantGeneralTerm}%
\begin{split}
d_{j_{Y}i_{X}i} = \frac{\partial^{2} H_{i}(\bar{s}_{i})}{\partial s_{jY}
\partial s_{iX}}= \int_{V_{i} \cap V_{j}} R(x) \prod_{\substack{k \in B_{i}\\k
\neq j}}\left[  1-p_{k}(x,s_{k})\right] \\
\cdot\frac{dp_{i}(x,s_{i})}{dD_{i}(x)}\frac{(x-s_{i})_{Y}}{\Vert x-s_{i}
\Vert} \cdot\frac{-dp_{j}(x,s_{j})}{dD_{j}(x)}\frac{(x-s_{j})_{X}}{\Vert
x-s_{j} \Vert} dx.
\end{split}
\end{equation}
It can be proven that each $d_{j_{Y}i_{X}i}$ term
involved in \eqref{Eq:LocalLipschitzConstant} can be evaluated at agent $i$ using only $\bar{s}_{i}$. Therefore, this analysis yields an easier, accurate,
and distributed way to compute Lipschitz constants $\{K_{1i}:i\in
\mathcal{V}\}$ which will also be utilized later in Section
\ref{SubSec:CoverageControlDBSConvergence}.

\subsection{Designing boosting functions}
\label{SubSec:CoverageControlBoostingFunctions} 

As discussed in Section \ref{SubSec:CoverageControlProblemFormulation}, the coverage control objective function in \eqref{Eq:GlobalObjectiveFunction} is non-convex. Therefore, the gradient-based technique proposed in Section
\ref{SubSec:CoverageControlDistributedOptimization} (i.e., agents following \eqref{Eq:GradientAscentGeneral}) will always face the problem of
converging to a local optimum. As a means of escaping such a local optimum
(upon convergence to it) and search for a better local optimum solution, the boosting function approach is used. Therefore, during boosting sessions, agents have to use the boosted gradient $\hat{d}_{i,k}$ in \eqref{Eq:BoostedGradientAscent} (i.e., while in boosting mode). This subsection mainly focuses on constructing an appropriate expression for the boosted gradient $\hat{d}_{i,k}$ for the coverage control problem.

\subsubsection{Boosted gradient expression construction}
\label{SubSubSec:CovConConstructionOfBoostedGradientExpression}
When constructing a closed-form expression for the boosted gradient $\hat
{d}_{i,k}$, the key is to identify the components of the normal gradient
expression which control its \textit{direction} and \textit{magnitude}. In
analyzing \eqref{Eq:CompleteGradientX} we already saw that $d_{i}=d_{i}%
(w_{i1},w_{i2},w_{i3})$ where each weight function $w_{ij}=w_{ij}(x,\bar
{s}_{i})$ represents the magnitude component of each of three infinitesimal
forces, $j=1,2,3$, acting on agent $i$ generated at a point $x\in V_{i}$. In
addition, $\hat{d}_{i}$ should satisfy $\hat{d}_{i,k}\neq0$ whenever
$d_{i,k}=0$. Note that $d_{i,k}=0$ occurs when all the aforementioned
infinitesimal forces add up to a resultant force with zero magnitude (see also Section \ref{SubSubSec:TerminationConditionsForImplementations}). Avoiding
such equilibrium configurations, an expression for $\hat{d}_{i,k}$ can be
constructed by appropriately transforming the weight functions $\{w_{ij}%
(x,\bar{s}_{i}):j=1,2,3\}$. In this paper, we consider weight function
transformations given by
\begin{equation}
\hat{w}_{ij}(x,\bar{s}_{i})=\alpha_{ij}(x,\bar{s}_{i})w_{ij}(x,\bar{s}%
_{i})+\eta_{ij}(x,\bar{s}_{i}),\text{ \ \ }j=1,2,3.\label{Eq:WeightTransforms}%
\end{equation}
Here, both $\alpha_{ij},\eta_{ij}:\mathbb{R}^{2}\times\mathbb{R}^{2|\bar
{B}_{i}|}\rightarrow\mathbb{R}$ are known as \textit{transformation
functions}. Therefore, the boosted gradient $\hat{d}_{i,k}$ takes the form
\begin{equation}
\label{Eq:BoostedGradientGeneral}
\hat{d}_{i,k}=d_{i,k}(\hat{w}_{i1},\hat{w}_{i2},\hat{w}_{i3}). 
\end{equation}

In order to make sure that using the boosted gradient direction $\hat{d}%
_{i,k}$ is an \textquotedblleft intelligent\textquotedblright\ choice
(compared to just using a random direction), each agent $i$ should choose the
transformation functions $\alpha_{ij},\eta_{ij},$ $j=1,2,3$, so as to trigger
a systematic exploration of the mission space, as discussed next.

The expression for the boosted gradient $\hat{d}_{i,k}$ in terms of
transformation functions $\alpha_{ij}(x,\bar{s}_{i})$ and $\eta_{ij}(x,\bar{s}_{i})$ can be obtained by combining \eqref{Eq:CompleteGradientX} and \eqref{Eq:WeightTransforms}. The $X$ component of the resulting boosted gradient $\hat{d}_{iX}$ takes the
form
\begin{equation}
\label{Eq:CompleteGradientX2}%
\begin{split}
\hat{d}_{iX}  &  = \int_{V_{i}}\alpha_{i1}(x,\bar{s}_{i})w_{i1}(x,\bar{s}%
_{i})\frac{(x-s_{i})_{X}}{\Vert x-s_{i} \Vert}dx\\
&  + \int_{V_{i}}\eta_{i1}(x,\bar{s}_{i})\frac{(x-s_{i})_{X}}{\Vert x-s_{i}
\Vert}dx\\
&  + \sum_{\Gamma_{ij} \in\Gamma_{i}} sgn(n_{ijX})\frac{\sin{\theta_{ij}}%
}{\Vert v_{ij}-s_{i} \Vert} \int_{0}^{Z_{ij}}\alpha_{i2}(\rho_{ir}(r),\bar
{s}_{i})\cdot\\
&  \ w_{i2}(\rho_{ir}(r),\bar{s}_{i})rdr\\
&  + \sum_{\Gamma_{ij} \in\Gamma_{i}} sgn(n_{ijX})\frac{\sin{\theta_{ij}}%
}{\Vert v_{ij}-s_{i} \Vert} \int_{0}^{Z_{ij}}\eta_{i2}(\rho_{ir}(r),\bar
{s}_{i})rdr\\
&  + \sum_{\Theta_{ij} \in\Theta_{i}} \delta_{i}\cos{\theta}\int_{\theta
_{ij1}}^{\theta_{ij2}}\alpha_{i3}(x,\bar{s}_{i})w_{i3}(\rho_{i\theta}%
(\theta),\bar{s}_{i})d\theta\\
&  + \sum_{\Theta_{ij} \in\Theta_{i}} \delta_{i}\cos{\theta}\int_{\theta
_{ij1}}^{\theta_{ij2}}\eta_{i3}(\rho_{i\theta}(\theta),\bar{s}_{i})d\theta.
\end{split}
\end{equation}
Similarly, the expression for $\hat{d}_{iY}$ can be derived by transforming the expression for $d_{iY}$ using \eqref{Eq:WeightTransforms}.

\subsubsection{Boosting function families}
\label{SubSubSec:DifferentBoostingFunctionFamilies}

A \textit{boosting function family} is characterized by the form of the
transformation functions $\alpha_{ij}(x,\bar{s}_{i}),$ $\eta_{ij}(x,\bar
{s}_{i}),$ $j=1,2,3$. Therefore, different boosting function families exhibit
different properties. We will review three such boosting function families
proposed in \cite{Sun2014} and will introduce two new ones with properties
that specifically address the presence of obstacles (more generally,
constraints) in \eqref{Eq:OptimizationProblemCoverageControl}.

The underlying rationale behind constructing a boosting function lies in the
answer to the question: \textquotedblleft Once an agent converges under the
normal gradient-based mode, how can the agent escape the achieved equilibrium
towards a `meaningful' direction? \textquotedblright\ Here, a
`meaningful'\ direction choice is a one that encourages the agent to explore
the mission space giving a high priority to points which are likely to achieve
a higher objective function value than the current local optimum. 

To answer
this question in the context of coverage control, consider a situation where
an agent $i$ has converged to $s_{i}^{1}$ at update step $k=k_{1}$ after
following the normal mode. To define appropriate $\alpha_{ij}(x,\bar{s}_{i}),$
$\eta_{ij}(x,\bar{s}_{i}),$ $j=1,2,3$, in (\ref{Eq:WeightTransforms}), the
information available to agent $i$ consists of: (\romannum{1}) The
neighborhood state $\bar{s}_{i}$, (\romannum{2}) The local objective function
$H_{i}(\cdot)$, (\romannum{3}) The neighboring mission space topological
information contained in $\Gamma_{i}\ $and$\ \Theta_{i}$ (see Fig.
\ref{Fig:MissionSpaceGeometry}), (\romannum{4}) Past state trajectory
information $\{s_{i,k}:\ k<k_{1}\}$. 

In order to construct a boosting function
family, one or more of these forms of local information are used. The three
boosting function families proposed in \cite{Sun2014} use $\bar{s}_{i}$ and
$H_{i}(\cdot)$. In contrast, the new boosting function families proposed in
this paper make use of $\Gamma_{i},\Theta_{i}$ and $\{s_{i,k}:\ k<k_{1}\}$ in
addition to $\bar{s}_{i}$ and $H_{i}(\cdot)$.

In what follows, we refer to the setting where $\alpha_{ij}(x,\bar{s}_{i})=1,$
$\eta_{ij}(x,\bar{s}_{i})=0,$ $j=1,2,3$, as the \textit{default configuration}
in (\ref{Eq:WeightTransforms}). In defining boosting function families, we
will use $\kappa$ and $\gamma$ as two positive \emph{gain parameters}.

\paragraph*{\textbf{$\Phi$-Boosting} \cite{Sun2014}}

This method uses 
\begin{eqnarray}
\alpha_{i1}(x,\bar{s}_{i}) &=& \kappa\Phi_{i}(x)^{\gamma},\\
\eta_{i1}(x,\bar{s}_{i}) &=& 0,
\end{eqnarray}
where $\Phi_{i}(x)$ in \eqref{Eq:Phi} indicates
the extent to which point $x\in V_{i}$ is \emph{not} covered by neighbors in
$B_{i}$. Thus, the effect of $\Phi$-Boosting is to force agent $i$ to move
towards regions of $V_{i}$ which are less covered by its neighbors.

\paragraph*{\textbf{$P$-Boosting} \cite{Sun2014}}
In this method 
\begin{eqnarray}
\alpha_{i1}(x,\bar{s}_{i}) &=& \kappa \lbrack P(x,\mb{s}) \rbrack ^{-\gamma},\\
\eta_{i1}(x,\bar{s}_{i}) &=& 0
\end{eqnarray}
are used, where $P(x,\mb{s})$ in \eqref{Eq:P} indicates the extent to which point $x \in \Omega$ is covered by all the agents in $\mathcal{V}$. However, when evaluating the boosted gradient:  $x \in V_i \subseteq \Omega$. Therefore, $P$-Boosting assigns higher weights to points $x\in V_{i}$ which are less covered by the closed neighborhood $\bar{B}_{i}$.

\paragraph*{\textbf{Neighbor-Boosting} \cite{Sun2014}}

This boosting function family uses 
\begin{eqnarray}
\alpha_{i1}(x,\bar{s}_{i}) &=& 1,\\
\eta_{i1}(x,\bar{s}_{i})&=&\sum_{j\in B_{i}}1_{\{x=s_{j}\}}\cdot\frac
{\kappa\cdot1_{\{s_{j}\in V_{i}\}}}{\Vert s_{i}-x\Vert^{\gamma}},
\end{eqnarray}
where, $1_{\{\cdot\}}$ represents the indicator function. As a result
of this boosting method, agent $i$ gets repelled from its neighbors who are
also in its visibility region $V_{i}$.

Note that these boosting methods are limited to transforming the first
integral term of the gradient expression in \eqref{Eq:CompleteGradientX},
i.e., only the weight $w_{i1}(x,\bar{s}_{i})$ through $\alpha_{i1}(x,\bar
{s}_{i}),\ \eta_{i1}(x,\bar{s}_{i})$ is transformed (while $\alpha_{ij}%
(x,\bar{s}_{i}),\ \eta_{ij}(x,\bar{s}_{i}),\ j=2,3$, are set to their default
configuration). Next, we present two new boosting function families.

\paragraph*{\textbf{V-Boosting}}

The intuition behind the V-Boosting function family is to use the information
of obstacle vertices ($v_{ij}\in\Gamma_{ij}$) which lie inside $V_{i}$ so as
to aid agent $i$ to navigate around obstacles. Recall that the second integral
term in \eqref{Eq:CompleteGradientX} represents the effect of obstacles
$\Gamma_{i}$ in $V_{i}$ on agent $i$. Therefore, in V-Boosting, this second
integral term is modified by transforming $w_{i2}(x,\bar{s}_{i})$ via the
$\eta_{i2}(x,\bar{s}_{i})$ term, in addition to transforming $w_{i1}(x,\bar
{s}_{i})$.

Specifically, the V-Boosting function family uses
\begin{align}
\alpha_{i1}(x,\bar{s}_{i}) &  =\kappa_{1}\Phi_{i}(x)^{\gamma_{1}}%
(1-p_{i}(x,s_{i})),\label{Eq:VBoostingAlpha}\\
\eta_{i2}(x,\bar{s}_{i}) &  =1_{\{x=Z_{ij}\}}\cdot\kappa_{2}\Vert x-s_{i}%
\Vert^{\gamma_{2}}.\label{Eq:VBoostingEta}%
\end{align}

The transformation in (\ref{Eq:VBoostingAlpha}) forces agent $i$ to move
toward less covered areas while the transformation in (\ref{Eq:VBoostingEta})
acts as an attraction force directed towards $Z_{ij}\in\Gamma_{ij}$ (same as
in the direction of obstacle vertex $v_{ij}$). The combination of these two
influences facilitates agent $i$ to navigate around obstacles aiming to expand
the mission space exploration.


\paragraph*{\textbf{Arc-Boosting}}
This boosting function family is particularly effective when there are
multiple obstacles/constraints in the vicinity of agent $i$. Similar to the
way V-Boosting uses the information in $\Gamma_{i}$ to transform the weight
function $w_{i2}(x,\bar{s}_{i})$, in Arc-Boosting, the information in
$\Theta_{i}$ is utilized to transform the weight function $w_{i3}(x,\bar
{s}_{i})$.

Recall that $\{\theta_{ij1},\theta_{ij2}\}=\Theta_{ij}\in\Theta
_{i}$ represents a circular boundary segment (also called an \textit{arc}) due
to the finite nature of the sensing range. An agent can have multiple arcs in
its boundary set $\partial V_{i}$ depending on how the agent is located in the
mission space relative to obstacles. For example, for the agent in Fig.
\ref{Fig:MissionSpaceGeometry}, there are three such arcs. Under the
Arc-Boosting method, first, each arc segment $\Theta_{ij}\in\Theta_{i}$ is
classified into one of three disjoint sets: (\romannum{1}) Attractive Arcs
$\Theta_{i}^{+}$, (\romannum{2}) Repulsive Arcs $\Theta_{i}^{-}$, and
(\romannum{3}) Neutral Arcs $\Theta_{i}^{0}$. 

This classification is based on the metric $A(\Theta_{ij})$:
\begin{equation*}
A(\Theta_{ij})=\frac{1}{(\theta_{ij2}-\theta_{ij1})}\int_{\theta_{ij1}%
}^{\theta_{ij2}}(1-\prod_{k\in\bar{B}_{i}}(1-\hat{p}_{k}(\rho_{i\theta}%
(\theta),s_{k})))d\theta,
\end{equation*}
which measures the mean coverage level on the arc segment $\Theta_{ij}$ by the
agents forming the closed neighborhood $\bar{B}_{i}$. The arc with the maximum
$A(\Theta_{ij})$ value is assigned to be a repulsive arc (i.e., in the set
$\Theta_{i}^{+}$), while the arc with the minimum $A(\Theta_{ij})$ value is
assigned to be an attractive arc (i.e., in the set $\Theta_{i}^{-}$). The
remaining arcs are labeled as neutral (i.e., in the set $\Theta_{i}^{0}$).

However, it is possible that an equilibrium\ occurs (i.e., $A(\Theta_{ij})$
are identical for all $j$), which may happen when $B_{i}=\emptyset$. In this
case, we use a recent state $s_{i,k-K}$, where $K\geq1$ is a parameter of the
Arc-Boosting method, selected from the agent's past state trajectory.
Specifically, the arc which is in the direction of $s_{i,k-K}$ (from point
$s_{i}$) is regarded as a repulsive arc while all other arcs are labeled as attractive.

The arc partition consisting of sets $\Theta_{i}^{+},\Theta_{i}^{-}$ and
$\Theta_{i}^{0}$ is used to define the Arc-Boosting function family by
transforming the weight function $w_{i3}(x,\bar{s}_{i})$ using
\begin{align}
\alpha_{i3}(x,\bar{s}_{i})=  &  1_{\{\Theta_{ij}\in\Theta
_{i}^{0}\}},\label{Eq:alpha3}\\
\eta_{i3}(x,\bar{s}_{i})=  &  [1_{\{\Theta_{ij}\in\Theta_{i}^{+}%
\}}-1_{\{\Theta_{ij}\in\Theta_{i}^{-}\}}]\cdot F_{c}(\kappa,\gamma).\label{Eq:eta3}
\end{align}
In \eqref{Eq:eta3}, the value of the term in brackets is either $1,-1$ or $0$
depending on whether $\Theta_{ij}$ belongs to $\Theta_{i}^{+},\Theta_{i}^{-}$
or $\Theta_{i}^{0}$ respectively. The term $F_{c}(\kappa,\gamma)$ is a gain
factor which depends on the usual gain parameters $\kappa$ and $\gamma$ used
before; a typical choice is of the form $F_{c}(\kappa,\gamma)=\kappa
e^{\gamma}$.

The intuition behind this method is to encourage agent $i$ to: (\romannum{1})
Move away from repulsive arcs (i.e., from highly covered regions),
(\romannum{2}) Move towards attractive arcs (i.e., towards less covered
regions), and (\romannum{3}) Move continuously towards unexplored regions
(i.e., towards an opposing direction to the already visited point $s_{i,k-K}%
$). The Arc-Boosting family has been found to be the most effective in
handling the presence of multiple obstacles/constraints within $V_{i}$.



\subsection{Proposing distributed boosting scheme}
\label{SubSec:CoverageControlBoostingSchemes} 

In order to deploy any of the
discussed boosting function families, a \textit{boosting scheme}) is also required. As
discussed in Section \ref{Sec:GeneralProblemFormulation}, the objective of a
boosting scheme is to define when each agent should switch between normal and
boosting modes. A typical \textit{centralized boosting scheme} (CBS) first
proposed in the work \cite{Sun2014} is shown in Figure
\ref{Fig:CentralizedBoostingBlockDiagram} and was discussed in Section
\ref{Sec:GeneralProblemFormulation}. As an improvement, this subsection intends to introduce a novel \textit{distributed boosting scheme} (DBS).


As discussed in the Section \ref{Sec:GeneralProblemFormulation}, in contrast
to a CBS, agents governed by a DBS acts independently and therefore carries out mode
changes asynchronously. As a result, at a given time instant different agents can be in different modes. Further, agents under a DBS are deprived of the global
information (such as $\mathbf{s}$ and $H(\cdot)$), and, only allowed to use
the locally available information: $\bar{s}_{i}$ and $H_{i}(\cdot)$. However,
as pointed out in the Remark 
\ref{Rm:NeighborhoodObjectiveForDistributedBoosting}, in any boosting scheme,
it is important to have a technique to measure the effect of the boosting
stage on the global objective function. So, in the proposing DBS, as a means
of \textit{locally} tracking the effect of the boosting stage on the
\textit{global} coverage objective function $H(\mathbf{s})$, the
\textit{neighborhood coverage objective} function $\bar{H}_{i}(\bar{s}_{i})$
is used. It is given by
\begin{equation}
\label{Eq:NeighborhoodCoverageObjective}\bar{H}_{i}(\bar{\mathbf{s}_{i}}) =
\int_{\bar{V}_{i}}R(x)[1-\prod_{j\in\bar{B}_{i}}(1-\hat{p}_{j}(x,s_{j}))] dx,
\end{equation}
where, $\bar{V}_{i} = \cup_{j\in\bar{B}_{i}}V_{j}$.

\begin{remark}
\label{Rm:NeighborhoodObjectiveForCoverageControlDBS} Note the difference
between the construction of proposed neighborhood coverage objective $\bar{H}_{i}%
(\bar{\mathbf{s}_{i}})$ in \eqref{Eq:NeighborhoodCoverageObjective} and the
used generalized neighborhood objective $\tilde{H}_{i}(\tilde{s}_{i})$ in
\eqref{Eq:NeighborhoodCostFunctionTheory}. This is justifiable because each of
those functions are intended to serve two different purposes: $\bar{H}%
_{i}(\bar{\mathbf{s}_{i}})$ is used to locally assess the global effect of a
boosting stage, and, $\tilde{H}_{i}(\tilde{s}_{i})$ is just used as a stepping
stone in the convergence analysis. So, it is not necessary for $\bar{H}%
_{i}(\bar{s}_{i})$ and $\tilde{H}_{i}(\tilde{s}_{i})$ to have equivalent
expressions. However, as pointed out in the Remarks
\ref{Rm:WightedNeighborhoodObjective} and
\ref{Rm:NeighborhoodObjectiveForDistributedBoosting}, for certain cooperative
multi-agent optimization problems, these two functions can take an identical
form as well.
\end{remark}

To accurately describe the proposing DBS, a few clone modes to the
\textit{normal mode} are introduced. They are typically labelled as
`NM\textsuperscript{m}' and referred to as `\textit{Normal Mode-m}' where m
$\in$ \{0,1,2,3\}. Whenever an agent is in a such normal mode
NM\textsuperscript{m}, it follows the usual state update step given in
\eqref{Eq:GradientAscentGeneral}. Further, the \textit{boosting mode} is
labeled as BM, and, an agent in BM follows \eqref{Eq:BoostedGradientAscent}.
Furthermore, in order to represent the termination of the complete
optimization algorithm, a mode labeled as `FM' which is referred to as the
'Final Mode' is also introduced. Once an agent reaches the FM, it terminates
any further state updates. So, when all the agents have reached the FM, the
optimization algorithm is considered to be terminated.

\begin{figure*}[tbh]
\centering
\includegraphics[width=6.5in]{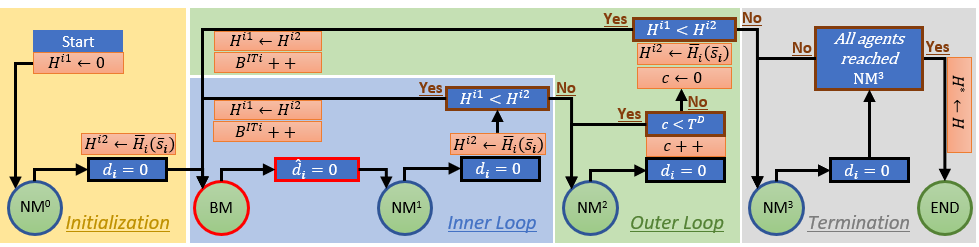}
\caption{The distributed boosting scheme (DBS) used by agent $i\in\{1,2,\ldots,N\}$ (see also Section \ref{SubSubSec:TerminationConditionsForImplementations})}%
\label{Fig:DecentralizedBoostingBlockDiagramCoverage}%
\vspace{-5mm}
\end{figure*}

The proposing novel DBS is described in the Algorithm
\ref{Alg:DecentralizedBoosting} and it is outlined in the Figure
\ref{Fig:DecentralizedBoostingBlockDiagramCoverage}. In this distributed setting, a
\textit{boosting iteration} of an agent is defined as the total time period
spent in the consecutive modes: BM, NM\textsuperscript{1} and/or
NM\textsuperscript{2} in a single run. The counting variable $B^{ITi}$ is used
to track the number of boosting iterations used by the agent $i$. In this DBS,
$H^{i1}$ and $H^{i2}$ are used as running variables to keep track of the
improvement occurred in the neighborhood coverage objective $\bar{H}_{i}%
(\bar{s}_{i})$ introduced in \eqref{Eq:NeighborhoodCoverageObjective} of the
agent $i$, during its boosting iteration. The corresponding neighborhood states of $H^{i1}$ and $H^{i2}$ are stored in $\bar{s}^{i1}$ and $\bar{s}^{i2}$ respectively.

\begin{algorithm}[!h]
\caption{The distributed boosting scheme (DBS) used by agent $i \in \{1,2,\ldots,N\} $}
\label{Alg:DecentralizedBoosting}
For each agent $i$ follow:
\begin{algorithmic}[1]
\State $B^{ITi}\leftarrow0$,\ $H^{i1}\leftarrow0$,\ $s_i\leftarrow[0,0]^T$, and start agent $i$ in NM\tsup{0}.
\State Wait till $d_{i}=0$ at some $\bar{s}_i=\bar{s}^{i2}$, assign \ $H^{i2} \leftarrow  \bar{H}_i(\bar{s}^{i2})$.
\While{$H^{i1} < H^{i2}$}
\While{$H^{i1} < H^{i2}$}
\State $H^{i1} \leftarrow H^{i2}$,\ $B^{ITi}\leftarrow B^{ITi}+1$.
\State Switch to BM and wait till $\hat{d}_{i}=0;$
\State \multiline{Switch to NM\tsup{1} and wait till $d_{i}=0$ at some\\ $\bar{s}_i=\bar{s}^{i2}$, and assign\ $H^{i2} \leftarrow  \bar{H}_i(\bar{s}^{i2})$}
\EndWhile
\State \multiline{Switch to NM\tsup{2}. Wait for $T^D$ update steps while \\ $d_{i}=0$ is occurred.}
\State \multiline{Wait for next $d_{i}=0$ instant at some $\bar{s}_i=\bar{s}^{i2}$, and assign \ $H^{i2} \leftarrow  \bar{H}_i(\bar{s}^{i2})$.}
\EndWhile
\State Switch to NM\tsup{3}. Wait for others to reach NM\tsup{3}.
\end{algorithmic}
Once all the agents are in NM\tsup{3} and $d_{i}=0; \ \forall i$, Switch all agents to END. Then,
\textbf{Return} $\mb{s}^{*} \leftarrow \mb{s}$, $H^{*} \leftarrow H(\mb{s})$.
\end{algorithm}

\begin{remark}\label{Rm:StateResetStepLackingInDBS}
The CBS showed in Fig. \ref{Fig:CentralizedBoostingBlockDiagram} uses a global
state reset (see the $H^{*} \leftarrow H$ step in the Figure
\ref{Fig:CentralizedBoostingBlockDiagram} which accompanies $\mb{s}^\ast \leftarrow \mb{s}$) whenever it finds out that the
boosting iteration failed to increase the global objective. However, such a state reset step is not being used in the proposed DBS. This is mainly due to the distributed nature of the proposed boosting scheme: Even if an agent
detects that no neighborhood coverage cost improvement occurred during a
certain boosting iteration, resetting its state $s_{i}$ to the last known best position (i.e., $\bar{s}^{i1}$'s component corresponds to the agent $i$) does not make any sense as agent $i$ has no control over its neighbor states
$\{s_{j}; j\in B_{i}\}$. 
\end{remark}

\begin{remark}\label{Rm:ExtraLoopInDBS}
Compared to the CBS shown in Fig. \ref{Fig:CentralizedBoostingBlockDiagram}, the proposed DBS uses an extra \textit{outer loop} to keep agents from exiting the boosting iterations (see `Outer Loop' section in Fig. \ref{Fig:DecentralizedBoostingBlockDiagramCoverage} and lines 3,9 - 10 in Algorithm \ref{Alg:DecentralizedBoosting}). This extra loop comes into
action when an agent $i$ detects no improvement in $\bar{H}_{i}(\bar{s}_{i})$
after its BM and NM\textsuperscript{1} stages. Once an agent $i$ is in that
outer loop, it delays assessing the improvement in $\bar{H}_{i}(\bar{s}_{i})$
by a minimum $T^{D}$ number of state update steps. The need for this `delay'
stage can be justified by the following reasoning. Due to the distributed
nature of the proposed boosting scheme, whenever an agent $i$ goes through
its BM, it will indirectly cause its neighbors to go through a transient
period. Note that $d_{i}=0$ (which is the condition that leads to the
improvement assessment of $\bar{H}_{i}(\bar{s}_{i})$) can occur even when
neighbors in $B_{i}$ have transient states. So, measuring the improvement
occurred in $\bar{H}_{i}(\bar{s}_{i})$ when its variables $\{s_{j},j \in
B_{i}\}$ are in a transient situation does not yield an accurate assessment.
So, the added outer loop gives an extra time period for the agents in $B_{i}$
to settle down before agent $i$ evaluates the improvement occurred in $\bar
{H}_{i}(\bar{s}_{i})$ again at the end of NM\textsuperscript{2}.
\end{remark}

\begin{remark}\label{Rm:TerminationOfTheDBS}
In the proposed DBS, note the block \textquotedblleft All agents reached NM\tsup{3}\textquotedblright\ at the termination stage. The underlying objective of it is to make each agent stop from updating their local state (by changing their modes to the mode END) at once, after all agents have finished their boosting iterations and achieved $d_i = 0, \forall i \in \mathcal{V}$. Therefore, all the agents will transition into the END mode synchronously, when the last remaining agent transitioned into the NM\tsup{3} mode. Although this step appears to be a global (i.e., a centralized) step, it can be easily achieved locally via a widely popular distributed binary consensus algorithm \cite{Abdaoui2013}. In such a scheme, the binary variable which the multi-agent network $\mathcal{V}$ should come to a consensus is \textquotedblleft $Y_i = \{\mbox{(Agent }i\mbox{ in NM\tsup{3}) AND (} d_i = 0 \mbox{)}\}$.\textquotedblright
\end{remark}

\begin{remark}
\label{Rm:NeighborhoodChangeCoverageControl} Note that $\bar{H}_{i}(\bar
{s}_{i})$ is strictly dependent on the current neighbor set $B_{i}$ of the
agent $i$. In the coverage control application, the neighbor set $B_{i}$ of an
agent $i$ can sometimes change during a simulation (specifically during the
early transient stages). This change is a result of the sensing capabilities based
`Neighbor' definition used in the considered coverage control problem: $j \in
B_{i} \iff V_{i} \cap V_{j} \neq\emptyset$. In order to handle such
neighborhood changes during the DBS, an additional subroutine given in the
Algorithm \ref{Alg:DBSNeighborhoodChange} should be evaluated in parallel to
the main routine in Algorithm \ref{Alg:DecentralizedBoosting}. This special
subroutine mainly ensures that both $H^{i1}$ and $H^{i2}$ measures used in
Algorithm \ref{Alg:DecentralizedBoosting} are computed based on a fixed
neighborhood so that comparing them is valid. However, in a situation where
the concept of `neighbors' is defined based on a fixed set of communication
links \cite{Sun2016}, such an extra subroutine is not required.

\begin{algorithm}[!h]
\caption{DBS - Reacting to neighborhood $B_i$ changes}\label{Alg:DBSNeighborhoodChange}
While executing the Algorithm \ref{Alg:DecentralizedBoosting}, if for any agent $i$, the neighborhood $B_i$ changed at time step $k$, then follow:
\begin{algorithmic}[1]
\State Recall $H^{i1}$ from Algorithm \ref{Alg:DecentralizedBoosting}
\State Take $\bar{s}^{i2}\leftarrow \bar{s}_{i,k-1}$,\ $H^{i2} \leftarrow  \bar{H}_i(\bar{s}^{i2})$
\If{Agent $i$ in BM}
\If{$H^{i1}<H^{i2}$}
\State $\bar{s}^{i1}\leftarrow \bar{s}_{i,k}$,\ $H^{i1} \leftarrow \bar{H}_i(\bar{s}^{i1})$.
\State Remain in BM. \textbf{Return}.
\Else
\State $\bar{s}^{i1}\leftarrow \bar{s}_{i,k}$,\ $H^{i1} \leftarrow \bar{H}_i(\bar{s}^{i1})$.
\State Switch to NM\tsup{1}. \textbf{Return}.
\EndIf
\ElsIf{Agent $i$ is in NM\tsup{1} or NM\tsup{2}}
\State $\bar{s}^{i1}\leftarrow \bar{s}_{i,k}$,\ $H^{i1} \leftarrow \bar{H}_i(\bar{s}^{i1})$.
\State Remain in same mode (NM\tsup{1} or NM\tsup{2}). \textbf{Return}.
\EndIf
\end{algorithmic}
\end{algorithm}

\end{remark}

\subsection{Convergence of the DBS}\label{SubSec:CoverageControlDBSConvergence}

When all the agents have reached the mode END, the DBS is considered to be
converged. However, to reach the END, two conditions should be satisfied: 1)
all the agents first need to reach the NM\textsuperscript{3}, and then, 2)
they should achieve $d_{i}=0;\ \forall i$. In order to guarantee the latter
condition, it is required to ensure each agent $i$ has the capability to converge locally (i.e., $\lim_{k\rightarrow\infty} d_{i,k}=0$) when all of its
neighbors are in a normal mode (as in NM\textsuperscript{3}). However, to
guarantee the first condition, it is required to ensure that any agent can
reach the NM\textsuperscript{3} irrespective of its neighbors. This is because
modes NM\textsuperscript{3} and FM are absorbing with respect to other modes.


For a fixed neighborhood $B_{i}$ of an agent $i$, the neighborhood cost
function $\bar{H}_{i}(\bar{s}_{i})$ is a non negative function with a finite
upper-bound. Thus, boosting iterations cannot improve $\bar{H}_{i}(\bar{s}%
_{i})$ indefinitely. As a consequence, an agent $i$ is guaranteed to reach
NM\textsuperscript{3} if it can always escape modes: 1) BM by reaching
$\hat{\mathbf{d}}_{i}=0$ and 2) NM\textsuperscript{0},NM\textsuperscript{1} or
NM\textsuperscript{2} by reaching $\mathbf{d}_{i}=0$, irrespective of the
modes of its neighbors. In essence, to guarantee the convergence of the
proposed DBS, it is required to establish the same convergence criteria given
in \eqref{Eq:RequiredConvergence1}-\eqref{Eq:RequiredConvergence3}, where
$\mathcal{B}$ stands for the set of agents who are in the boosting mode and
$\mathcal{N}$ stands for the set of agents who are in a normal mode.

The information presented so far in this Section
\ref{Sec:CoverageControlApplication} confirms the fact that coverage control
problem falls directly under the general class of cooperative multi-agent
optimization problems discussed in Section \ref{Sec:GeneralProblemFormulation}. 
As a result, the developed general variable step size scheme presented in
Section \ref{Sec:VariableStepSizeScheme} can be considered as an available
avenue for guaranteeing the convergence of the proposed DBS. However, in
order to use this specific variable step size scheme (i.e. the step sizes
given by Theorem \ref{Th:ConvergenceHat}), coverage control problems should
satisfy the underlying assumptions of Theorem \ref{Th:ConvergenceHat}:
Assumptions \ref{As:PropertiesofH_i},\ref{As:LocalAvailabilityOfd_ij}%
,\ref{As:PositivityOfQTilde} and \ref{As:PsiHat}.

The Assumption \ref{As:PropertiesofH_i} holds for the coverage control problem
due to two reasons: 1) Section \ref{SubSubSec:StepSizeConstantDerivation}
already discussed a methodology for computing the Lipshitz constant $K_{1i}$
of $\nabla H_{i}(\bar{s}_{i})$ - locally. From
\eqref{Eq:LocalLipschitzConstant} and
\eqref{Eq:LocalLipschitzConstantGeneralTerm} it is clear that whenever the
sensing capabilities are smooth (i.e. $p_{i}(x,s_{i})$ is differentiable w.r.t
$D_{i}(x)$) the computed $K_{1i}$ value will be always finite. 2) $H_{UB} =
\int_{V_{i}}R(x)dx$ is a typical upper bound for $H_{i}(\bar{s}_{i})$ as
$\int_{\Omega}R(x)dx<\infty$ is already enforced in subsection
\ref{SubSec:CoverageControlProblemFormulation}.

The Assumption \ref{As:LocalAvailabilityOfd_ij} holds for coverage control
problem because information sharing capability is already assumed in the basic
coverage control problem framework \cite{Zhong2011,Sun2014}. However, the following lemma is useful to convince that no additional communication
bandwidth is required to satisfy this assumption.

\begin{lemma}
\label{Lm:LocalAvailabilityd_ij} For the class of coverage control problems,
any agent $i\in\mathcal{V}$ can locally compute $d_{ij}=\frac{\partial
H_{j}(\bar{s}_{j})}{\partial s_{i}}$ value $\forall j\in\bar{B}_{i}$.
\end{lemma}

\emph{Proof: } By taking the partial derivative of
\eqref{Eq:LocalObjectiveCoverageControl} (written for agent $j$) w.r.t. the
local state $s_{i}$ yields
$$
d_{ij} = -\int_{V_{j}}R(x)p_{j}(x,s_{j})\prod_{l\in B_{j} -\{i\}}%
(1-p_{l}(x,s_{l}))\frac{dp_{i}(x,s_{i})}{ds_{i}}dx.
$$
Now, note that $\forall x \not \in V_{i},\ \frac{-dp_{i}(x,s_{i})}{ds_{i}}=0$,
and, $\forall l\not \in B_{i} \cap B_{j}, \forall x \in V_{i} \cap
V_{j},\ p_{l}(x,s_{l})=0$. By incorporating these relationships into the
obtained expression for $d_{ij}$ gives a locally computable (at agent $i$)
expression for $d_{ij}$ as
\begin{equation}
d_{ij} = -\int_{V_{i} \cap V_{j}}R(x)p_{j}(x,s_{j})\prod_{l\in B_{i} \cap
B_{j}}(1-p_{l}(x,s_{l}))\frac{dp_{i}(x,s_{i})}{ds_{i}}dx.
\end{equation}
\hfill $\blacksquare$

\begin{figure}[t]
\centering
\includegraphics[width=2.8in]{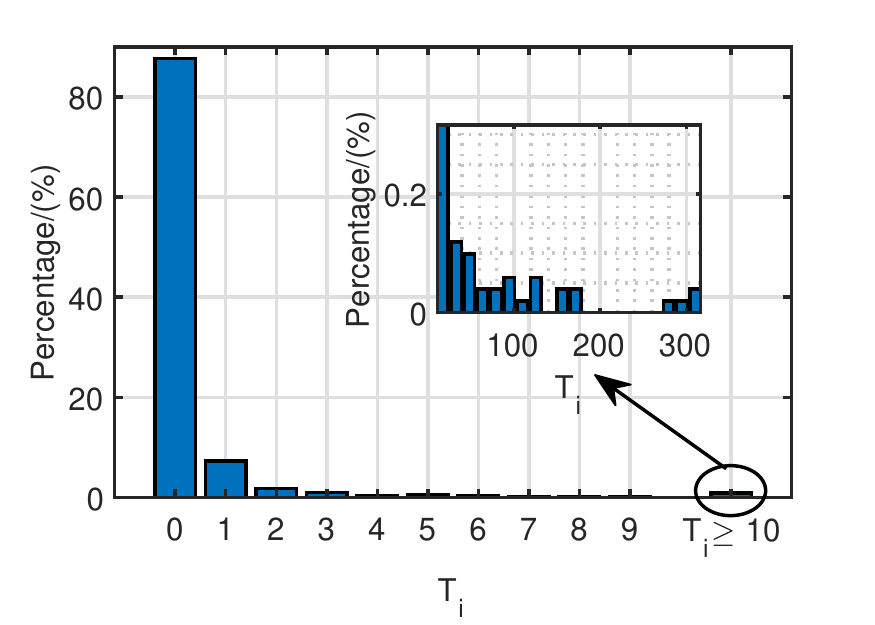}\caption{Percentage occurrence of
different $T_{i}$ values (Regarding Assumption \ref{As:PositivityOfQTilde} for
the simulation which produced the result shown in Fig. \ref{fig:MDA10}).}%
\label{Fig:PositivityOfQ_i}%
\end{figure}

The Assumption \ref{As:PositivityOfQTilde} has been previously justified for
general applications using Lemma \ref{Lm:QProperties1} and
\ref{Lm:QProperties2}. Further, to ensure that Assumption
\ref{As:PositivityOfQTilde} is satisfied by coverage control problem, the
parameter $T_{i}$ was observed during all the simulations (presented in
Section \ref{SubSec:SimulationResults}) for all the agents. In all occasions,
$T_{i}$ was found to be a finite value, implying that Assumption
\ref{As:PositivityOfQTilde} is valid. One such observed $T_{i}$ value
distribution is given in Fig. \ref{Fig:PositivityOfQ_i}, where $T_{i}$ lied
below 10 for 99.1\% of the time.

As justified in the Section \ref{Sec:VariableStepSizeScheme}, the Assumption
\ref{As:PsiHat} is trivial and it will hold for any general cooperative multi-agent problem including coverage control problems.

In conclusion, Assumptions \ref{As:PropertiesofH_i},
\ref{As:LocalAvailabilityOfd_ij}, \ref{As:PositivityOfQTilde} and
\ref{As:PsiHat} holds for the class of coverage control problems. Thus, the
variable step size scheme proposed in Theorem \ref{Th:ConvergenceHat} can be
used for the class of coverage control problems to ensure its convergence when
the proposed novel distributed boosting scheme is used.

\subsection{Simulation Results}
\label{SubSec:SimulationResults}

\begin{table}[b]
\caption{Boosting function parameters used in simulation results}%
\label{Tab:BoostingParamters}%
\centering
\begin{tabular}
[c]{|l|l|}\hline
\multicolumn{1}{|c|}{\textbf{Boosting Method}} &
\multicolumn{1}{c|}{\textbf{Associated Default Parameters}}\\\hline
$P$-Boosting & $\kappa=1,\ \gamma=1$\\\hline
Neighbor-Boosting & $\kappa=10000,\ \gamma=1$\\\hline
$\Phi$-Boosting & $\kappa=4,\ \gamma=2$\\\hline
V-Boosting & $\kappa_{1}=10,\ \kappa_{2}=5,\ \gamma_{1}=1,\ $and, $\gamma
_{2}=1$\\\hline
Arc-Boosting & $\kappa=1,\ \gamma=1,\ K=50,\ T_{D} = 5$\\\hline
\end{tabular}
\end{table}

\begin{figure}[t]
\centering
\includegraphics[width=3in]{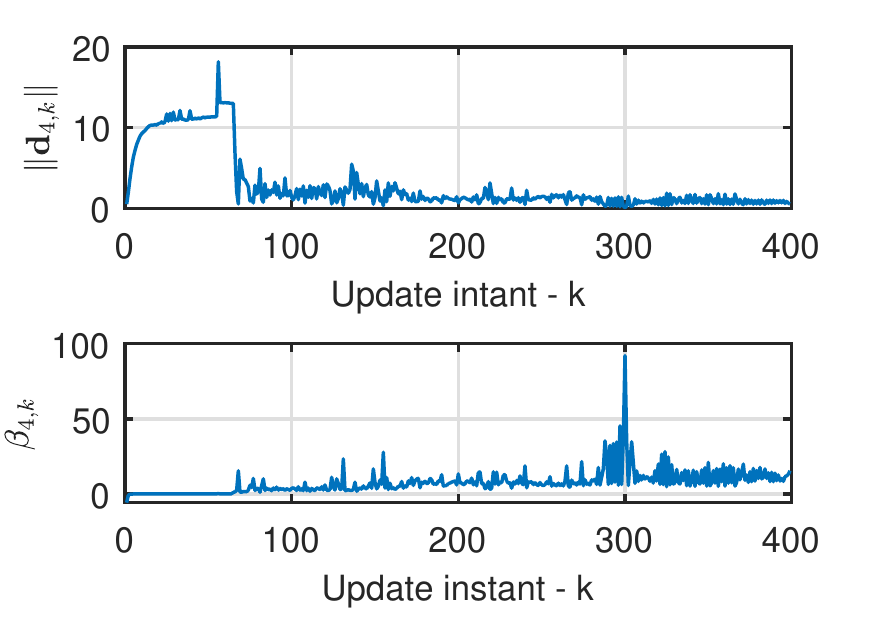}\caption{Variation of gradient
magnitude and the step size for the agent $i=4$ during the simulation which
yielded Fig. \ref{fig:NDV10}.}%
\label{Fig:StepSizeVariationExample}%
\end{figure}

\begin{figure}[!t]
\centering
\begin{subfigure}{0.22\columnwidth}
\includegraphics[width=\textwidth]{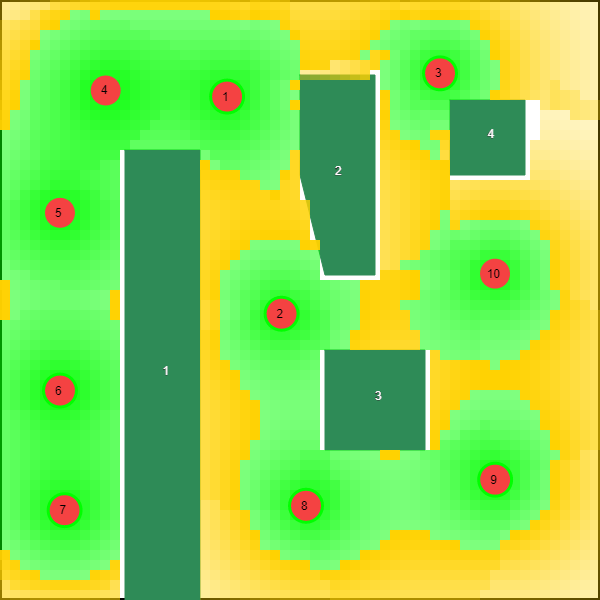}
\caption{\\GA:\ $158,821$}
\label{fig:GDZ10}
\end{subfigure}
\begin{subfigure}{0.22\columnwidth}
\includegraphics[width=\textwidth]{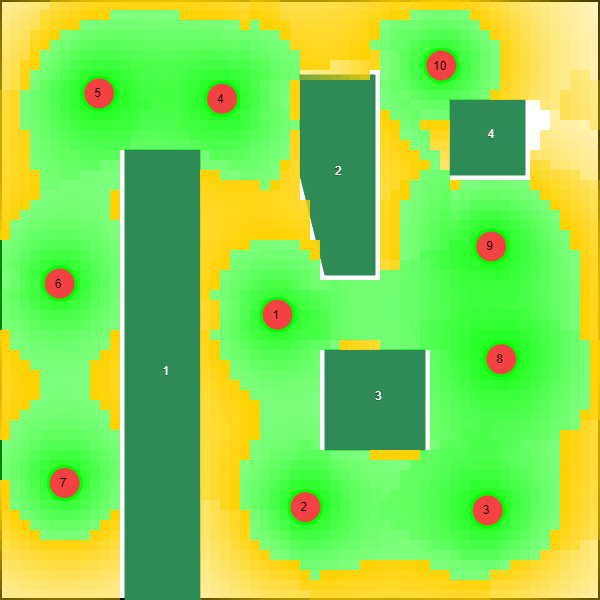}
\caption{$+ 2.31\%$\\$\Phi$B:\ $162,495$}
\label{fig:GDA10}
\end{subfigure}
~ \begin{subfigure}{0.22\columnwidth}
\includegraphics[width=\textwidth]{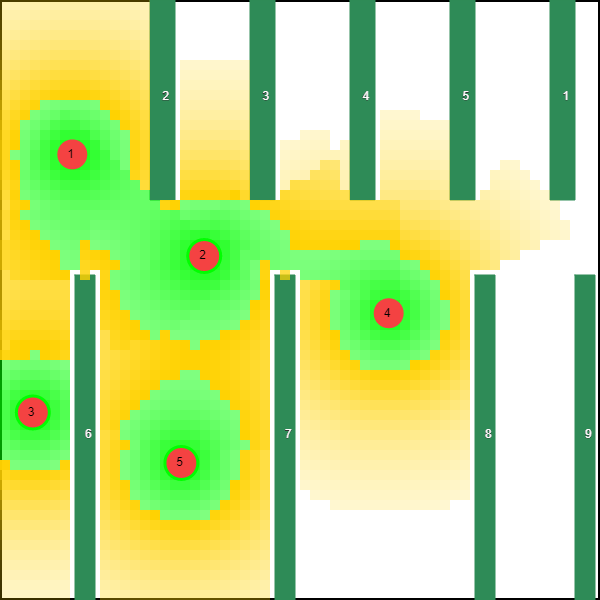}
\caption{\\GA:\ $86,638$}
\label{fig:RDZ10}
\end{subfigure}
\begin{subfigure}{0.22\columnwidth}
\includegraphics[width=\textwidth]{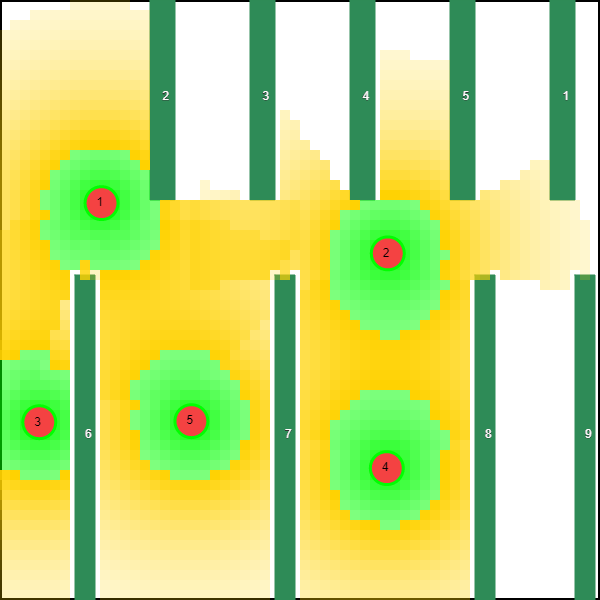}
\caption{$+ 2.89\%$\\$\Phi$B:\ $89,146$}
\label{fig:RDA10}
\end{subfigure}
~ \begin{subfigure}{0.22\columnwidth}
\includegraphics[width=\textwidth]{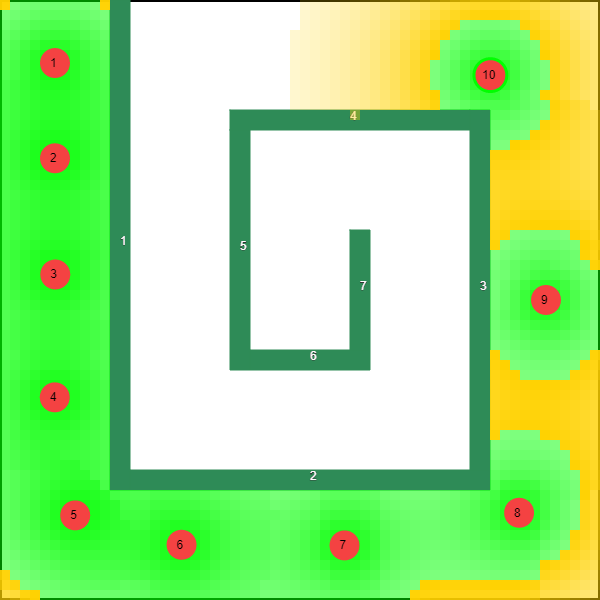}
\caption{\\GA:\ $120,343$}
\label{fig:MDZ10}
\end{subfigure}
\begin{subfigure}{0.22\columnwidth}
\includegraphics[width=\textwidth]{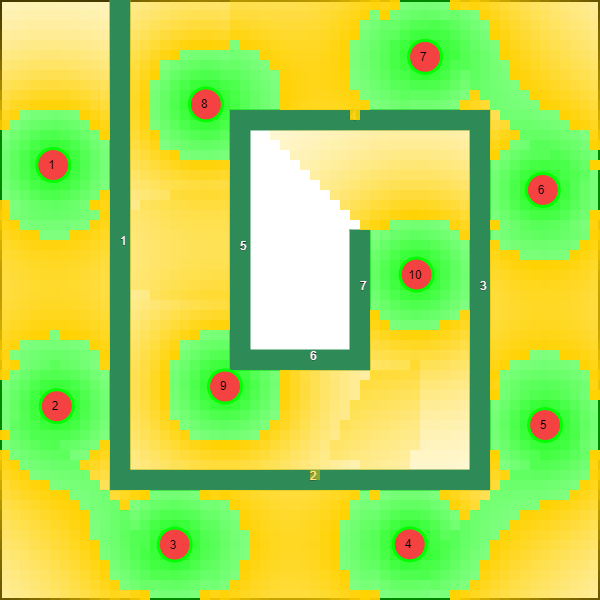}
\caption{$+ 19.8\%$\\$\Phi$B:\ $144,211$}
\label{fig:MDA10}
\end{subfigure}
~ \begin{subfigure}{0.22\columnwidth}
\includegraphics[width=\textwidth]{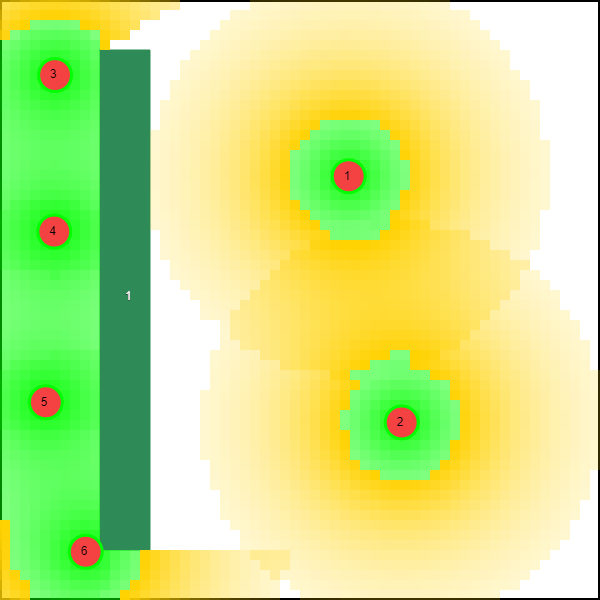}
\caption{\\GA:\ $101,976$}
\label{fig:NDZ10}
\end{subfigure}
\begin{subfigure}{0.22\columnwidth}
\includegraphics[width=\textwidth]{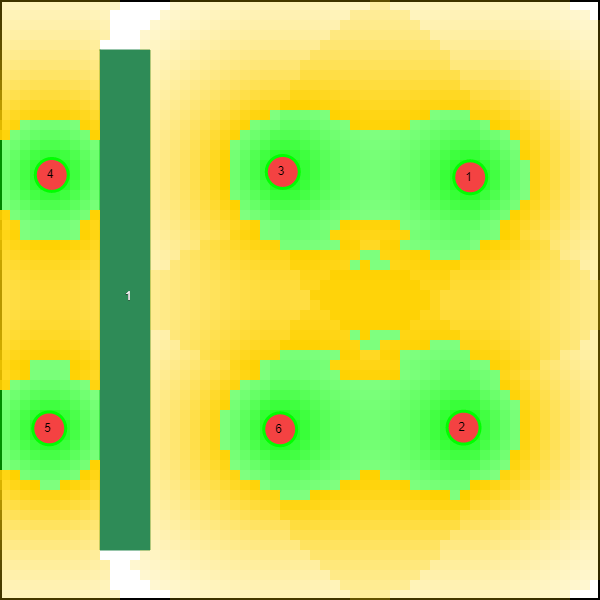}
\caption{$+ 27.0\%$\\$\Phi$B:\ $129,542$}
\label{fig:NDV10}
\end{subfigure}
\caption{Coverage improvement due to distributed $\Phi$-Boosting for
$N=10,5,6$ situations.}%
\label{Fig:Results0}%
\end{figure}

\begin{table}[!b]
\caption{Coverage objective value for cases with $N=1,2$ with decentralized
boosting}%
\label{Tab:Results2}
\centering
\begin{tabular}
[c]{|l|l|r|r|r|}\hline
\multicolumn{2}{|l|}{Configuration} &
\multicolumn{1}{c|}{\multirow{2}{*}{\begin{tabular}[c]{@{}c@{}}Gradient \\ Descent\end{tabular}}}
&
\multicolumn{1}{c|}{\multirow{2}{*}{\begin{tabular}[c]{@{}c@{}}Decentralized \\ V-Boosting\end{tabular}}}
&
\multicolumn{1}{c|}{\multirow{2}{*}{\begin{tabular}[c]{@{}c@{}}Decentralized \\ Arc-Boosting\end{tabular}}}%
\\\cline{1-2}%
Obstacles & N & \multicolumn{1}{c|}{} & \multicolumn{1}{c|}{} &
\multicolumn{1}{c|}{}\\\hline
General & 1 & \textbf{20,494} & 20,404 & \textbf{23,193}\\\hline
Maze & 1 & \textbf{14,759} & 14,774 & \textbf{17,090}\\\hline
Narrow & 1 & \textbf{13,669} & \textbf{30,259} & 30,178\\\hline
Narrow & 2 & \textbf{26,258} & \textbf{58,693} & 58,681\\\hline
\end{tabular}
\end{table}

\begin{table*}[t]
\caption{Coverage objective value \textbf{increment} (+/-) achieved by
different boosting schemes}%
\label{Tab:Results1}
\centering
\resizebox{\textwidth}{!}{\begin{tabular}{|l|l|r|r|r|r|r|r|r|r|r|r|r|r|r|}
\hline
\multicolumn{2}{|l|}{} & \multicolumn{1}{l|}{\textbf{\begin{tabular}[c]{@{}c@{}}Reference\\ Level $H(\mb{s}^1)$\end{tabular}}} & \multicolumn{12}{c|}{{\ul \textbf{Coverage objective value increment occurred with respect to the `Reference Level $H(\mb{s}^1)$'}}} \\ \hline
\multicolumn{2}{|l|}{Configuration} & \multicolumn{1}{c|}{\multirow{2}{*}{\begin{tabular}[c]{@{}c@{}}Gradient\\ Ascent (GA)\end{tabular}}} & \multicolumn{2}{c|}{Random Pert.} & \multicolumn{2}{c|}{$P$-Boosting} & \multicolumn{2}{c|}{Neighbor Boo.} & \multicolumn{2}{c|}{$\Phi$-Boosting ($\Phi$B)} & \multicolumn{2}{c|}{V-Boosting (VB)} & \multicolumn{2}{c|}{Arc-Boosting (AB)} \\ \cline{1-2} \cline{4-15}
Obstacles & N & \multicolumn{1}{c|}{} & \multicolumn{1}{l|}{Centr.} & \multicolumn{1}{l|}{Decen.} & \multicolumn{1}{l|}{Centr.} & \multicolumn{1}{l|}{Decen.} & \multicolumn{1}{l|}{Centr.} & \multicolumn{1}{l|}{Decen.} & \multicolumn{1}{l|}{Centr.} & \multicolumn{1}{l|}{Decen.} & \multicolumn{1}{l|}{Centr.} & \multicolumn{1}{l|}{Decen.} & \multicolumn{1}{l|}{Centr.} & \multicolumn{1}{l|}{Decen.} \\ \hline
General & 10 & \textbf{158,821} & +233 & +409 & +235 & +3684 & +235 & +3676 & +243 & +3674 & +2453 & +3621 & +3553 & \textbf{+3739} \\ \hline
Room & 10 & \textbf{143,583} & +1366 & +484 & +1578 & +2680 & +2374 & +968 & +1578 & +2626 & +1739 & +2455 & +1578 & \textbf{+2768} \\ \hline
Maze & 10 & \textbf{120,343} & +20037 & +19409 & +25937 & +25897 & +19443 & +25895 & +26952 & +23868 & +19970 & +25702 & +25945 & \textbf{+27142} \\ \hline
Narrow & 10 & \textbf{169,793} & +150 & +8781 & +9204 & +8835 & +15258 & +9391 & +15008 & +9376 & +14969 & \textbf{+15286} & +15238 & +15120 \\ \hline
\end{tabular}}\end{table*}

\begin{figure}[t]
\centering
\begin{subfigure}{0.22\columnwidth}
\includegraphics[width=\textwidth]{figures/GDZ10.png}
\caption{\ \\GA:\ $158,821$}
\label{fig:GDZ10}
\end{subfigure}
\begin{subfigure}{0.22\columnwidth}
\includegraphics[width=\textwidth]{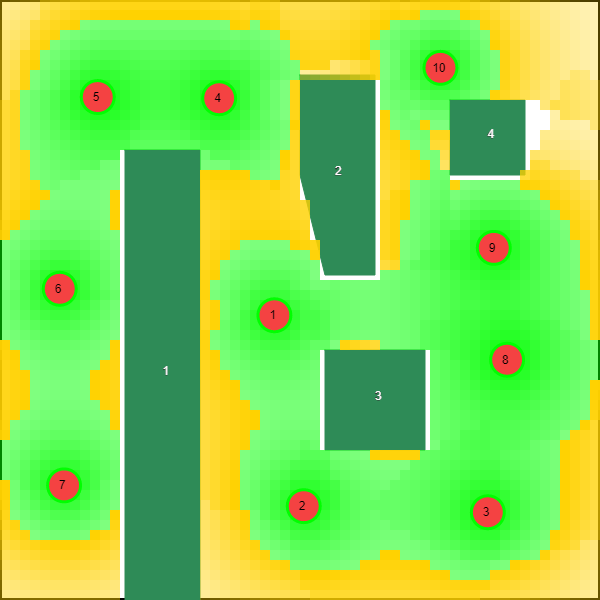}
\caption{$+ 2.354\%$\\AB:\ $162,560$}
\label{fig:GDA10}
\end{subfigure}
~ \begin{subfigure}{0.22\columnwidth}
\includegraphics[width=\textwidth]{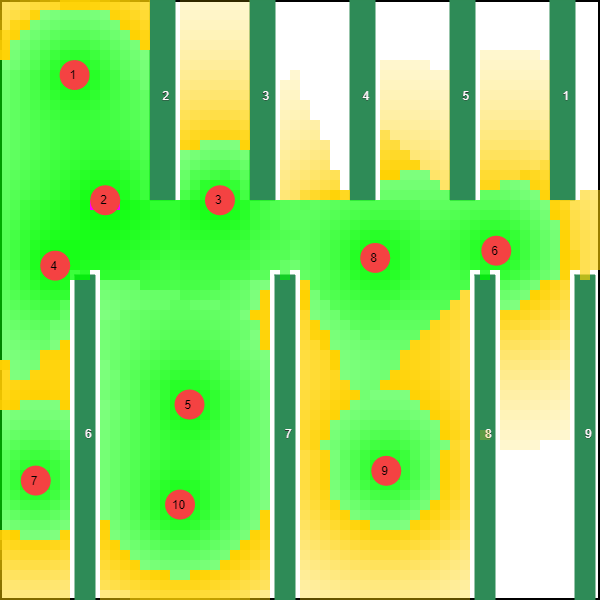}
\caption{\ \\GA:\ $143,583$}
\label{fig:RDZ10}
\end{subfigure}
\begin{subfigure}{0.22\columnwidth}
\includegraphics[width=\textwidth]{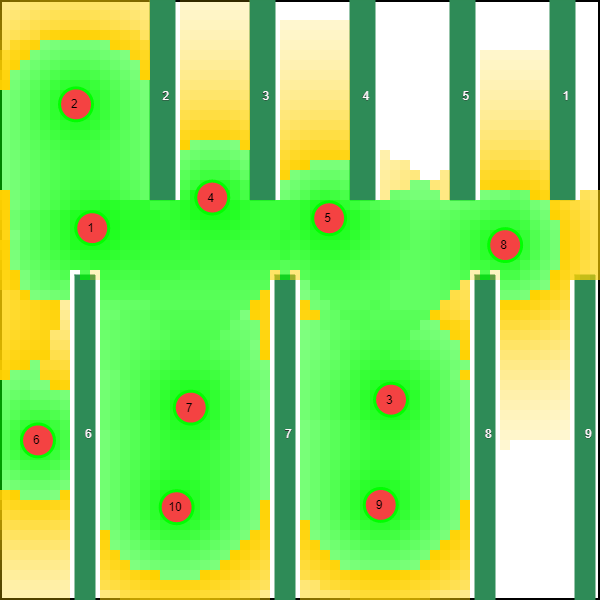}
\caption{$+ 1.928\%$\\AB:\ $146,351$}
\label{fig:RDA10}
\end{subfigure}
~ \begin{subfigure}{0.22\columnwidth}
\includegraphics[width=\textwidth]{figures/MDZ10.png}
\caption{\ \\GA:\ $120,343$}
\label{fig:MDZ10}
\end{subfigure}
\begin{subfigure}{0.22\columnwidth}
\includegraphics[width=\textwidth]{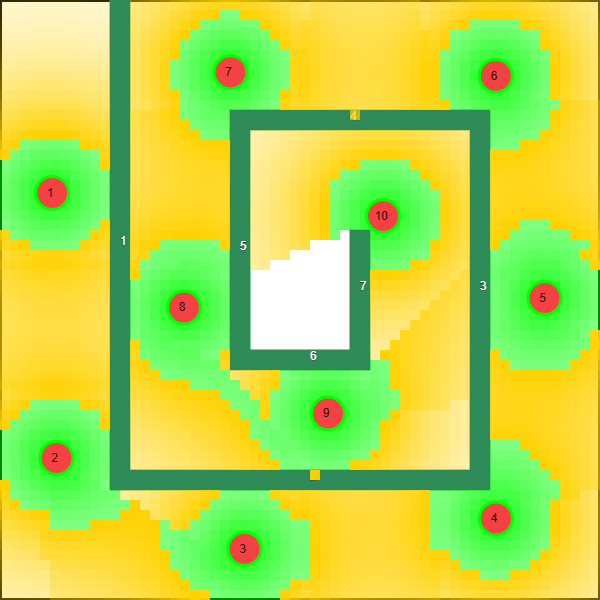}
\caption{$+ 22.55\%$\\AB:\ $147,485$}
\label{fig:MDA10}
\end{subfigure}
~ \begin{subfigure}{0.22\columnwidth}
\includegraphics[width=\textwidth]{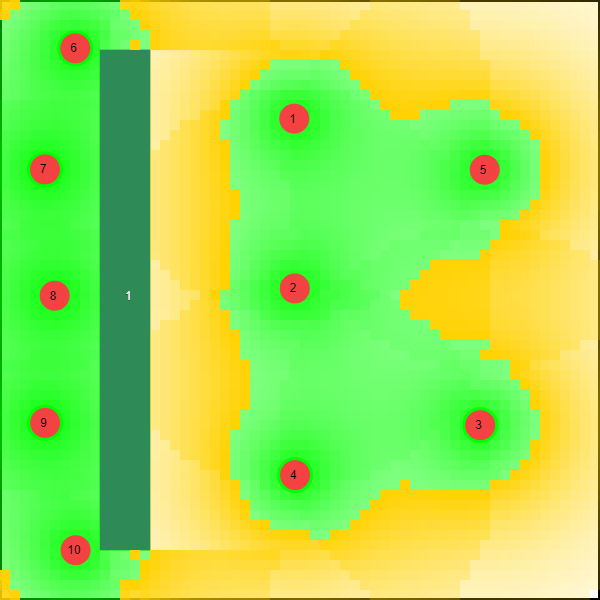}
\caption{\ \\GA:\ $169,793$}
\label{fig:NDZ10}
\end{subfigure}
\begin{subfigure}{0.22\columnwidth}
\includegraphics[width=\textwidth]{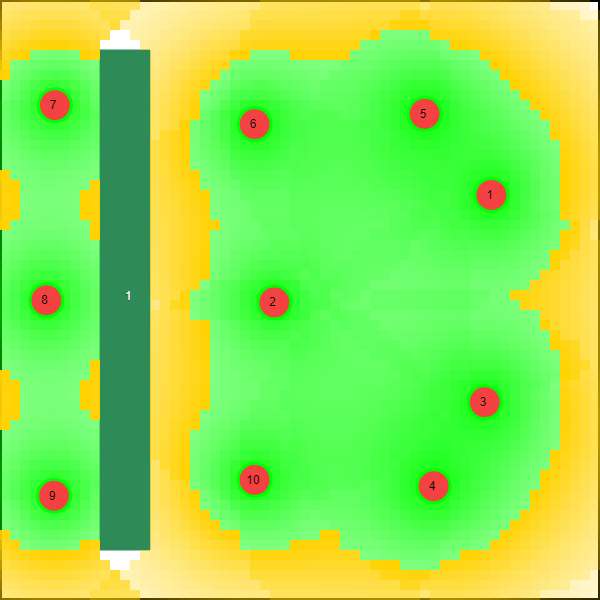}
\caption{$+ 9.003\%$\\VB:\ $185,079$}
\label{fig:NDV10}
\end{subfigure}
\caption{Maximum coverage improvement achieved due to boosting for $N=10$}%
\label{Fig:Results1}%
\vspace{-5mm}
\end{figure}

\begin{figure}[t]
\centering
\begin{subfigure}{0.22\columnwidth}
\includegraphics[width=\textwidth]{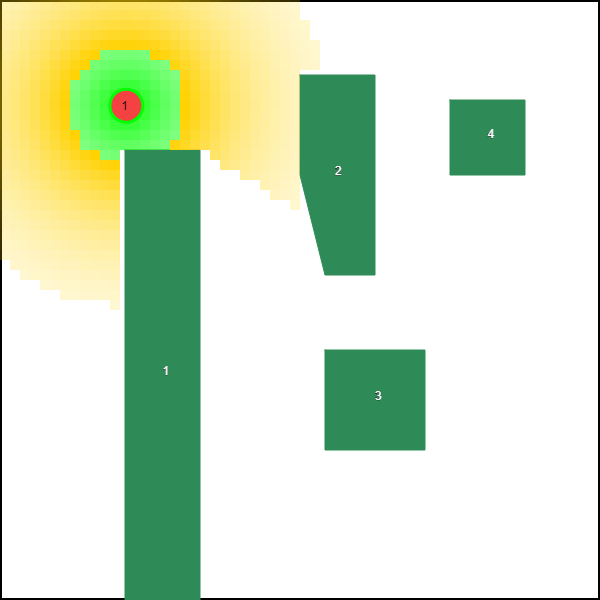}
\caption{\ \\GA:\ $20,494$}
\label{fig:GDZ1}
\end{subfigure}
\begin{subfigure}{0.22\columnwidth}
\includegraphics[width=\textwidth]{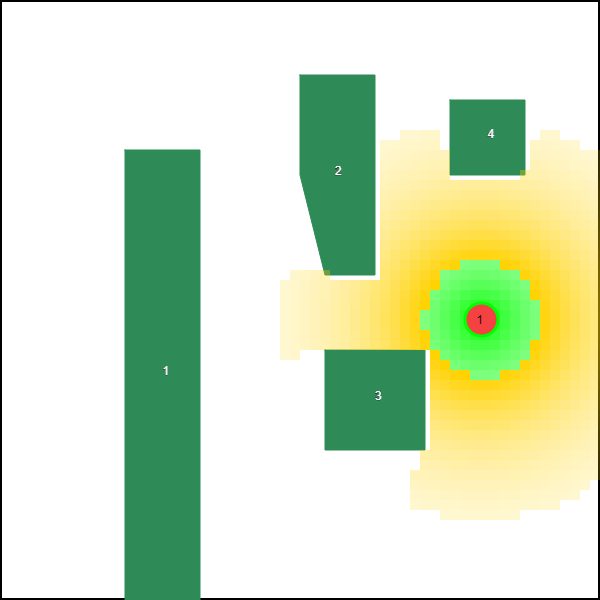}
\caption{$+ 13.17\%$\\AB:\ $23,193$}
\label{fig:GDA1}
\end{subfigure}
~ \begin{subfigure}{0.22\columnwidth}
\includegraphics[width=\textwidth]{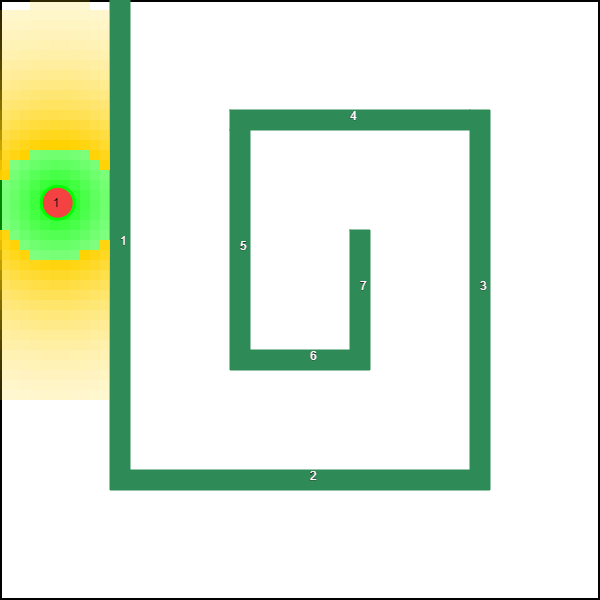}
\caption{\ \\GA:\ $147,59$}
\label{fig:MDZ1}
\end{subfigure}
\begin{subfigure}{0.22\columnwidth}
\includegraphics[width=\textwidth]{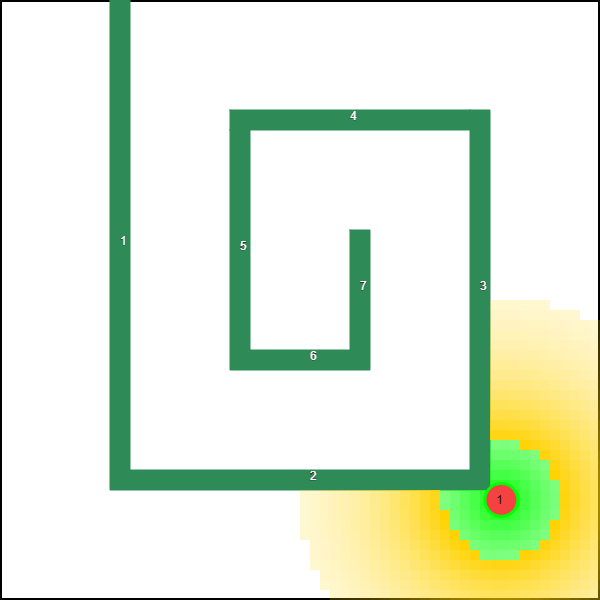}
\caption{$+ 15.79\%$\\AB:\ $17,090$}
\label{fig:MDA1}
\end{subfigure}
~ \begin{subfigure}{0.22\columnwidth}
\includegraphics[width=\textwidth]{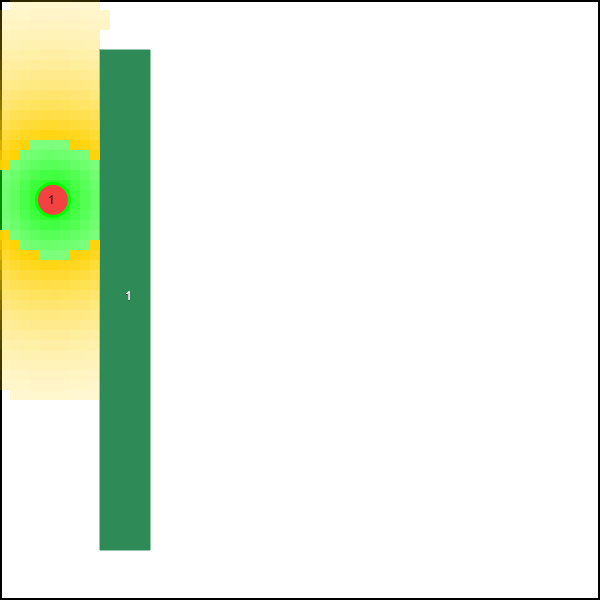}
\caption{\ \\GA:\ $13,669$}
\label{fig:NDZ1}
\end{subfigure}
\begin{subfigure}{0.22\columnwidth}
\includegraphics[width=\textwidth]{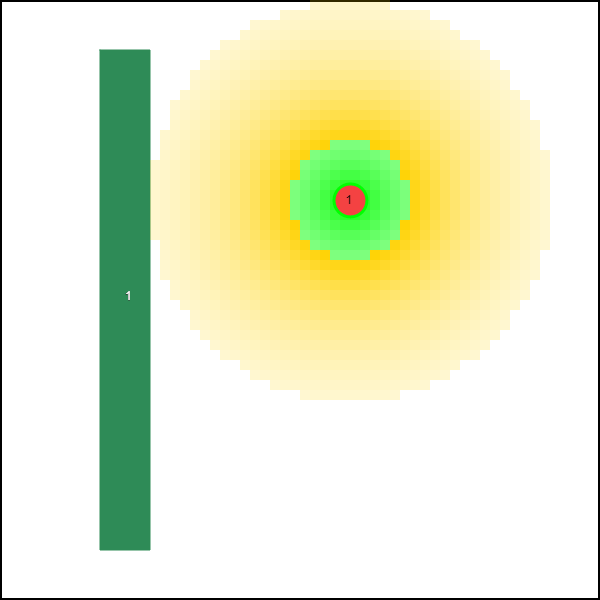}
\caption{$+ 121.4\%$\\VB:\ $30,259$}
\label{fig:NDV1}
\end{subfigure}
~ \begin{subfigure}{0.22\columnwidth}
\includegraphics[width=\textwidth]{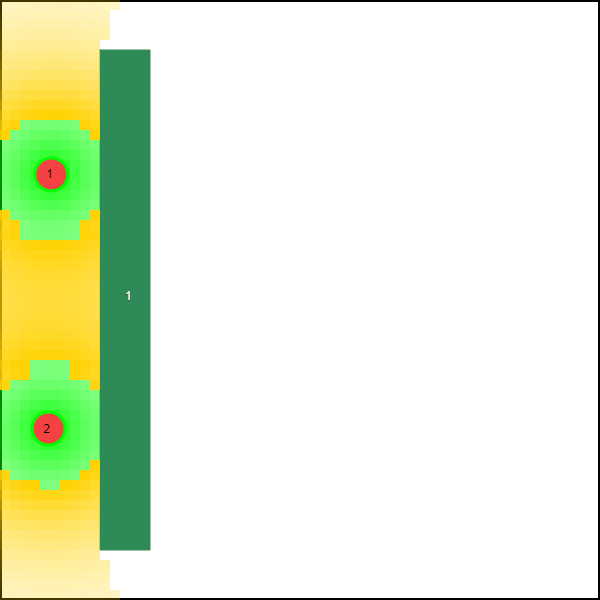}
\caption{\ \\GA:\ $26,258$}
\label{fig:NDZ2}
\end{subfigure}
\begin{subfigure}{0.22\columnwidth}
\includegraphics[width=\textwidth]{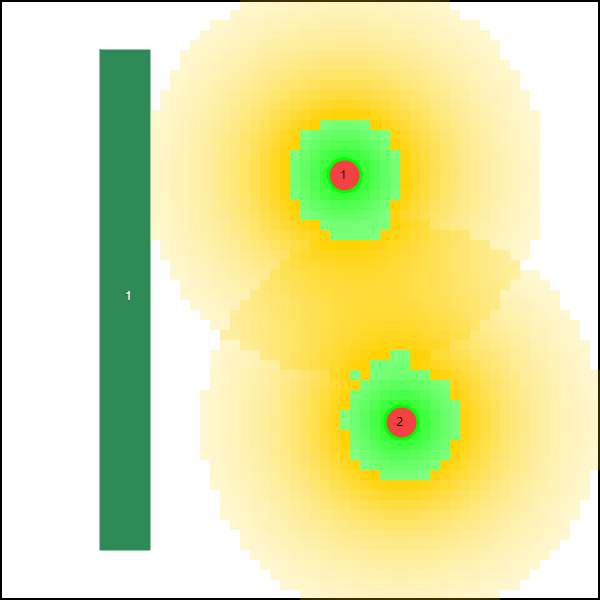}
\caption{$+ 123.5\%$\\VB:\ $58,693$}
\label{fig:NDV2}
\end{subfigure}
\caption{Maximum coverage improvement achieved due to boosting for $N=1,2$}%
\label{Fig:Results2}%
\end{figure}

\begin{table}[t]
\caption{Coverage objective value for cases with $N=5,6$ with decentralized
boosting (See Fig. \ref{Fig:Results3})}%
\label{Tab:Results3}
\centering
\begin{tabular}
[c]{|l|l|r|r|r|}\hline
\multicolumn{2}{|l|}{Configuration} &
\multicolumn{1}{c|}{\multirow{2}{*}{\begin{tabular}[c]{@{}c@{}}Gradient \\ Ascent (GA)\end{tabular}}}
&
\multicolumn{1}{c|}{\multirow{2}{*}{\begin{tabular}[c]{@{}c@{}}Decentralized \\ V-Boosting\end{tabular}}}
&
\multicolumn{1}{c|}{\multirow{2}{*}{\begin{tabular}[c]{@{}c@{}}Decentralized \\ Arc-Boosting\end{tabular}}}%
\\\cline{1-2}%
Obstacles & N & \multicolumn{1}{c|}{} & \multicolumn{1}{c|}{} &
\multicolumn{1}{c|}{}\\\hline
General & 5 & \textbf{93,637} & \textbf{97,214} & 96,832\\\hline
Maze & 6 & \textbf{90,953} & 94,026 & \textbf{94,436}\\\hline
Room & 5 & \textbf{86,638} & 89,078 & \textbf{89,088}\\\hline
Narrow & 6 & \textbf{101,976} & 116,481 & \textbf{129,476}\\\hline
\end{tabular}
\vspace{-5mm}
\end{table}

\begin{figure}[t]
\centering
\begin{subfigure}{0.22\columnwidth}
\includegraphics[width=\textwidth]{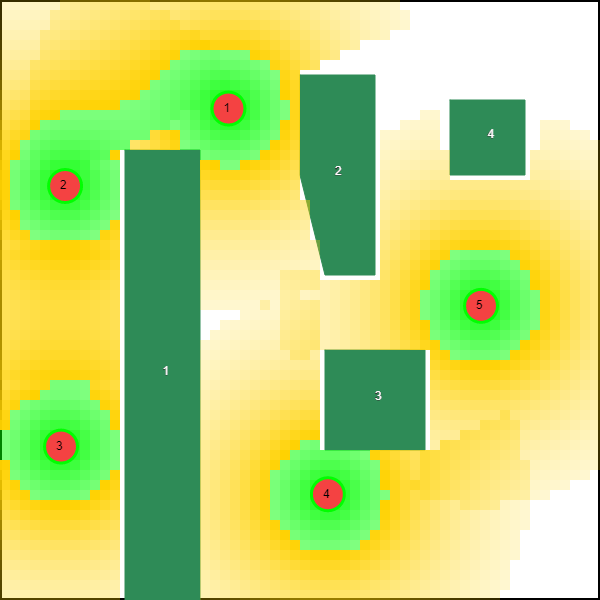}
\caption{\ \\GA:\ $93,637$}
\label{fig:GDZ5}
\end{subfigure}
\begin{subfigure}{0.22\columnwidth}
\includegraphics[width=\textwidth]{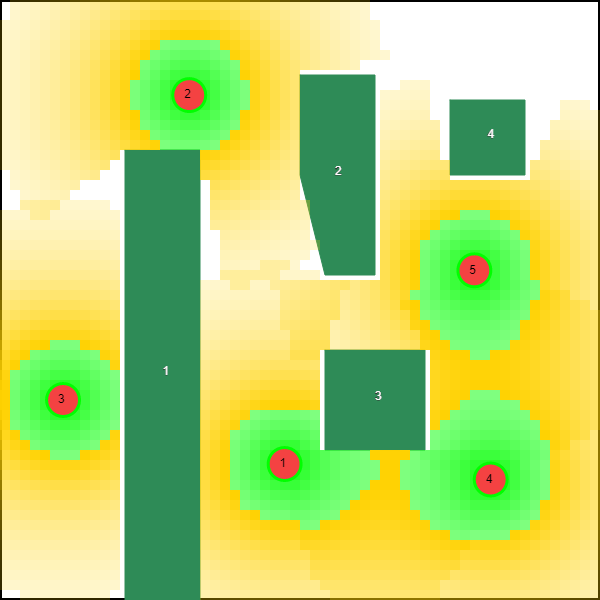}
\caption{$+ 3.820\%$\\VB:\ $97,214$}
\label{fig:GDV5}
\end{subfigure}
~ \begin{subfigure}{0.22\columnwidth}
\includegraphics[width=\textwidth]{figures/RDZ5.png}
\caption{\ \\GA:\ $86,638$}
\label{fig:RDZ5}
\end{subfigure}
\begin{subfigure}{0.22\columnwidth}
\includegraphics[width=\textwidth]{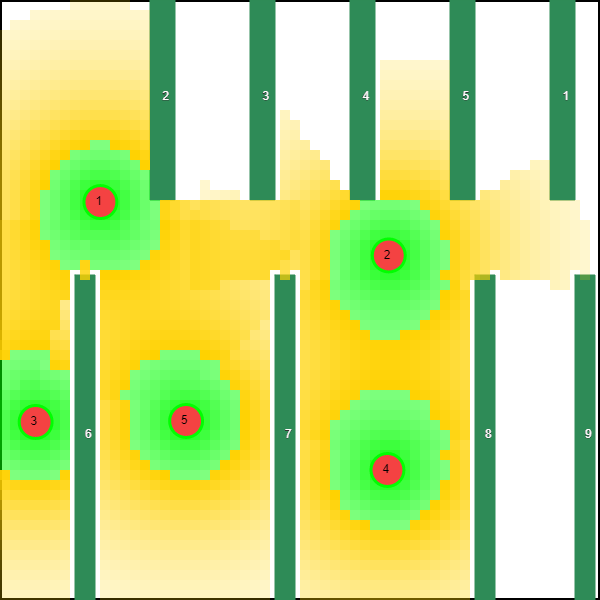}
\caption{$+ 2.828\%$\\AB:\ $89,088$}
\label{fig:RDA5}
\end{subfigure}
~ \begin{subfigure}{0.22\columnwidth}
\includegraphics[width=\textwidth]{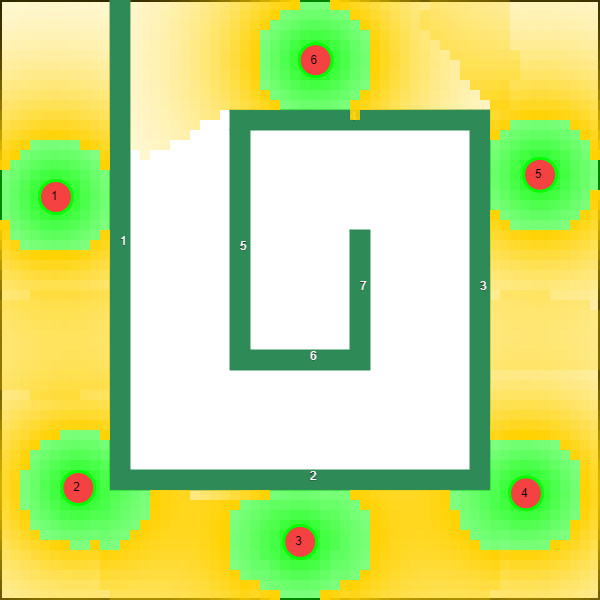}
\caption{\ \\GA:\ $90,953$}
\label{fig:MDZ6}
\end{subfigure}
\begin{subfigure}{0.22\columnwidth}
\includegraphics[width=\textwidth]{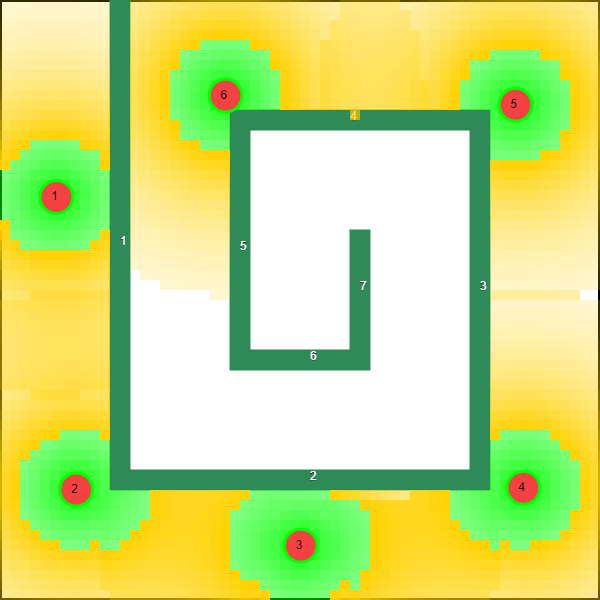}
\caption{$+ 3.829\%$\\AB:\ $94,436$}
\label{fig:MDA6}
\end{subfigure}
~ \begin{subfigure}{0.22\columnwidth}
\includegraphics[width=\textwidth]{figures/NDZ6.png}
\caption{\ \\GA:\ $101,976$}
\label{fig:NDZ6}
\end{subfigure}
\begin{subfigure}{0.22\columnwidth}
\includegraphics[width=\textwidth]{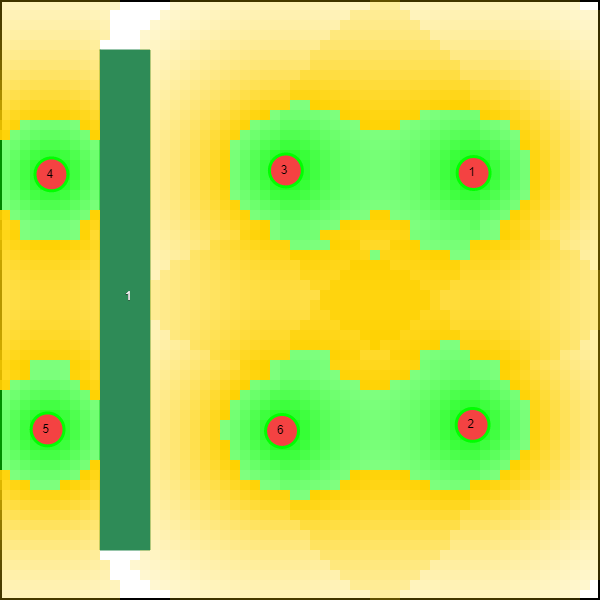}
\caption{$+ 26.97\%$\\AB:\ $129,476$}
\label{fig:NDA6}
\end{subfigure}
\caption{Maximum coverage improvement achieved due to boosting for $N=5,6$ (See Tab. \ref{Tab:Results3})}%
\label{Fig:Results3}
\vspace{-5mm}
\end{figure}

As the final step, obtained simulation results for the coverage control problem are presented which highlights the impact of the main contributions
of this work: 
(\romannum{1}) The generalized distributed multi-agent optimization problem solving technique based on boosting functions approach, 
(\romannum{2}) The convergence guaranteeing optimal step size selection method for the use of distributed boosting schemes, 
(\romannum{3}) The two new boosting function families (V-Boosting and Arc-Boosting) for the coverage control application, 
(\romannum{4}) The developed distributed boosting scheme for the coverage control application, and, 
(\romannum{5}) The application of convergence guaranteed optimal step sizes for the coverage control application. 

The proposed distributed coverage control algorithm (i.e., the DBS) including
the methods proposed in \cite{Zhong2011,Sun2014} were implemented in a
JavaScript-based simulator which is available at
\href{http://www.bu.edu/codes/simulations/shiran27/CoverageFinal/}{http://www.bu.edu/ codes/simulations/shiran27/CoverageFinal/}. The reader is invited to reproduce the reported results using the interactive interface and to explore the performance of the proposed method under diverse mission space environments and operating conditions. The source code of the simulator is also available at \href{https://github.com/shiran27/CoverageControl}{https://github.com/shiran27/CoverageControl}. The boosting function parameters used in generating the results reported next (i.e., gain parameters $\kappa,\gamma$) are listed in Table
\ref{Tab:BoostingParamters}.


\begin{remark}\label{Rm:BoostingFunctionScalingParameters}
The exact numerical values suitable for the gain parameters in different boosting function families (i.e., $\kappa,\gamma$) are application dependent. However, it is advisable to select those gain parameters such that the magnitudes of resulting boosted gradients (i.e. $\Vert \hat{d}_{i} \Vert$) and normal gradients (i.e., $\Vert d_{i} \Vert$) are in the same order (the initially computed values).
\end{remark}

In the simulations, four different mission space arrangements named
`General',`Room',`Maze' and `Narrow' are considered based on the obstacle
arrangement of each mission space. As the first step, the conventional
distributed gradient ascent method proposed in \cite{Zhong2011} was applied
for each of those mission spaces with 10 agents (i.e. $N=10$) to get the final
solutions shown in figures \ref{fig:GDZ10},\ref{fig:RDZ10},\ref{fig:MDZ10},
and, \ref{fig:NDZ10} respectively. The corresponding objective function values
are listed in Table \ref{Tab:Results1} under the column: `Reference Level
$H(\mathbf{s}^{1})$'. Also note that, as another baseline for the proposed boosting methods, a random gradient perturbation method is also implemented
which uses $\hat{d}_{i,k}=d_{i,k}+\kappa\zeta_{i,k}$ during the boosting
sessions. Here, $\kappa=5$ and $\zeta_{i,k}\in \R^{2}$ is a two-dimensional
random vector, independently generated from a standard uniform distribution at
each time step.

Then, as the next step, the effect of different boosting function families
proposed in Section  \ref{SubSubSec:DifferentBoostingFunctionFamilies} were
explored under the CBS \cite{Sun2014} and under the novel DBS. The increment achieved in the coverage objective value (with respect to the
reference level $H(\mathbf{s}^{1})$) by each of those boosting methods are
tabulated in Table \ref{Tab:Results1}. The cases with the highest coverage
objective value increments are shown in bold letters and they are illustrated
in the figures \ref{fig:GDA10}, \ref{fig:RDA10}, \ref{fig:MDA10} and
\ref{fig:NDV10}. Note that in all figures
(i.e., the subfigures under Fig. \ref{Fig:Results0}, Fig. \ref{Fig:Results1} and Fig.
\ref{Fig:Results2}) agent locations are shown in red-colored dots, and they have been initialized at the top left corner of the mission space. Further,
darker green colored areas indicate higher coverage levels.

The results in Table \ref{Tab:Results1} shows that the distributed Arc-Boosting (labeled "AB") and distributed V-Boosting (labeled "VB") schemes outperform all the other methods for all the tested obstacle configurations when $N=10$. Some results obtained with the distributed $\Phi$-Boosting (labeled "$\Phi$B") are also shown in Fig. \ref{Fig:Results0}. Moreover, simulation results obtained for cases with ~$N=1,2$ shown in Table \ref{Tab:Results2} and Fig. \ref{Fig:Results2} also leads to the same conclusion. Therefore, to further investigate the performance of the distributed V-Boosting and Arc-Boosting methods, simulation results were generated with moderate $N$ values, where $N=5,6$. The corresponding results are shown in Table \ref{Tab:Results3} and Fig. \ref{Fig:Results3}.

From the presented results, it is clear that boosting functions approach can
successfully evade the local optima given by the conventional
gradient ascent based method. Further, the systematic gradient modification process
achieved via boosted gradients and the distributed boosting scheme enables such approaches to obtain superior objective function values compared to conventional gradient ascent based methods as well as compared to random gradient perturbation based techniques.

\begin{table}[t]
\caption{Comparison of coverage objective and convergence time
values observed for the DBS with fixed and variable steps.}%
\label{Tab:FixedVsVariableForDBS}%
\centering
\resizebox{\columnwidth}{!}{\begin{tabular}{|c|l|r|r|r|r|}
\hline
\multirow{2}{*}{\begin{tabular}[c]{@{}c@{}}Boosting \\ Method\end{tabular}} & \multicolumn{1}{c|}{($N=8$)} & \multicolumn{2}{c|}{$H(s^*)$} & \multicolumn{2}{c|}{Convergence Time} \\ \cline{2-6}
& \multicolumn{1}{c|}{Configuration} & \multicolumn{1}{c|}{\begin{tabular}[c]{@{}c@{}}Fixed\\ steps\end{tabular}} & \multicolumn{1}{c|}{\begin{tabular}[c]{@{}c@{}}Variable\\ steps\end{tabular}} & \multicolumn{1}{c|}{\begin{tabular}[c]{@{}c@{}}Fixed\\ steps\end{tabular}} & \multicolumn{1}{c|}{\begin{tabular}[c]{@{}c@{}}Variable\\ steps\end{tabular}} \\ \hline
\multirow{4}{*}{\rotatebox[origin=c]{90}{\begin{tabular}[c]{@{}c@{}} V- \\ Boosting\end{tabular}}} & General & 140,592 & 140,649 & 550.7 & 91.3 \\ \cline{2-6}
& Room & 127,557 & 127,517 & 613.5 & 140.3 \\ \cline{2-6}
& Maze & 120,832 & 121,231 & 302.2 & 134.1 \\ \cline{2-6}
& Narrow & 163,478 & 155,528 & 415.7 & 161.8 \\ \hline
\multirow{4}{*}{\rotatebox[origin=c]{90}{\begin{tabular}[c]{@{}c@{}}Arc-\\ Boosting\end{tabular}}} & General & 140,615 & 140,542 & 80.3 & 104.9 \\ \cline{2-6}
& Room & 127,647 & 127,455 & 390.0 & 158.1 \\ \cline{2-6}
& Maze & 119,967 & 121,231 & 151.7 & 125.0 \\ \cline{2-6}
& Narrow & 155,641 & 155,485 & 127.3 & 88.1 \\ \hline
\multicolumn{2}{|c|}{Average:} & 137,041 & 136,205 & 328.9 & 125.4 \\ \hline
\end{tabular}}
\vspace{-5mm}
\end{table}

\subsubsection{Discussion}

\paragraph*{\textbf{Different boosting function families}}
Different boosting function families can be useful in different mission space
configurations. For example, when there are one or no obstacles in the mission
space, the V-Boosting method performed the best. However, when there are
multiple obstacles in the environment, the Arc-Boosting method gave the best result.

\paragraph*{\textbf{Effect of decentralization}}

Due to decentralization, overall the simulations carried out for $N=10$, on average (per simulation) the convergence time to the final optimal solution is improved (i.e., reduced) by $39.97\%$ (i.e., by approximately $165.2\ s$), on an Intel\textsuperscript{\tiny \textregistered} Core\texttrademark\ i7-8700
CPU @ 3.20 GHz Processor with a 32 GB RAM. Further, on average (per
simulation), the final coverage cost achieved is increased by $0.381\%$ (approximately $451$ units) due to decentralization. Furthermore, inherently, decentralization reduces the associated communication cost compared to a centralized implementation. Therefore, the proposed DBS clearly outperforms the
CBS in every aspect.

\paragraph*{\textbf{Variable step size scheme and Convergence}}

In the simulations, whenever the decentralized boosting scheme (Proposed in
Section \ref{SubSec:CoverageControlBoostingSchemes}) is used, the variable
step size method proposed in Section \ref{Sec:VariableStepSizeScheme} (in
Theorem \ref{Th:ConvergenceHat}) was used to guarantee the convergence. When
the methods proposed in \cite{Zhong2011,Sun2014} are simulated, to preserve
the typical operating conditions, the modified conventional step size
selection method described in subsection
\ref{SubSubSec:StepSizeConstantDerivation} was used. However, under each method, the exact same terminal condition was used to terminate the simulation (i.e., to determine the final convergence). Fig. \ref{Fig:StepSizeVariationExample} shows an example step size sequence and the associated gradient sequence with the agent $i=4$ during the simulation, which leads to the result shown in \ref{fig:NDV10}. Moreover, Table
\ref{Tab:FixedVsVariableForDBS} provides a comparison of coverage objective
and convergence time values observed for the DBS when fixed and variable step
sizes are used. Note that the use of variable step sizes has improved (i.e.,
reduced) the convergence time by $61.9\%$ (i.e., by $203.5 s$).

\section{conclusion}

\label{Sec:Colclusion}

The concept of boosting provides systematic ways to overcome the problem of multiple local optima arising in cooperative multi-agent optimization problems with non-convex objective functions. An optimal step size selection scheme is developed to guarantee convergence in a distributed (or centralized) boosting framework for such general multi-agent optimization problems. The application of boosting functions is illustrated using the class of cooperative multi-agent coverage control problems, where two novel boosting function families are developed and applied successfully. A new distributed framework is also proposed for the use of boosted gradients for optimal coverage control problems. Simulation results are used to illustrate the effectiveness of the proposed boosting functions and the distributed boosting framework. Ongoing research aims to explore the generality of boosting functions to be used regardless of the intended application. Also, concerning coverage control applications, current research aims to explore the possibility of combining multiple methods used for computing the boosted gradients, to reap the aggregate benefits of all such methods.













\bibliographystyle{IEEEtran}
\bibliography{References}

\end{document}